\documentclass[12pt,a4paper,english]{article}
\usepackage{amsmath,amssymb}
\usepackage{babel}
\usepackage{eucal}

\DeclareMathVersion{scriptstyle}

\usepackage[latin1]{inputenc}\usepackage[T1]{fontenc}
\usepackage{theorem}
\theoremstyle{change}
{\theorembodyfont{\slshape} 
  \newtheorem{thm}{Theorem.}[section]
  \newtheorem{lemma}[thm]{Lemma.}
  \newtheorem{prop}[thm]{Proposition.}
  \newtheorem{cor}[thm]{Corollary.}
}
{\theorembodyfont{\rmfamily} 

  \newtheorem{remark}[thm]{Remark.}

}

\numberwithin{equation}{section} 
\def\labelenumi{(\roman{enumi})} 

\newcommand{\supp}{\mathrm{supp}}
\newcommand{\dist}{\mathrm{dist}}
\newcommand{\Prob}{\mathrm{Prob}}
\newcommand{\cA}{{\cal A}}

\newcommand{\cJ}{{\cal J}}
\newcommand{\sa}{\mathrm{sa}}
\newcommand{\GUE}{\mathrm{GUE}}
\newcommand{\bR}{{\mathbb R}}
\newcommand{\bZ}{{\mathbb Z}}
\newcommand{\bH}{{\mathbb H}}
\newcommand{\bV}{{\mathbb V}}
\newcommand{\bE}{{\mathbb E}}
\newcommand{\bC}{{\mathbb C}}
\newcommand{\bN}{{\mathbb N}}
\newcommand{\Ext}{\mathrm{Ext}}
\newcommand{\red}{\mathrm{red}}
\newcommand{\proof}[1][Proof. ]{{\it#1}}

\newcommand{\eproof}{\hfill $\Box$\\}

\def\endproof{{\nobreak\qquad{\scriptstyle \blacksquare}}}

\def\C{{\mathbb C}}

\def\E{{\mathbb E}}
\def\F{{\mathbb F}}

\def\N{{\mathbb N}}
\def\R{{\mathbb R}}

\def\Z{{\mathbb Z}}
\def\V{{\mathbb V}}
\def\Q{{\mathbb Q}}

\def\CA{{\mathcal A}}
\def\CB{{\mathcal B}}

\def\CD{{\mathcal D}}
\def\CE{{\mathcal E}}
\def\CF{{\mathcal F}}

\def\CI{{\mathcal I}}
\def\CJ{{\mathcal J}}

\def\CO{{\mathcal O}}
\def\CP{{\mathcal P}}

\def\CR{{\mathcal R}}

\def\CU{{\mathcal U}}

\def\unit{{\bf 1}}
\def\i{{\rm i}}
\def\e{{\rm e}}
\def\d{{\rm d}}
\def\eps{\varepsilon}
\def\<{{\langle}}
\def\>{{\rangle}}

\def\id{{\text{\rm id}}}
\def\sumr{\sum_{i=1}^r}
\def\sumrs{\sum_{i=1}^{r+s}}
\def\tr{{\text{\rm tr}}}
\def\Tr{{\text{\rm Tr}}}
\def\GRM{{\text{\rm GRM}}}
\def\SGRM{{\text{\rm SGRM}}}
\def\GRMR{{\text{\rm GRM}}^{\R}}
\def\GOE{{\text{\rm GOE}}}
\def\GOES{{\text{\rm GOE}}^\ast}
\def\GSE{{\text{\rm GSE}}}
\def\GSES{{\text{\rm GSE}}^\ast}
\def\im{{\text{\rm Im}}}
\def\re{{\text{\rm Re}}}
\def\supp{{\text{\rm supp}}}
\def\grad{{\text{\rm grad}}}
\def\Mnsa{M_n(\C)_{sa}}
\def\diff{\frac{\rm d}{{\rm d}t}\Big|_{t=0}}

\def\cc{^*}

\def\Distc{\CD_c'(\R)}

\def\Ccinf{C_c^\infty(\R)}

\def\unitH{\begin{pmatrix}
1 & 0\\
0 & 1
\end{pmatrix}}
\def\j{\begin{pmatrix}
\i & 0\\
0 & -\i
\end{pmatrix}}
\def\k{\begin{pmatrix}
0 & 1\\
-1 & 0
\end{pmatrix}}
\def\l{\begin{pmatrix}
0 & \i\\
\i & 0
\end{pmatrix}}

\begin{document}

\newcommand{\firstnote}{%
Affiliated with MaPhySto - A network in Mathematical Physics and Stochastics, which is funded by a great from the Danish National Research Foundation.
}

\newcommand{\secondnote}{%
Partially supported by MaPhySto and by the PhD-school OP-ALG-TOP-GEO, which is funded by the Danish Training Research Counsil.
}

\title{\bf A Random Matrix Approach to\\ the Lack of Projections in $C^*_\red(\F_2)$}

\author{Uffe Haagerup\footnote{\firstnote},\; Hanne Schultz\footnote{\secondnote}\; and Steen Thorbjørnsen\footnotemark[1]}

\date{}
 
\maketitle

\begin{abstract}
\noindent In 1982 Pimsner and Voiculescu computed the $K_0$- and
$K_1$-groups of the reduced group $C^*$-algebra $C^*_\red(F_k)$ of the free
group $F_k$ on $k$ generators and settled thereby a long standing
conjecture: $C^*_\red(F_k)$ has no projections except for the trivial projections 0 and 1. Later
simpler proofs of this conjecture were found by methods from K-theory or
from non-commutative differential geometry. In this paper we provide a new proof of the fact that
$C^*_\red(F_k)$ is projectionless. The new proof is based on random
matrices and is obtained by a refinement of the methods recently used by
the first and the third named author to show that the semigroup
$\Ext(C^*_\red(F_k))$ is not a group for $k\ge 2$. By the same type of
methods we also obtain that two phenomena proved by Bai and Silverstein for
certain classes of random matrices: ``no eigenvalues outside (a small
neighbourhood of) the support of the limiting distribution'' and ``exact
separation of eigenvalues by gaps in the limiting distribution'' also hold
for arbitrary non-commutative selfadjoint polynomials of independent GUE,
GOE or GSE random matrices with matrix coefficients.
\end{abstract}

\section{Introduction.}
\label{sec1}
In \cite{HT} the first and the third named author proved the following
extension of Voiculescu's random matrix model for a semicircular system:

Let $X_1^{(n)},\dots,X_r^{(n)}$ be $r$ independent selfadjoint $n\times n$
random matrices from Gaussian unitary ensembles (GUE) and with the scaling
used in Voiculescu's paper \cite{V1}. Moreover, let $x_1,\dots,x_r$ be a
semicircular system in a $C^*$-probability space $(A,\tau)$ with $\tau$
faithful. Then, for every polynomial $p$ in $r$ non-commuting variables,
\begin{equation}
\label{eq1-1}
\lim_{n\to\infty} \| p(X_1^{(n)},\dots,X_r^{(n)})\| = \| p(x_1,\dots,x_r)\|
\end{equation}
holds almost surely.\newpage The main steps in the proof of \eqref{eq1-1} were:

{\bf STEP 1 (Linearization trick).} In order to prove \eqref{eq1-1}, it is sufficient to show that for every
$m\in\N$, every $\varepsilon >0$ and every selfadjoint polynomial $q$ in
$r$ non-commuting variables with coefficients in $M_m(\C)$ and with
deg$(q)=1$, 
\begin{equation}
\label{eq1-2}
\sigma(q(X_1^{(n)},\dots,X_r^{(n)}))\subseteq
\sigma(q(x_1,\dots,x_r))+(-\varepsilon,\varepsilon) 
\end{equation}
eventually as $n\to\infty$ (almost surely). Here $\sigma(\cdot)$ denotes
the spectrum of a matrix or of an element in a $C^*$-algebra.

{\bf STEP 2 (Mean value estimate).} If $q$ is a polynomial of first degree with matrix coefficients as in step 1,
then for every $\varphi\in C_c^\infty(\bR,\bR)$ 
\begin{eqnarray}
\label{eq1-3}
\bE\{(\tr_m\otimes\tr_n)\varphi(q(X_1^{(n)},\dots,X_r^{(n)}))\} =
 (\tr_m\otimes\tau)\varphi(q(x_1,\dots,x_r))+O(\tfrac{1}{n^2}) 
\end{eqnarray}
where $\tr_m=\frac1m \Tr_m$ is the normalized trace on $M_m(\bC)$.

{\bf STEP 3 (Variance estimates).} If $q$ is a polynomial of first degree with matrix coefficients as in step 1,
then for every $\varphi\in C_c^\infty (\bR,\bR)$,
\begin{equation}
\label{eq1-4}
\bV\{(\tr_m\otimes
\tr_n)\varphi(q(X_1^{(n)},\dots,X_r^{(n)}))\}=O(\tfrac{1}{n^2}).
\end{equation}
Moreover, if $\varphi'=\frac{d\varphi}{dx}$ vanishes in a neighbourhood of
$\sigma(q(x_1,\dots,x_r))$, then
\begin{equation}
\label{eq1-5}
\bV\{(\tr_m\otimes\tr_n)\varphi(q(X_1^{(n)},\dots,X_r^{(n)}))\}=O(\tfrac{1}{n^4}).
\end{equation}
A standard application of the Borel-Cantelli lemma and the Chebychev
inequality to \eqref{eq1-3} and \eqref{eq1-5} gives that if $\varphi'$
vanishes on a neighbourhood of $\sigma(q(x_1,\dots,x_r))$, then
\begin{eqnarray}
\label{eq1-6}
(\tr_m\otimes\tr_n)\varphi(q(X_1^{(n)},\dots,X_r^{(n)})) = (\tr_m\otimes\tau)\varphi(q(x_1,\dots,x_r))+O(n^{-\frac43})
\end{eqnarray}
holds almost surely, and from this \eqref{eq1-2} easily follows
(cf. \cite[proof of Theorem~6.4]{HT}).

In \cite{S}, the second named author generalized the above to real and
symplectic Gaussian random matrices (the GOE- and GSE-cases). The main new
problem in these two cases is that \eqref{eq1-3} no longer holds. However,
the following formula holds (cf. \cite[Theorem~5.6]{S}):
\begin{eqnarray}
\label{eq1-7} \bE\{(\tr_m\otimes\tr_n)\varphi(q(X_1^{(n)},\dots,X_r^{(n)}))\}=
 (\tr_m\otimes\tau)\varphi(q(x_1,\dots,x_r))+\textstyle{\frac1n}
 \Lambda(\varphi)+O(\tfrac{1}{n^2}),
\end{eqnarray}
where $\Lambda\colon C_c^\infty(\bR)\to\bC$ is a distribution (in the sense
of L.~Schwartz) depending on the polynomial $q$ and on the scalar field
($\bR$ or $\bH$). Moreover, by \cite[Lemma~5.5]{S},
\begin{equation}
\label{eq1-8}
\supp(\Lambda)\subseteq\sigma(q(x_1,\dots,x_r)).
\end{equation}
Using \eqref{eq1-7} and \eqref{eq1-8} instead of \eqref{eq1-3} the proofs
of \eqref{eq1-2} and \eqref{eq1-1} could be completed essentially as in the
GUE-case.

The proof of the linearization trick relied on $C^*$-algebra techniques,
namely on Stinespring's theorem and on Arveson's extension theorem for
completely 
positive maps. In the present paper we give a purely algebraic proof of the
linearization trick which in turn allows us to work directly with
polynomials of degree greater than 1. As a result, we prove in
Section~\ref{sec4b} and Section~\ref{sec8} that \eqref{eq1-3} for the
GUE-case (resp. \eqref{eq1-7} and 
\eqref{eq1-8} for the GOE- and GSE-cases) holds for selfadjoint
polynomials $q$ of {\it any} degree with coefficients in $M_m(\bC)$. Also
\eqref{eq1-4}, \eqref{eq1-5} and \eqref{eq1-6} hold in this
generality. Consequently, \eqref{eq1-2} holds for all such polynomials $q$
in all three cases (GUE, GOE and GSE). This is the phenomenon ``no
eigenvalues outside (a small neighbourhood of) the support of the limiting
distribution'', which Bai and Silverstein obtained in \cite {BS1} for a
different class of selfadjoint random matrices.

Let us next discuss the application to projections in $C^*_\red(\F_r)$:
Recall from \cite[Lemma~8.1]{HT} that $C^*_\red(\F_r)$ has a unital, trace
preserving embedding in $C^*(x_1,\dots,x_r,\unit)$, where $x_1,\dots,x_r$ is a
semicircular system. If $e$ is a projection in $M_m(C^*(x_1,\dots,x_r,\unit))$,
then by standard $C^*$-algebra techniques (cf. Section~\ref{sec5}) there
exists a 
projection $f$ in $M_m(C^*(x_1,\dots,x_r,\unit))$ such that $\|e-f\|<1$ and such
that $f$ takes the form
\[
f=\varphi(q(x_1,\dots,x_r)),
\]
where $q$ is a selfadjoint polynomial in $r$ non-commuting variables with
coefficients in $M_m(\C)$, and
$\varphi$ is a $C^\infty$-function with compact support, such that
$\varphi$ only takes the values 0 and 1 is some neighbourhood of
$\sigma(q(x_1,\dots,x_r))$.

Consider now random matrices $X_1^{(n)},\dots,X_r^{(n)}$ as in the GUE-case
described above. By \eqref{eq1-6} we have that 
\begin{eqnarray*}
(\tr_m\otimes\tr_n)\varphi(q(X_1^{(n)},\dots,X_r^{(n)})) &=&
(\tr_m\otimes\tau)\varphi(q(x_1,\dots,x_r))+O(n^{-\frac43})\\
&=& (\tr_m\otimes\tau)(f) + O(n^{-\frac43})
\end{eqnarray*}
holds almost surely and hence the corresponding unnormalized trace
satisfies
\begin{equation}
\label{eq1-9}
(\Tr_m\otimes\Tr_n)\varphi(q(X_1^{(n)},\dots,X_r^{(n)})) =
n(\Tr_m\otimes\tau)(f)+O(n^{-\frac13}).
\end{equation}
Using that the left hand side of \eqref{eq1-9} is an
integer for all large $n\in\bN$, it is not hard to prove that
$(\Tr_m\otimes\tau)(f)$ is an integer (cf. section~7 for details). Moreover, since $\|e-f\|<1$ implies
that $e=ufu^*$ for a unitary $u\in M_m(C^*(x_1,\dots,x_r,1))$, we also have
\begin{equation}
\label{eq1-10}
(\Tr_m\otimes\tau)(e)\in\bZ.
\end{equation}
Hence, using the existence of a unital trace-preserving embedding of
$C^*_\red(\F_r)$ into\\ $C^*(x_1,\dots,x_r,\unit)$, it follows that:
\begin{equation}
\label{eq1-11}
(\Tr_m\otimes\tau)(e)\in\bZ\quad\mbox{for all projections} \quad e\in
M_m(C^*_\red(\F_r)).
\end{equation}
In particular:
\begin{equation}
\label{eq1-12}
C^*_\red(\F_r)\quad\mbox{has no projections except 0 and 1}.
\end{equation}

The two statements \eqref{eq1-11} and \eqref{eq1-12} were first obtained by
Pimsner and Voiculescu in 1982 (cf. \cite{PV}) by proving that
$K_0(C^*_\red(\F_r))=\bZ$, where the $K_0$-class $[\unit]$ of the unit in
$C^*_\red(\F_r)$ corresponds to $1\in\bZ$. Simpler proofs of
$K_0(C^*_\red(\F_r))=\bZ$ were later obtained by Cuntz \cite{Cu1} and Lance
\cite{L}. Connes gave in \cite[pp.~269--272]{Co} a more direct proof of
\eqref{eq1-12} based on Fredholm modules. Connes' argument can be further
simplified to a short selfcontained proof without explicit mentioning of
Fredholm modules (\cite{Cu2}, \cite{CF}).

It is an elementary consequence of \eqref{eq1-11} that for every selfadjoint
polynomial $q$ in $r$ non-commuting variables and with coefficients in
$M_m(\bC)$, the spectrum of $q(x_1,\dots,x_r)$ has at most $m$ connected
components, that is 
\[
\sigma(q(x_1,\dots,x_r)) = I_1\cup\dots\cup I_j,\qquad ({\rm disjoint}\;{\rm union}),
\]
where each $I_i$ is a compact interval or a one-point set
(cf. Proposition~\ref{sec6prop1}), and $j\leq m$. Let
\[
\unit =e_1+\dots + e_j
\]
be the corresponding decomposition of the unit in
$M_m(C^*(x_1,\dots,x_r,\unit ))$ into ortogonal projections, and put
\[
k_i = (\Tr_m\otimes\tau)(e_i)\in\bZ.
\]
Let now $0<\varepsilon<\frac13\varepsilon_0$ where $\varepsilon_0$ is the
smallest length of the gaps between the sets $I_1,\dots,I_j$. In
Section~\ref{sec6} and Section~\ref{sec9} we prove that the number of eigenvalues of
$q(X_1^{(n)},\dots,X_r^{(n)})$ in the open intervals
$I_i+(-\varepsilon,\varepsilon)$, $i=1,\dots,j$, are exactly $nk_i$
eventually as $n\to\infty$ (almost surely) in all three cases (GUE, GOE and
GSE). This is the phenomenon ``Exact separation of eigenvalues'' by the
gaps in the support of the limiting distribution which Bai and Silverstein
obtained in \cite{BS2} for the class of selfadjoint random matrices previously
studied in \cite{BS1}.

\section{Matrix Results.}\label{sec2}

\begin{prop}\label{faktorisering}
Let $d,m$ and $m'$ be positive integers and let $p$ be a polynomial in
$M_{m,m'}\otimes\C\<X_1,\ldots,X_r\>$ of degree $d$. Then there exist
positive integers $m_1,m_2,\ldots,m_{d+1}$ and polynomials
\[
u_j\in M_{m_j,m_{j+1}}(\C)\otimes\C\<X_1,\ldots,X_r\>, \qquad
(j=1,2,\ldots,d),
\]
such that

\begin{itemize}

\item[(i)] $m_1=m$ and $m_{d+1}=m'$,

\item[(ii)] $\deg(u_j)\le1$ for all $j$ in $\{1,2,\ldots,r\}$,

\item[(iii)] $p=u_1u_2\cdots u_d$.

\end{itemize}
\end{prop}

\proof The proof proceeds by induction on $d$. Noting that the case
$d=1$ is trivial, we assume that $d\ge2$ and that the proposition has been
verified for all polynomials in $M_{m,m'}(\C)\otimes\C\<X_1,\ldots,X_r\>$
og degree at most $d-1$. Given then a polynomial $p$ in
$M_{m,m'}(\C)\otimes\C\<X_1,\ldots,X_r\>$ of degree $d$, we may, setting
$X_0=\unit_{\C\<X_1,\ldots,X_r\>}$, write $p$ in the form
\[
p(X_1,\ldots,X_r)=\sum_{0\le i_1,i_2,\ldots,i_d\le r}c(i_1,i_2,\ldots,i_d)\otimes
X_{i_1}X_{i_2}\cdots X_{i_d},
\]
for suitable matrices
\[
c(i_1,i_2,\ldots,i_d)\in M_{m,m'}(\C), \qquad (i_1,i_2,\ldots,i_d\in\{0,1,\ldots,r\}).
\]
For any $i_1$ in $\{0,1,\ldots,r\}$, we put
\[
Y(i_1)=\sum_{0\le i_2,\ldots,i_d\le r}c(i_1,i_2,\ldots,i_d)\otimes
X_{i_2}\cdots X_{i_d}.
\]
Note then that
\[
\begin{split}
p(X_1,\ldots,X_r)
&=\sum_{i_1=0}^r(\unit_m\otimes X_{i_1})\Big(\sum_{0\le i_2,\ldots,i_d\le
  r}c(i_1,i_2,\ldots,i_d)\otimes
X_{i_2}\cdots X_{i_d}\Big) \\[.2cm]
&=\sum_{i_1=0}^r(\unit_{m}\otimes X_{i_1})Y(i_1) \\[.2cm]
&=
\begin{pmatrix}
\unit_m\otimes X_0 & \unit_m\otimes X_1 & \cdots & \unit_m\otimes X_r
\end{pmatrix}\cdot
\begin{pmatrix}
Y(0) \\
Y(1) \\
\vdots \\
Y(r)
\end{pmatrix} \\[.2cm]
&=u_1(X_1,\ldots,X_r)\cdot p'(X_1,\ldots,X_r),
\end{split}
\]
where
\[
u_1(X_1,\ldots,X_r):=
\begin{pmatrix}
\unit_m\otimes X_0 & \unit_m\otimes X_1 & \cdots & \unit_m\otimes X_r
\end{pmatrix}
\in M_{m,(r+1)m}(\C)\otimes\C\<X_1,\ldots,X_r\>,
\]
and
\[
p'(X_1,\ldots,X_r):=
\begin{pmatrix}
Y(0) \\
Y(1) \\
\vdots \\
Y(r)
\end{pmatrix}
\in M_{(r+1)m,m'}(\C)\otimes\C\<X_1,\ldots,X_r\>.
\]
We note that $\deg(u_1)=1$ and that $\deg(p')\le d-1$. By the
induction hypothesis, there are positive integers
$m_2,m_3,\ldots,m_{d+1}$ and polynomials $u_j$ in
$M_{m_j,m_{j+1}}(\C)\otimes\C\<X_1,\ldots,X_r\>$ ($j=2,3,\ldots,d$)
such that

\begin{itemize}

\item[(i$'$)] $m_2=(r+1)m$ and $m_{d+1}=m'$,

\item[(ii$'$)] $\deg(u_j)\le 1$,  for all $j$ in $\{2,3,\ldots,d\}$,

\item[(iii$'$)] $p'=u_2u_3\ldots u_d$.

\end{itemize}

Now, $p=u_1p'=u_1u_2\cdots u_d$ and we have the desired decomposition.
$\endproof$

\begin{remark}\label{eksplicite udtryk}
By inspection of the proof of Proposition~\ref{faktorisering}, it is
apparent that the decomposition, implicitly given in that proposition, is
explicitly given as follows:
\[
\begin{split}
m_1&=m, \ m_2=(r+1)m, \ m_3=(r+1)m_2=(r+1)^2m, \ \ldots \ ,
m_d=(r+1)^{d-1}m, \\[.2cm]
m_{d+1}&=m' \\[.2cm]
u_j&=
\begin{pmatrix}
\unit_{m_j}\otimes X_0 & \unit_{m_j}\otimes X_1 & \cdots & \unit_{m_j}\otimes X_r
\end{pmatrix},
\qquad (j=1,2,\ldots,d-1), \\[.2cm]
u_d&=
\begin{pmatrix}
\sum_{i_d=0}^rc(i_1,\ldots,i_{d-1},i_d)\otimes X_{i_d}
\end{pmatrix}_{0\le i_1,i_2,\ldots,i_{d-1}\le r},
\end{split}
\]
where $u_d$ should be thought of as a block column matrix with
(block) rows indexed by the tuples
$(i_1,i_2,\ldots,i_{d-1})$ in a certain order.

Note in particular that the polynomials $u_1,u_2,\ldots,u_{d-1}$ are
basicly canonical, in the sense that they only depend on $p$ through
the degree $d$ and the dimension $m$. Conversely, the
polynomial $u_d$ basicly contains all information about $p$.
\end{remark}

\begin{prop}\label{inv-formel}
Let $\CA$ be an algebra with unit $\unit_{\CA}$ and let
$d,m,m_1,m_2,\ldots,m_{d+1}$ be positive integers such that
$m_1=m=m_{d+1}$. Put $k=\sum_{j=1}^dm_j$.

Consider further for each $j$ in $\{1,2,\ldots,d\}$ a
matrix $u_j$ from $M_{m_j,m_{j+1}}(\CA)$, and note that $u_1u_2\cdots
u_d\in M_m(\CA)$. For each $\lambda$ in $M_m(\CA)$, define the matrix
$A(\lambda)$ in $M_{k}(\CA)$ by
\begin{equation}
A(\lambda)=
\begin{pmatrix}
\lambda & -u_1 & 0 & 0 & \cdots & 0    \\
0 & \unit_{m_2} & -u_2 & 0 & \cdots & 0 \\
0  &  0 & \unit_{m_3} & -u_3 & \cdots & 0 \\
\vdots & \vdots & \vdots & \ddots & \ddots & \vdots \\
0 & 0 & 0 & \cdots & \unit_{m_{d-1}} & -u_{d-1} \\
-u_d & 0 & 0 & \cdots & 0 & \unit_{m_d} 
\end{pmatrix},
\label{e1.1}
\end{equation}
where $\unit_{m_j}$ denotes the unit in $M_{m_j}(\CA)$.
For any $\lambda$ in $M_m(\C)$ we then have
\[
\lambda-u_1u_2\cdots u_d \quad \mbox{is invertible in $M_m(\CA)$}
\iff A(\lambda) \quad \mbox{is invertible in $M_k(\CA)$},
\]
in which case
\[
A(\lambda)^{-1}=B(\lambda)+C,
\]
where
\[
B(\lambda)=
\begin{pmatrix}
\unit_m \\
u_2u_3\cdots u_d \\
u_3u_4\cdots u_d \\
u_4\cdots u_d \\
\vdots \\
u_d
\end{pmatrix}
\big(\lambda-u_1u_2\cdots u_d\big)^{-1}
\begin{pmatrix}
\unit_m & u_1 & u_1u_2 & u_1u_2u_3 & \cdots & u_1u_2\cdots u_{d-1}
\end{pmatrix}
\]
and
\[
C=
\begin{pmatrix}
0 & 0 & 0 & 0 & 0 & \cdots & 0 \\
0 & \unit_{m_2} & u_2 & u_2u_3 & u_2u_3u_4 & \cdots & u_2u_3\cdots
u_{d-1} \\
0& 0 & \unit_{m_3} & u_3 & u_3u_4 & \cdots & u_3u_4\cdots u_{d-1} \\
0 & 0 & 0 & \unit_{m_4} & u_4 & \cdots & u_4u_5\cdots u_{d-1} \\
\vdots & \vdots & \vdots & \vdots & \ddots & \ddots & \vdots   \\
0 & 0 & 0 & 0 & \cdots & \unit_{m_{d-1}} & u_{d-1} \\
0 & 0 & 0 & 0 & \cdots & 0 & \unit_{m_d} \\
\end{pmatrix}.
\]
Note, in particular, that $(\lambda-u_1u_2\cdots u_d)^{-1}$ is the
(block-) entry at position $(1,1)$ of $A(\lambda)^{-1}$.
\end{prop}

\proof At first assume that $A(\lambda)$ is invertible with inverse
$F(\lambda)$. We write $F(\lambda)$ in block matrix form as
\[
F(\lambda)=\big(f_{i,j}(\lambda)\big)_{1\le i,j\le d},
\]
corresponding to the block matrix form of $A(\lambda)$:
\[
A(\lambda)=\big(a_{i,j}(\lambda)\big)_{1\le i,j\le d},
\]
specified above.

>From the equality $\unit_k=A(\lambda)F(\lambda)$, we get, in particular, the
identities
\begin{equation}
\sum_{i=1}^da_{i,j}(\lambda)f_{i,1}(\lambda)={\bf \delta}_{j,1}, \qquad (j=1,2,\ldots,n),
\label{e1.2}
\end{equation}
where
\[
{\bf\delta}_{i,j}=
\begin{cases}
\unit_{m_i}, &\mbox{if $i=j$}, \\
{\pmb 0}_{m_i\times m_j} &\mbox{if $i\ne j$}.
\end{cases}
\]
For $j=1$, \eqref{e1.2} becomes
\begin{equation}
\lambda f_{1,1}(\lambda)-u_1f_{2,1}(\lambda)=\unit_{m_1},
\label{e1.3}
\end{equation}
and for $j$ in $\{2,3,\ldots,d-1\}$, we get
\begin{equation}
f_{j,1}(\lambda)-u_jf_{j+1,1}(\lambda)={\pmb 0}_{m_j\times m_1}, 
\quad \mbox{i.e.,} \quad f_{j,1}(\lambda)=u_jf_{j+1,1}(\lambda).
\label{e1.4}
\end{equation}
Finally, for $j=d$, \eqref{e1.2} yields
\begin{equation}
-u_df_{1,1}(\lambda)+f_{d,1}(\lambda)={\pmb 0}_{m_d\times m_1}
\quad \mbox{i.e.,} \quad f_{d,1}(\lambda)=u_df_{1,1}(\lambda).
\label{e1.5}
\end{equation}
Then, by successive applications of the formulae \eqref{e1.4} and
\eqref{e1.5}, we find that
\[
f_{2,1}(\lambda)=u_2f_{3,1}(\lambda)=u_2u_3f_{4,1}(\lambda)
=\cdots=u_2u_3\cdots u_{d-1}f_{d,1}(\lambda)
=u_2u_3\cdots u_df_{1,1}(\lambda).
\]
Inserting this in \eqref{e1.4}, we obtain
\[
\unit_{m_1}=\lambda f_{1,1}(\lambda)-u_1\big(u_2u_3\cdots
u_df_{1,1}(\lambda)\big)=\big(\lambda-u_1u_2\cdots
u_d\big)f_{1,1}(\lambda).
\]
To verify that also $f_{1,1}(\lambda)(\lambda-u_1u_2\cdots
u_d\big)=\unit_{m_1}$, we consider the equality
$F(\lambda)A(\lambda)=\unit_k$, from which
\[
\sum_{i=1}^df_{1,i}(\lambda)a_{i,j}(\lambda)=\delta_{1,j}, \qquad
(j=1,2,\ldots d).
\]
For $j=1$ we obtain
\begin{equation}
f_{1,1}(\lambda)\lambda-f_{1,d}(\lambda)u_d=\unit_{m_1},
\label{e1.6}
\end{equation}
and for $j$ in $\{2,3,\ldots,d\}$,
\begin{equation}
f_{1,j-1}(\lambda)u_{j-1}+f_{1,j}(\lambda)={\pmb 0}_{m_1\times m_j}
\quad \mbox{i.e.,} \quad f_{1,j}(\lambda)=f_{1,j-1}(\lambda)u_{j-1}.
\label{e1.7}
\end{equation}
By successive applications of \eqref{e1.7},
\[
f_{1,d}(\lambda)=f_{1,d-1}(\lambda)u_{d-1}=f_{1,d-2}(\lambda)u_{d-2}u_{d-1}
=\cdots=f_{1,1}(\lambda)u_1u_2\cdots u_{d-1},
\]
and inserting this in \eqref{e1.6}, we obtain
\[
\unit_{m_1}=f_{1,1}(\lambda)\lambda-\big(f_{1,1}(\lambda)u_1u_2\cdots
u_{d-1}\big)u_d=f_{1,1}(\lambda)\big(\lambda-u_1u_2\cdots u_d\big),
\]
as desired.

Assume next that $(\lambda-u_1u_2\cdots u_d)$ is invertible in
$M_m(\CA)$ and consider the matrices $B(\lambda)$ and $C$
introduced in Proposition~\ref{inv-formel}. At first we show that
\[
A(\lambda)\big(B(\lambda)+C\big)=\unit_k.
\]
It is easily seen that 
\[
A(\lambda)
\begin{pmatrix}
\unit_m \\
u_2u_3\cdots u_d \\
u_3u_4\cdots u_d \\
u_4\cdots u_d \\
\vdots \\
u_d
\end{pmatrix}
=
\begin{pmatrix}
\lambda-u_1u_2\cdots u_d \\
0 \\
0 \\
0 \\
\vdots \\
0
\end{pmatrix}
\]
so that 
\[
\begin{split}
A(\lambda)B(\lambda)&=
\begin{pmatrix}
\unit_m\\
0 \\
0 \\
\vdots \\
0
\end{pmatrix}
\begin{pmatrix}
1 & u_1 & u_1u_2 & u_1u_2u_3 & \cdots & u_1u_2\cdots u_{d-1}
\end{pmatrix} \\[.2cm]
&=
\begin{pmatrix}
1 & u_1 & u_1u_2 & u_1u_2u_3 & \cdots & u_1u_2\cdots u_{d-1} \\
0 & 0 & 0 & 0 & \cdots & 0 \\
\vdots & \vdots & \vdots & \vdots & \ddots & \vdots \\
0 & 0 & 0 & 0 & \cdots & 0 \\
\end{pmatrix}.
\end{split}
\]
It thus remains to verify that
\begin{equation}\label{A(lambda)C}
A(\lambda)C
=
\begin{pmatrix}
0 & -u_1 & -u_1u_2 & -u_1u_2u_3 & \cdots & -u_1u_2\cdots u_{d-1} \\
0 & \unit_{m_2} & 0 & 0 & \cdots & 0 \\
0 & 0 & \unit_{m_3} & 0 & \cdots & 0 \\
0 & 0 & 0 & \unit_{m_4} & \cdots & 0 \\
\vdots & \vdots & \vdots & \vdots & \ddots & \vdots \\
0 & 0 & 0 & 0 & \cdots & \unit_{m_d} \\
\end{pmatrix}.
\end{equation}
To this end, note that the first column in $A(\lambda)C$
consists entirely of zeroes and that the second column in
$A(\lambda)C$ equals that of $A(\lambda)$.

Note next that
for $j$ in $\{3,4,\ldots,d\}$, the entry at position $(1,j)$ is
\[
\big[A(\lambda)C\big]_{1,j}=
\begin{pmatrix}
\lambda & -u_1 & 0 & \cdots & 0
\end{pmatrix}
\begin{pmatrix}
0 \\
u_2u_3\cdots u_{j-1} \\
u_3u_4\cdots u_{j-1} \\
\vdots \\
u_{j-1} \\
\unit_{m_j} \\
0 \\
\vdots \\
0
\end{pmatrix}
=-u_1u_2\cdots u_{j-1}.
\]
Next, if $i\in\{2,3,\ldots,d-1\}$ and $j\in\{3,4,\ldots,d\}$, then the
entry of $A(\lambda)C$ at position $(i,j)$ is
\[
\begin{split}
\big[A(\lambda)C\big]_{i,j}&=
\begin{pmatrix}
0 & \cdots & 0 & \unit_{m_i} & -u_i & 0 & \cdots & 0
\end{pmatrix}
\begin{pmatrix}
0 \\
u_2u_3\cdots u_{j-1} \\
u_3u_4\cdots u_{j-1} \\
\vdots \\
u_{j-1} \\
\unit_{m_j} \\
0 \\
\vdots \\
0
\end{pmatrix} \\[.2cm]
&=
\begin{cases}
0, &\mbox{if $i>j$}, \\
\unit_{m_i}, &\mbox{if $i=j$}, \\
u_iu_{i+1}\cdots u_{j-1}-u_iu_{i+1}\cdots u_{j-1}=0, &\mbox{if $i\le
  j-1$}.
\end{cases}
\end{split}
\]
Finally, for $j$ in $\{3,4,\ldots,k\}$, the entry at position $(d,j)$
is
\[
\big[A(\lambda)C\big]_{d,j}=
\begin{pmatrix}
-u_d & 0 & \cdots & 0 & \unit_{m_d}
\end{pmatrix}
\begin{pmatrix}
0 \\
u_2u_3\cdots u_{j-1} \\
u_3u_4\cdots u_{j-1} \\
\vdots \\
u_{j-1} \\
\unit_{m_j} \\
0 \\
\vdots \\
0
\end{pmatrix}=
\begin{cases}
0, &\mbox{if $j<d$}, \\
\unit_{m_d} &\mbox{if $j=d$}.
\end{cases}
\]
Hence \eqref{A(lambda)C} holds, and therefore
$A(\lambda)(B(\lambda)+C)=\unit_k$. 

To verify that also
\[
\big(B(\lambda)+C\big)A(\lambda)=\unit_k,
\]
at first note that
\[
\begin{pmatrix}
\unit_m & u_1 & u_1u_2 & u_1u_2u_3 & \cdots & u_1u_2\cdots u_{d-1}
\end{pmatrix}
A(\lambda)=
\begin{pmatrix}
\lambda-u_1u_2\ldots u_d & 0 & \cdots & 0
\end{pmatrix}
\]
so that
\[
B(\lambda)A(\lambda)=
\begin{pmatrix}
\unit_m \\
u_2u_3\cdots u_d \\
u_3u_4\cdots u_d \\
\vdots \\
u_d
\end{pmatrix}
\begin{pmatrix}
\unit_m & 0 & \cdots & 0
\end{pmatrix}=
\begin{pmatrix}
\unit_m & 0 & \cdots & 0 \\
u_2u_3\cdots u_d & 0 & \cdots & 0 \\
u_3u_4\cdots u_d & 0 & \cdots & 0\\
\vdots & \vdots & \ddots & \vdots \\
u_d & 0 & \cdots & 0
\end{pmatrix}.
\]
It thus remains to show that
\[
CA(\lambda)=
\begin{pmatrix}
0 & 0 & 0 & \cdots & 0 \\
-u_2u_3\cdots u_d & \unit_{m_2} & 0 & \cdots & 0 \\
-u_3u_4\cdots u_d & 0 & \unit_{m_3} & \cdots & 0\\
\vdots & \vdots & \vdots & \ddots & \vdots \\
-u_d & 0 & 0 & \cdots & \unit_{m_d}
\end{pmatrix}.
\]
This follows easily by considerations similar to those described above.
$\endproof$

\begin{cor}\label{invertibilitet}
Let $p$ be a polynomial in $ M_m(\C)\otimes
\C\<X_1,\ldots,X_r\>$ of degree $d$, and let $x_1, \ldots, x_r$ be elements
in a unital algebra $\CA$. As in Proposition~\ref{faktorisering}, choose a
factorization of $p$ into polynomials of first degree,
\[
p=u_1 u_2\cdots u_d.
\]
Put
\[
v_j=u_j(x_1,\ldots,x_r),\qquad (j=1,2,\ldots,d),
\]
and let $\lambda\in M_m(\C)$. Then $\lambda\otimes\unit_\CA-p(x_1,
\ldots,x_r)$ is invertible in $M_m(\CA)$ iff the matrix
\[
A(\lambda,v_1,\ldots,v_r)=
\begin{pmatrix}
\lambda\otimes\unit_n & -v_1 & 0 & 0 & \cdots & 0    \\
0 & \unit_{m_2}\otimes\unit_n & -v_2 & 0 & \cdots & 0 \\
0  &  0 & \unit_{m_3}\otimes\unit_n & -v_3 & \cdots & 0 \\
\vdots & \vdots & \vdots & \ddots & \ddots & \vdots \\
0 & 0 & 0 & \cdots & \unit_{m_{d-1}}\otimes\unit_n & -v_{d-1} \\
-v_d & 0 & 0 & \cdots & 0 & \unit_{m_d}\otimes\unit_n 
\end{pmatrix},
\]
is invertible in $M_k(\CA)$. 
\end{cor}

By application of Corollary~\ref{invertibilitet}, one can give a purely
algebraic proof of the following ``linearization trick'' which was obtained in
\cite{HT} by use of Stinespring's Theorem and Arveson's Extension Theorem
for completely positive maps:

\begin{cor}\cite[Theorem~2.2]{HT} Let $\CA$ and $\CB$ be unital $C\cc$-algebras, and let $x_1,
  \ldots, x_r\in \CA_{sa}$, $y_1, \ldots, y_r\in \CB_{sa}$. If for all
  $m\in\N$ and all $a_0, a_1,\ldots, a_r\in M_m(\C)_{sa}$ we have that
  \begin{equation}\label{eq2.10}
    \sigma\Big(a_0\otimes\unit_\CA + \sum_r a_i\otimes x_i\Big)\supseteq
    \sigma\Big(a_0\otimes\unit_\CB + \sum_r a_i\otimes y_i\Big),
  \end{equation}
  then there exists a unital $\ast$-homomorphism
  \[
  \phi : C\cc(\unit_\CA, x_1, \ldots, x_r)\rightarrow C\cc(\unit_\CB, y_1,
  \ldots, y_r),
  \]
  such that $\phi(x_i)=y_i$ for $i=1, \ldots, r$.
\end{cor}

\proof As in step~I of the proof of \cite[Theorem~2.2]{HT}, a simple
$2\times 2$-matrix argument shows that if \eqref{eq2.10} holds, then it
also holds for arbitrary elements $a_0, a_1,\ldots, a_r\in M_m(\C)$. That
is, for every polynomial $q$ of degree at most 1 in $M_m(\C)\otimes
\C\<X_1, \ldots, X_r\>$ one has that
\begin{equation}\label{eq2.11}
  \sigma(q(x_1, \ldots, x_r))\supseteq\sigma(q(y_1, \ldots, y_r)).
\end{equation}
Now, let $p\in\C\<X_1, \ldots, X_r\>$ be a polynomial of degree $d\geq 1$,
and as in Proposition~\ref{faktorisering} (with $m=1$), choose a factorization
\[
p=u_1u_2\cdots u_d.
\]
For $j=1, \ldots, d$, put
\[
v_j=u_j(x_1, \ldots, x_r)
\]
and
\[
w_j=u_j(y_1, \ldots, y_r).
\]
Then, with the notation of Corollary~\ref{invertibilitet}, for
$\lambda\in\C$ we have that
\[
\lambda\unit_\CA- p(x_1, \ldots,x_r)\in \CA_{inv} \Leftrightarrow
A(\lambda, v_1, \ldots, v_d)\in M_k(\CA)_{inv}
\]
and
\[
\lambda\unit_\CB- p(y_1, \ldots,y_r)\in \CB_{inv} \Leftrightarrow
A(\lambda, w_1, \ldots, w_d)\in M_k(\CB)_{inv}.
\]
Since the $u_i$'s have degree 1,
\[
q:= A(\lambda, u_1, \ldots, u_d)
\]
is a polynomial of degree 1 in $M_k(\C)\otimes\C\<X_1, \ldots, X_r\>$. Note
that $A(\lambda, v_1, \ldots, v_d)=q(x_1, \ldots, x_r)$ and $A(\lambda,
w_1, \ldots, w_d)= q(y_1, \ldots, y_r)$. Hence, by \eqref{eq2.11},
\[
\sigma(A(\lambda, v_1, \ldots, v_d))\supseteq \sigma(A(\lambda,
w_1, \ldots, w_d)).
\]
In particular,
\[
A(\lambda, v_1, \ldots, v_d)\in M_k(\CA)_{inv} \Rightarrow A(\lambda,
w_1, \ldots, w_d)\in M_k(\CB)_{inv}.
\]
Altogether, we have shown that
\[
\lambda\unit_\CA- p(x_1, \ldots,x_r)\in \CA_{inv} \Rightarrow
\lambda\unit_\CB- p(y_1, \ldots,y_r)\in \CB_{inv},
\]
i.e.
\begin{equation}\label{eq2.12}
  \sigma(p(x_1, \ldots, x_r))\supseteq \sigma(p(y_1, \ldots, y_r))
\end{equation}
holds for all polynomials $p\in \C\<X_1, \ldots, X_r\>$. In particular, the
spectral radii, $r(p(x_1, \ldots, x_r))$ and $r(p(y_1, \ldots, y_r))$ satisfy
\begin{equation}\label{eq2.13}
  r(p(x_1, \ldots, x_r))\geq r(p(y_1, \ldots, y_r)).
\end{equation}
Applying \eqref{eq2.13} to the self-adjoint polynomial $p\cc p$, we get
that
\[
\|p(x_1, \ldots, x_r)\|^2\geq \|p(y_1, \ldots, y_r)\|^2.
\]
Hence, the map
\[
\phi_0 : p(x_1, \ldots, x_r)\mapsto p(y_1, \ldots, y_r), \qquad (p\in
\C\<X_1, \ldots, X_r\>),
\]
is well-defined and extends by continuity to a unital $\ast$-homomorphism
$\phi$ from  $C\cc(\unit_\CA, x_1, \ldots, x_r)$ into $ C\cc(\unit_\CB, y_1,
\ldots, y_r)$ with $\phi(x_i)=y_i$, $i=1, \ldots, r$. $\endproof$ 

\section{Norm estimates.}
\label{sec3}

In this section we consider a fixed self-adjoint polynomial $p$ in $r$ non-commuting variables with
coefficients in $M_m(\C)$, i.e. $p\in (M_m(\C)\otimes\C\<X_1, \ldots,
X_r\>)_{sa}$, and for each $n\in\N$, we let $X_1^{(n)}, \ldots,X_r^{(n)}$ be
stochastically independent random matrices from $\SGRM(n, \frac1n)$. Define
self-adjoint random matrices $(Q_n)_{n=1}^\infty$ by
\begin{equation}\label{defQ_n}
Q_n(\omega)=p(X_1^{(n)}(\omega), \ldots,X_r^{(n)}(\omega)), \qquad (\omega\in\Omega),
\end{equation}

where $(\Omega,\CF,P)$ denotes the underlying probability space. 

With $d={\rm deg}(p)$ we may, according to Proposition~\ref{faktorisering}, choose $m_1,\ldots, m_{d+1}\in \N$ with
$m=m_1=m_{d+1}$, and polynomials $u_j\in M_{m_j,m_{j+1}}(\C)\otimes
\C\<X_1, \ldots, X_r\>$ of first degree, $ j=1,\ldots, d$, such that $p=u_1u_2\cdots u_d$. For each $n\in\N$ define random matrices $u_j^{(n)}$, $j=1,\ldots, d$, by
\[
u_j^{(n)}(\omega) = u_j(X_1^{(n)}(\omega), \ldots,X_r^{(n)}(\omega)), \qquad (\omega\in\Omega).
\]

For $\lambda\in M_m(\C)$ we put $\im\lambda =
\frac{1}{2\i}(\lambda-\lambda\cc)$ as in \cite[Section~3]{HT}.

Since $Q_n(\omega)$ is self-adjoint, $\lambda\otimes\unit_n - Q_n(\omega)$
is invertible for every $\lambda\in M_m(\C)$ with $\im\lambda$ positive definite
(cf. \cite[Lemma~3.1]{HT}). Then, according to Corollary~\ref{invertibilitet}, the random matrix
\begin{equation}\label{3.1}
A_n(\lambda) = \begin{pmatrix} 
\lambda\otimes\unit_n & -u_1^{(n)} & 0 & 0 & \cdots & 0\\
0 & \unit_{m_2}\otimes\unit_n & -u_2^{(n)} & 0 & \cdots & 0\\
0 & 0 & \unit_{m_3}\otimes\unit_n & -u_3^{(n)} & \cdots & 0\\
\vdots & \vdots & \vdots & \ddots & \ddots & \vdots\\
0 & 0 & 0 & \cdots & \unit_{m_{d-1}}\otimes\unit_n & -u_{d-1}^{(n)}\\
-u_{d}^{(n)} & 0 & 0 & \cdots & 0 & \unit_{m_d}\otimes\unit_n
\end{pmatrix}
\end{equation}

is (point-wise) invertible in $M_k(\C)$, where $k=\sum_{i=1}^d m_i$.

\begin{lemma}\label{normlemma} For every $p\in\N$ there exist constants
  $C_{1,p}, C_{2,p}\geq 0$, such that for all $m\in\N$ and for all
  $\lambda\in M_m(\C)$ with $\im\lambda$ positive definite,
\[
\sup_{n\in\N}\,\E\{\|A_n(\lambda)^{-1}\|^p\}\leq C_{1,p}+ C_{2,p}\|(\im\lambda)^{-1}\|^p.
\]
\end{lemma}

\proof Let $p\in\N$. According to Proposition~\ref{inv-formel} we may write 
\[
A_n(\lambda)^{-1} = C_n +  B_n^{(1)} (\lambda\otimes\unit_n -Q_n)^{-1}B_n^{(2)},
\]

where
\[
C_n=\begin{pmatrix}
0 & 0 & 0 & 0 & 0 & \hdots & 0\\
0 & \unit_{m_2}\otimes\unit_n & u_{2}^{(n)} & u_{2}^{(n)}u_{3}^{(n)} & u_{2}^{(n)}u_{3}^{(n)}u_{4}^{(n)} & \hdots & u_{2}^{(n)}u_{3}^{(n)}\cdots u_{d-1}^{(n)}\\
0 & 0 & \unit_{m_3}\otimes\unit_n & u_{3}^{(n)} & u_{3}^{(n)}u_{4}^{(n)} & \hdots & u_{3}^{(n)}u_{4}^{(n)}\cdots u_{d-1}^{(n)}\\
0 & 0 & 0 & \unit_{m_4}\otimes\unit_n & u_{4}^{(n)} & \hdots & u_{4}^{(n)}u_{5}^{(n)}\cdots u_{d-1}^{(n)}\\
\vdots & \vdots & \vdots & \vdots & \ddots & \ddots &\vdots\\
0 & 0 & 0 & 0 & \hdots & \unit_{m_{d-1}}\otimes\unit_n & u_{d-1}^{(n)}\\
0 & 0 & 0 & 0 & \hdots & 0 & \unit_{m_d}\otimes\unit_n
\end{pmatrix},
\]
\[
 B_n^{(1)} = \begin{pmatrix}
\unit_m\otimes\unit_n \\
u_{2}^{(n)}u_{3}^{(n)}\cdots u_{d}^{(n)}\\
u_{3}^{(n)}u_{4}^{(n)}\cdots u_{d}^{(n)}\\
\vdots\\
u_{d}^{(n)}
\end{pmatrix},
\]

and
\[
B_n^{(2)} = \begin{pmatrix}
\unit_m\otimes\unit_n & u_{1}^{(n)} & u_{1}^{(n)}u_{2}^{(n)} & \hdots & u_{1}^{(n)}u_{2}^{(n)}\cdots u_{d-1}^{(n)}\end{pmatrix}.  
\]

By \cite[Lemma~3.1]{HT}, $\|(\lambda\otimes\unit_n -Q_n)^{-1}\|\leq
\|(\im\lambda)^{-1}\|$, and therefore
\begin{eqnarray}\label{primary}
\E\{\|A_n(\lambda)^{-1}\|^p\}&\leq& \E\{(\|C_n\|+ \|B_n^{(1)}\|\|B_n^{(2)}\|\|(\im\lambda)^{-1}\|)^p\}\nonumber\\
&\leq & \E\{(2\,
\max\{\|C_n\|,\|B_n^{(1)}\|\|B_n^{(2)}\|\|(\im\lambda)^{-1}\|
\})^p\}\nonumber\\
& \leq & 2^p \, \E\{\|C_n\|^p + \|B_n^{(1)}\|^p\|B_n^{(2)}\|^p\|(\im\lambda)^{-1}\|^p\}.
\end{eqnarray}

With 
\begin{equation}\label{Kn}
K_n=\max\{1,\|u_{1}^{(n)}\|,\|u_{2}^{(n)}\|,\ldots, \|u_{d}^{(n)}\|\}
\end{equation}

one easily proves that
\[
\|C_n\|\leq d-1 + d^2 K_n^{d-2}\leq d^2(1+K_n^d),
\]

implying that
\begin{equation}\label{normBn}
\|C_n\|^p \leq 2^p d^{2p}(1+K_n^{dp})\leq 2^pd^{2p}(1+K_n^{2pd}).
\end{equation}

Moreover,
\begin{equation}\label{normCn}
\|B_n^{(1)}\|^p\|B_n^{(2)}\|^p\leq d^{2p}K_n^{2p(d-1)}\leq d^{2p}K_n^{2pd}.
\end{equation}

Now, according to (\ref{Kn}),
\begin{equation}\label{Kn2pd}
K_n^{2pd}\leq  1+\sum_{j=1}^d\|u_{j}^{(n)}\|^{2pd}.
\end{equation}

Since $u_{j}^{(n)}$ is of first degree, we may choose $a_0^{(j)}, \ldots, a_r^{(j)}\in M_{m_j,m_{j+1}}(\C)$ such that
\[
u_{j}^{(n)}= a_0^{(j)}\otimes\unit_n + \sum_{i=1}^r a_i^{(j)}\otimes X_i^{(n)}. 
\]

Hence 
\begin{equation}\label{normun}
\|u_{j}^{(n)}\|^{2pd}\leq (1+r)^{2pd}\, \max\{\|a_0^{(j)}\|,
\|a_1^{(j)}\|\| X_1^{(n)}\|, \ldots ,\|a_r^{(j)}\|\| X_r^{(n)}\|\}^{2pd}.
\end{equation}

According to [S, Lemma~6.4],
\[
\sup_n(\E\{\| X_i^{(n)}\|^{2pd}\})<\infty,
\]

and combining this fact with (\ref{primary}), (\ref{normBn}),
(\ref{normCn}), (\ref{Kn2pd}) and (\ref{normun}) we obtain the desired estimate. $\endproof$



\begin{lemma}\label{estimat of A(lambda)}
Let $x_1,\ldots, x_r$ be a semicircular system in a $C\cc$-probability
space $(\CA,\tau)$, and let
\begin{equation}\label{defq}
q=p(x_1,\ldots,x_r).
\end{equation}
Then again, with
$u_j=u_j(x_1,\ldots, x_r)$, $j=1,\ldots, d$, the matrix
\[
A(\lambda)= \begin{pmatrix} 
\lambda\otimes\unit_\CA & -u_1 & 0 & 0 & \cdots & 0\\
0 & \unit_{m_2}\otimes\unit_\CA & -u_2 & 0 & \cdots & 0\\
0 & 0 & \unit_{m_3}\otimes\unit_\CA & -u_3 & \cdots & 0\\
\vdots & \vdots & \vdots & \ddots & \ddots & \vdots\\
0 & 0 & 0 & \cdots & \unit_{m_{d-1}}\otimes\unit_\CA & -u_{d-1}\\
-u_{d} & 0 & 0 & \cdots & 0 & \unit_{m_d}\otimes\unit_\CA
\end{pmatrix}
\]

is invertible for every $\lambda\in M_m(\C)$ with $\im\lambda>0$. Moreover,
for every $p\in\N$ there exist constants $C_{1,p}'$ and $C_{2,p}'$ such that 
\begin{equation}\label{normAinv}
\|A(\lambda)^{-1}\|^p\leq C_{1,p}'+C_{2,p}'\|(\im\lambda)^{-1}\|^p
\end{equation}

holds for every $\lambda\in M_m(\C)$ with $\im\lambda$ positive definite.
\end{lemma}

\proof 
It follows again form \cite[Lemma~3.1]{HT} that $\lambda\otimes\unit_\CA -
q$ is invertible. Hence, by Proposition~\ref{inv-formel}, $A(\lambda)$ is
invertible with
\begin{equation*}
A(\lambda)^{-1} =  C +  B^{(1)} (\lambda\otimes\unit_{\CA} -q)^{-1}B^{(2)},
\end{equation*}

where 
\[
C=\begin{pmatrix}
0 & 0 & 0 & 0 & 0 & \hdots & 0\\
0 & \unit_{m_2}\otimes\unit_{\CA} & u_{2} & u_{2}u_{3} & u_{2}u_{3}u_{4} & \hdots & u_{2}u_{3}\cdots u_{d-1}\\
0 & 0 & \unit_{m_3}\otimes\unit_{\CA} & u_{3} & u_{3}u_{4} & \hdots & u_{3}u_{4}\cdots u_{d-1}\\
0 & 0 & 0 & \unit_{m_4}\otimes\unit_{\CA} & u_{4} & \hdots & u_{4}u_{5}\cdots u_{d-1}\\
\vdots & \vdots & \vdots & \vdots & \ddots & \ddots &\vdots\\
0 & 0 & 0 & 0 & \hdots & \unit_{m_{d-1}}\otimes\unit_{\CA} & u_{d-1}\\
0 & 0 & 0 & 0 & \hdots & 0 & \unit_{m_d}\otimes\unit_{\CA}
\end{pmatrix},
\]
\[
 B^{(1)} = \begin{pmatrix}
\unit_m\otimes\unit_{\CA} \\
u_{2}u_{3}\cdots u_{d}\\
u_{3}u_{4}\cdots u_{d}\\
\vdots\\
u_{d}
\end{pmatrix},
\]

and 
\[
B^{(2)} = \begin{pmatrix}
\unit_m\otimes\unit_{\CA} & u_{1} & u_{1}u_{2} & \hdots & u_{1}u_{2}\cdots u_{d-1}\end{pmatrix}.  
\]

Then, since $\| (\lambda\otimes\unit_{\CA} -q)^{-1}\|\leq
\|(\im\lambda)^{-1}\|$, we have as in the proof of Lemma~\ref{normlemma} that
\[
\|A(\lambda)^{-1}\|^p\leq 2^p(\|C\|^p + \|B^{(1)}\|^p\|B^{(2)}\|^p \|(\im\lambda)^{-1}\|^p),
\]

and the claim follows. $\endproof$

\section{The master equation and master inequality.}\label{sec4}

In this section we prove generalizations of the master equation and master
inequality from \cite{HT}. These generalizations will allow us to handle
self-adjoint polynomials $p(X_1^{(n)}, \ldots, X_r^{(n)})$ of arbitrary
degree in $r$ independent random matrices from $\SGRM(n,\frac 1n)$. 

Let $r$ and $n$ be positive integers. As in \cite[Section~3]{HT}, we shall consider
the {\it real} vector space $(\Mnsa)^r$, which we denote by $\CE_{r,n}$.
We equip $\CE_{r,n}$ with the inner product $\<\cdot,\cdot\>_e$ given by
\[
\<(A_1,\ldots,A_r),(B_1,\ldots,B_r)\>_{e}
=\Tr_n\Big(\sum_{j=1}^rA_jB_j\Big), \qquad
((A_1,\ldots,A_r),(B_1,\ldots,B_r)\in\CE_{r,n}),
\]
and we denote the corresponding norm by $\|\cdot\|_e$.
Still following \cite{HT}, we consider the linear isomorphism
$\Psi_0$ between $\Mnsa$ and $\R^{n^2}$ given by
\begin{equation}
\Psi_0((a_{uv})_{1\le u,v\le n})=\big((a_{uu})_{1\le u\le
n},(\sqrt{2}\re(a_{uv}))_{1\le u<v\le n},
(\sqrt{2}\im(a_{ub}))_{1\le u<v\le n}\big),
\label{e2.0}
\end{equation}
for $(a_{uv})_{1\le u,v\le n}$ in $\Mnsa$.
We consider further the natural extension $\Psi\colon\CE_{r,n}\to\R^{rn^2}$
of $\Psi_0$ given by
\[
\Psi(A_1,\ldots,A_r)=(\Psi_0(A_1),\ldots,\Psi_0(A_r)), \qquad
(A_1,\ldots,A_r\in\Mnsa).
\]
We note that $\Psi$ is an isometry between $(\CE_{r,n},\|\cdot\|_e)$ and
$\R^{rn^2}$ equipped with its usual Hilbert space norm. Accordingly, we
shall identify $\CE_{r,n}$ with $\R^{rn^2}$ via $\Psi$.

In the following we consider a fixed self-adjoint polynomial $p$ from
$M_m(\C)\otimes\C\<X_1,\ldots,X_r\>$ of degree $d$ and the corresponding
polynomials
\[
u_j=u_j(X_1,\ldots,X_r)\in M_{m_j,m_{j+1}}(\C)\otimes\C\<X_1,\ldots,X_r\>,
\qquad (j=1,2,\ldots,d),
\]
introduced in Proposition~\ref{faktorisering}. We put $k=m_1+m_2+\cdots m_d$.

We consider further independent random matrices $X_1^{(n)},\ldots,X_r^{(n)}$ from
$\SGRM(n,\frac{1}{n})$ and a fixed matrix $\lambda$ from $M_m(\C)$, such
that $\im(\lambda)$ is positive definite. We may then consider the (random)
matrix $A(\lambda,X_1^{(n)},\ldots,X_r^{(n)})$ defined in
Corollary~\ref{invertibilitet}. Since the polynomials
$u_1,\ldots,u_d$ are of degree 1, we may write
\[
A(\lambda,X_1^{(n)},\ldots,X_r^{(n)})=
\begin{pmatrix}
\lambda & 0 \\
0 & \unit_{k-m}
\end{pmatrix}
\otimes\unit_n
-a_0\otimes \unit_n
-\sum_{j=1}^ra_i\otimes X_i^{(n)},
\]
for suitable matrices $a_0, a_1,\ldots,a_r$ in $M_k(\C)$. We put
\begin{equation}
\Lambda=
\begin{pmatrix}
\lambda & 0 \\
0 & \unit_{k-m}
\end{pmatrix}
\quad \mbox{and} \quad
S_n=a_0\otimes \unit_n + \sum_{j=1}^ra_j\otimes X_j^{(n)},
\label{e2.1a}
\end{equation}
so that
\[
A(\lambda,X_1^{(n)},\ldots,X_r^{(n)})=\Lambda\otimes\unit_n-S_n.
\]
According to Corollary~\ref{invertibilitet}, $\Lambda\otimes\unit_n-S_n$ is
invertible, and hence we may consider the $k\times k$ matrix
\begin{equation}
H_n(\lambda)=(\id_k\otimes\tr_n)\big[(\Lambda\otimes\unit_n-S_n)^{-1}\big].
\label{e2.1}
\end{equation}

The following lemma generalizes \cite[Lemma~3.5]{HT}:

\begin{lemma}\label{lemma til master eq} With $H_n(\lambda)$ as defined in
  \eqref{e2.1}, we have for any $j$ in $\{1,2,\ldots,r\}$ the formula
\[
\E\big\{H_n(\lambda)a_jH_n(\lambda)\big\}=
\E\big\{ (\id_k\otimes\tr_n)\big[(\unit_k\otimes 
X_j^{(n)})\cdot(\Lambda\otimes\unit_n-S_n)^{-1}\big]\big\}.
\]
\end{lemma}

\proof For any $v_1,\ldots,v_r\in\Mnsa$ we may consider the
matrix $A(\lambda,v_1,\ldots,v_r)$ described in
Corollary~\ref{invertibilitet}, and we clearly have that
\[
A(\lambda,v_1,\ldots,v_r)=\Lambda\otimes\unit_n--a_0\otimes\unit_n - \sum_{j=1}^ra_j\otimes
v_j,
\]
with $\Lambda,a_1,\ldots,a_r$ as above. According to
Corollary~\ref{invertibilitet}, we may then consider the mapping
$\tilde{F}\colon\CE_{r,n}\to M_k(\C)\otimes M_n(\C)$ given by
\[
\tilde{F}(v_1,\ldots,v_r)=\big((\Lambda-a_0)\otimes\unit_n-
\textstyle{\sum_{j=1}^ra_j\otimes v_j}\big)^{-1}, \qquad
((v_1,\ldots,v_r)\in\CE_{r,n}).
\]
We consider furthermore the mapping $F\colon\R^{rn^2}\to M_k(\C)\otimes
M_n(\C)$, given by
\[
F=\tilde{F}\circ\Psi^{-1}.
\]
Note then that
\begin{equation}
(\Lambda\otimes\unit_n-S_n)^{-1}= \tilde F(X_1^{(n)},\ldots,X_r^{(n)})
=F(\Psi(X_1^{(n)},\ldots,X_r^{(n)})),
\label{e2.2}
\end{equation}
where
$\Psi(X_1^{(n)},\ldots,X_r^{(n)})=(\gamma_1,\gamma_2,\ldots,\gamma_{rn^2})$
with $\gamma_1,\gamma_2,\ldots,\gamma_{rn^2}\sim \mbox{i.i.d.\ } \
N(0,\frac{1}{n})$. 

According to the proof of \cite[Lemma~3.1]{HT}, we have that
\begin{equation}
\big\|\tilde{F}(v_1,\ldots,v_r)\big\|\le h(\|v_1\|,\ldots,\|v_r\|)(1+\|(\im\lambda)^{-1}\|),
\label{e2.3}
\end{equation}
for some polynomial $h$ in $\C[X_1,\ldots,X_r]$. From \eqref{e2.2} and \eqref{e2.3} it follows firstly that that the
expectations in Lemma~\ref{lemma til master eq} are well-defined. In
addition, \eqref{e2.3} shows that the function $F$ is a polynomially
bounded function of $rn^2$ real variables. In order to apply
\cite[Lemma~3.3]{HT}, we need to check that the
partial derivatives of $F$ are polynomially bounded as well. To this end,
consider the standard orthonormal basis for $\Mnsa$:
\[
\begin{split}
&e_{u,u}^{(n)}, \qquad (1\le u\le n) \\[.2cm]
&f_{u,v}^{(n)}=\textstyle{\frac{1}{\sqrt{2}}}\big(e_{u,v}^{(n)}+e_{v,u}^{(n)}\big)
\qquad (1\le u<v\le n), \\[.2cm]
&g_{u,v}^{(n)}=\textstyle{\frac{\i}{\sqrt{2}}}\big(e_{u,v}^{(n)}-e_{v,u}^{(n)}\big)
\qquad (1\le u<v\le n),
\end{split}
\]
where $\{e_{u,v}^{(n)}\mid 1\le u,v\le n\}$ are the standard $n\times n$
matrix units. The correponding orthonormal basis for $\CE_{r,n}$ is
\[
\begin{split}
e_{j,u,u}^{(n)}&=(0,\ldots,0,e_{u,u}^{(n)},0,\ldots,0), \qquad (1\le j\le r, \
1\le u\le n) \\[.2cm] 
f_{j,u,v}^{(n)}&=(0,\ldots,0,f_{u,v}^{(n)},0,\ldots,0),
\qquad (1\le j\le r, \ 1\le u<v\le n), \\[.2cm]
g_{j,u,v}^{(n)}&=(0,\ldots,0,g_{u,v}^{(n)},0,\ldots,0),
\qquad (1\le j\le r, \ 1\le u<v\le n),
\end{split}
\]
with the non-zero entry in the $j$'th slot. Note that the images by $\Psi$ of these
basis vectors is exactly the standard orthonormal basis for
$\R^{rn^2}$. Hence, the partial derivatives of $F$ at a point $\xi$ in
$\R^{rn^2}$ are, setting $(v_1,\ldots,v_r)=\Psi^{-1}(\xi)$,
\begin{equation}
\begin{split}
\diff F\big(\xi+& t\Psi(e_{j,u,u}^{(n)})\big)\\
&=\diff \tilde{F}\big((v_1,\ldots,v_r)+te_{j,u,u}^{(n)}\big) \\[.2cm]
&=\diff\tilde{F}\big(v_1,\ldots,v_{j-1},v_j+te_{u,u}^{(n)},v_{j+1},\ldots,v_r\big)
\\[.2cm]
&=\diff \big((\Lambda-a_0)\otimes\unit_n-\textstyle{\sum_{i=1}^ra_i\otimes
  v_i}-t(a_j\otimes e_{u,u}^{(n)})\big)^{-1} \\[.2cm]
&=\big((\Lambda-a_0)\otimes\unit_n-\textstyle{\sum_{i=1}^ra_i\otimes
  v_i}\big)^{-1}\big(a_j\otimes e_{u,u}^{(n)}\big)
\big((\Lambda-a_0)\otimes\unit_n-\textstyle{\sum_{i=1}^ra_i\otimes
  v_i})\big)^{-1},
\end{split}
\label{e2.4}
\end{equation}
where the last equality uses \cite[Lemma 3.2]{HT}. We find similarly that
\[
\begin{split}
\diff F\big(\xi+&t\Psi(f_{j,u,v}^{(n)})\big)=\\
&\big((\Lambda-a_0)\otimes\unit_n-\textstyle{\sum_{i=1}^ra_i\otimes
  v_i}\big)^{-1}\big(a_j\otimes f_{u,v}^{(n)}\big)
\big((\Lambda-a_0)\otimes\unit_n-\textstyle{\sum_{i=1}^ra_i\otimes
  v_i}\big)^{-1}, \\[.2cm]
\diff F\big(\xi+&t\Psi(g_{j,u,v}^{(n)})\big)=\\
&\big((\Lambda-a_0)\otimes\unit_n-\textstyle{\sum_{i=1}^ra_i\otimes
  v_i}\big)^{-1}\big(a_j\otimes g_{u,v}^{(n)}\big)
\big((\Lambda-a_0)\otimes\unit_n-\textstyle{\sum_{i=1}^ra_i\otimes
  v_i}\big)^{-1}. \\[.2cm]
\end{split}
\]
Appealing once more to \eqref{e2.3}, it follows that the partial
derivatives of $F$ are polynomially bounded as well. Hence, we may apply
\cite[Lemma~3.3]{HT} to $F$ and the i.i.d.\ Gaussian variables
$\Psi(X_1^{(n)},\ldots,X_r^{(n)})=(\gamma_1,\gamma_2,\ldots,\gamma_{rn^2})$.
For any $j$ in $\{1,2,\ldots,r\}$ put
\[
\begin{split}
X_{j,u,u}^{(n)}&=(X_j^{(n)})_{kk}, \qquad (1\le u\le n), \\[.2cm]
Y_{j,u,v}^{(n)}&=\sqrt{2}\re(X_j^{(n)})_{u,v}, \qquad (1\le u<v\le n),
\\[.2cm] 
Z_{j,kl}^{(n)}&=\sqrt{2}\im(X_j^{(n)})_{u,v}, \qquad (1\le u<v\le n),
\end{split}
\]
and note that these random variables are the coefficients of
$\Psi(X_1^{(n)},\ldots,X_r^{(n)})$ w.r.t.\ the standard
orthonormal basis for $\R^{rn^2}$, in the sense that
\[
\Psi(X_1^{(n)},\ldots,X_r^{(n)})=\sum_{j=1}^r\Big(\sum_{u=1}^nX_{j,u,u}^{(n)}\Psi(e_{j,u,u}^{(n)})
+\sum_{1\le u<v\le n}Y_{j,u,v}^{(n)}\Psi(f_{j,u,v}^{(n)})
+\sum_{1\le u<v\le n}Z_{j,u,v}^{(n)}\Psi(g_{j,u,v}^{(n)})\Big).
\]
It follows thus from the Gaussian Poincar\'e Inequality and
\eqref{e2.4} that
\[
\begin{split}
\E\big\{X_{j,u,u}^{(n)}\big(\Lambda\otimes\unit_n-S_n\big)^{-1}\big\}&=
\E\big\{X_{j,u,u}^{(n)}F(\Psi(X_1^{(n)},\ldots,X_r^{(n)}))\big\} \\[.2cm]
&=\frac{1}{n}\E\Big\{\diff
F(\Psi(X_1^{(n)},\ldots,X_r^{(n)}))-t\Psi(e_{j,u,u})\Big\} \\[.2cm]
&=\frac{1}{n}\E\big\{\big(\Lambda\otimes\unit_n-S_n\big)^{-1}\big(a_j\otimes e_{u,u}^{(n)}\big)
\big(\Lambda\otimes\unit_n-S_n\big)^{-1}\big\},
\end{split}
\]
and similarly we get that
\[
\begin{split}
&\E\big\{Y_{j,u,v}^{(n)}\cdot\big(\Lambda\otimes\unit_n-S_n\big)^{-1}\big\}=\frac{1}{n}\E\big\{(\Lambda\otimes\unit_n-S_n)^{-1}(a_j\otimes
f_{u,v}^{(n)})(\Lambda\otimes\unit_n-S_n)^{-1}\big\}
\\[.2cm]
&\E\big\{Z_{j,u,v}^{(n)}\cdot\big(\Lambda\otimes\unit_n-S_n\big)^{-1}\big\}=\frac{1}{n}\E\big\{(\Lambda\otimes\unit_n-S_n)^{-1}(a_j\otimes
g_{u,v}^{(n)})(\Lambda\otimes\unit_n-S_n)^{-1}\big\}.
\end{split}
\]
>From this point, the proof is completed exactly as in the proof of
\cite[Lemma~3.5]{HT}. 
$\endproof$

Lemma~\ref{lemma til master eq} implies the following analogue of
\cite[Theorem~3.6]{HT}. The proof is the same as in  \cite{HT} and will
therefore be omitted.

\begin{thm}[Master equation]\label{master eq}
Let $\lambda$ be a matrix in $M_m(\C)$ such that
$\im(\lambda)$ is positive definite, and let $\Lambda$ and 
$S_n$ be the matrices introduced in \eqref{e2.1a}. Then with
\[
H_n(\lambda)=(\id_k\otimes\tr_n)
\big[(\Lambda\otimes\unit_n-S_n)^{-1}\big]
\]
we have the formula
\begin{equation}
\E\Big\{\sum_{i=1}^r a_iH_n(\lambda)a_iH_n(\lambda) +(a_0-\Lambda) H_n(\lambda)
+\unit_m\Big\}=0.
\label{e2.10}
\end{equation}
\end{thm}

We next prove the following analogue of \cite[Theorem~4.5]{HT}:

\begin{thm}[Master inequality]\label{master ineq}
Let $\lambda$ be a matrix in $M_m(\C)$ such that
$\im(\lambda)$ is positive definite, and let $\Lambda$ and 
$S_n$ be the matrices introduced in \eqref{e2.1a}. Then with
\[
H_n(\lambda)=(\id_k\otimes\tr_n)
\big[(\Lambda\otimes\unit_n-S_n)^{-1}\big]
\]
and
\[
G_n(\lambda)=\E\{H_n(\lambda)\},
\]
we have the estimate
\[
\Big\|\sum_{i=1}^ra_iG_n(\lambda)a_iG_n(\lambda)+(a_0-\Lambda) G_n(\lambda) +
\unit_k\Big\|\le\frac{C}{n^2}
\Big(C_{1,4}+C_{2,4}\big\|(\im\lambda)^{-1}\big\|^4\Big),
\]
where $C=k^3(\sum_{i=1}^r\|a_i\|^2)^2$ and 
$C_{1,4},C_{2,4}$ are the constants introduced in Lemma~\ref{normlemma}.
\end{thm}

\proof Setting
$K_n(\lambda)=H_n(\lambda)-G_n(\lambda)=H_n(\lambda)-\E\{H_n(\lambda)\}$,
we find exactly as in \cite[proof of Theorem~4.5]{HT} that
\[
\sum_{i=1}^ra_iG_n(\lambda)a_iG_n(\lambda)+(a_0-\Lambda) G_n(\lambda) +
\unit_k
=-\E\Big\{\sum_{i=1}^ra_iK_n(\lambda)a_iK_n(\lambda)\Big\},
\]
and from this
\begin{equation}
\begin{split}
\Big\|\sum_{i=1}^ra_iG_n(\lambda)a_iG_n(\lambda)+(a_0-\Lambda) G_n(\lambda) +
\unit_k\Big\|
&\le\sum_{i=1}^r\|a_i\|^2\E\big\{\|K_n(\lambda)\|^2\big\} \\[.2cm]
&\le\Big(\sum_{i=1}^r\|a_i\|^2\Big)
\sum_{u,v=1}^k\E\big\{|K_{n,u,v}(\lambda)|^2\big\} \\[.2cm]
&=\Big(\sum_{i=1}^r\|a_i\|^2\Big)
\sum_{u,v=1}^k\V\big\{H_{n,u,v}(\lambda)\big\},
\end{split}
\label{e2.5}
\end{equation}
where $K_{n,u,v}(\lambda)$ (resp.\ $H_{n,u,v}(\lambda)$), $1\le u,v\le
k$, are the entries of $K_n(\lambda)$ (resp.\ $H_n(\lambda)$). 
As in \cite[proof of Theorem~4.5]{HT} we note that
\[
H_{n,u,v}(\lambda)=f_{n,u,v}\big(X_1^{(n)},\ldots,X_r^{(n)}\big),
\]
where $f_{n,u,v}\colon\CE_{r,n}\to\C$ is the function given by
\[
f_{n,u,v}(v_1,\ldots,v_r)
=k(\tr_k\otimes\tr_n)\big[(e_{u,v}^{(k)}\otimes\unit_n)
\big((\Lambda-a_0)\otimes\unit_n-\textstyle{\sum_{i=1}^ra_i\otimes v_i}\big)^{-1}\big],
\]
for $v=(v_1,\ldots,v_r)\in\CE_{r,n}$. For any unit vector
$w=(w_1,\ldots,w_r)$ from $\CE_{r,n}$, we find as in \cite{HT} that
\[
\Big|\diff f_{n,u,v}(v+tw)\Big|
\le\frac{1}{n}\textstyle{\big\|\sum_{i=1}^ra_i\otimes
  w_i\big\|^2_{2,\Tr_k\otimes\Tr_n}
\big\|((\Lambda-a_0)\otimes\unit_n-\sum_{i=1}^ra_i\otimes v_i)^{-1}\big\|^4},
\]
and here, by arguing as in the proof of \cite[Lemma~4.4]{HT},
\[
\textstyle{\big\|\sum_{i=1}^ra_i\otimes
  w_i\big\|^2_{2,\Tr_k\otimes\Tr_n}
\le k\big\|\sum_{i=1}^ra_i^*a_i\big\|
\le k\sum_{i=1}^r\|a_i\|^2},
\]
so that
\[
\Big|\diff f_{n,u,v}(v+tw)\Big|
\le\frac{k}{n}\big(\textstyle{\sum_{i=1}^r\|a_i\|^2\big)
\big\|((\Lambda-a_0)\otimes\unit_n-\sum_{i=1}^ra_i\otimes v_i)^{-1}\big\|^4}.
\]
Consequently,
\[
\begin{split}
\big\|\grad f_{n,u,v}(v)\big\|^2
&=\max\Big\{\Big|\diff f_{n,u,v}(v+tw)\Big|^2 \Bigm | w\in\CE_{r,n}, \
\|w\|_e=1\Big\} \\[.2cm]
&\le\frac{k}{n}\big(\textstyle{\sum_{i=1}^r\|a_i\|^2\big)
\big\|((\Lambda-a_0)\otimes\unit_n-\sum_{i=1}^ra_i\otimes v_i)^{-1}\big\|^4}.
\end{split}
\]
Combining this with \eqref{e2.3}, it is clear
that $\grad f_{n,u,v}$ (as well as $f_{n,u,v}$ itself) is polynomially
bounded as a function of $rn^2$ real variables. Hence we may apply the
Gaussian Poincar\'e Inequality in the form of \cite[Corollary~4.2]{HT} as
follows:
\begin{equation}
\begin{split}
\V\big\{H_{n,u,v}(\lambda)\big\}&=
\V\big\{f_{n,u,v}(X_1^{(n)},\ldots,X_r^{(n)})\big\}
\le\frac{1}{n}\E\left\{\big\|\grad
f_{n,u,v}(X_1^{(n)},\ldots,X_r^{(n)})\big\|^2\right\} \\[.2cm]
&\le\frac{k}{n^2}\big(\textstyle{\sum_{i=1}^r\|a_i\|^2}\big)
\E\left\{\big\|\big((\Lambda-a_0)\otimes\unit_n-\textstyle{\sum_{i=1}^ra_i\otimes
X_i^{(n)}}\big)^{-1}\big\|^4\right\}
\\[.2cm]
&\le\frac{k}{n^2}\big(\textstyle{\sum_{i=1}^r\|a_i\|^2}\big)
\cdot\Big(C_{1,4}+C_{2,4}\big\|(\im\lambda)^{-1}\big\|^4\Big),
\end{split}
\label{e2.6}
\end{equation}
where $C_{1,4}$ and $C_{2,4}$ are the constants given in
Lemma~\ref{normlemma}. Since \eqref{e2.6} holds for all $u,v$ in
$\{1,2,\ldots,k\}$, we find in combination with \eqref{e2.5} that
\[
\Big\|\sum_{i=1}^ra_iG_n(\lambda)a_iG_n(\lambda)+(a_0-\Lambda) G_n(\lambda) +
\unit_k\Big\|\le\frac{k^3}{n^2}\Big(\sum_{i=1}^r\|a_i\|^2\Big)^2
\cdot\Big(C_{1,4}+C_{2,4}\big\|(\im\lambda)^{-1}\big\|^4\Big),
\]
and this is the desired estimate.
$\endproof$

\section{Estimation of $\|G_n(\lambda)-G(\lambda)\|$.}\label{sec4a}

As in the two previous sections, for each $n\in\N$ we consider
stochastically independent random matrices $X_1^{(n)}, \ldots, X_r^{(n)}$
from $\GUE(n, \frac 1n)$, and we let
\begin{equation}
  Q_n = p(X_1^{(n)}, \ldots, X_r^{(n)}),
\end{equation}
where $p$ is a fixed self-adjoint polynomial from $M_m(\C)\otimes \C\<X_1,
\ldots, X_r\>$. We let $A_n(\lambda)$ be given by \eqref{3.1}, where
$\lambda\in M_m(\C)$, and $\im\lambda$ is positive definite. Then
we may write
\begin{equation}
A_n(\lambda)= \Lambda\otimes\unit_n - S_n,
\end{equation}
where
\begin{equation}
  \Lambda = \begin{pmatrix} \lambda & 0\\ 0 & \unit_{k-m}
    \end{pmatrix},
\end{equation}
and
\begin{equation}
  S_n = a_0\otimes \unit_n + \sumr a_i\otimes X_i^{(n)}
\end{equation}
for suitable matrices $a_0, a_1, \ldots, a_r\in M_k(\C)$. Note that,
according to \eqref{3.1}, the $a_i$'s are block matrices of the form
\begin{equation}\label{5.4}
  a_i = \begin{pmatrix} 
0 & a_i^{(1)} & 0 & 0 & \cdots & 0\\
0 & 0 & a_i^{(2)} & 0 & \cdots & 0\\
0 & 0 & 0 & a_i^{(3)} & \cdots & 0\\
\vdots & \vdots & \vdots & \ddots & \ddots & \vdots\\
0 & 0 & 0 & \cdots & 0 & a_i^{(d-1)}\\
a_i^{(d)} & 0 & 0 & \cdots & 0 & 0
\end{pmatrix},
\end{equation}
where $a_i^{(j)}\in M_{m_j, m_{j+1}}(\C)$, $j=1, \ldots, d-1$, and
$a_i^{(d)}\in M_{m_d, m_{1}}(\C)$. As in the previous section, let
\begin{equation}
  H_n(\lambda)= (\id_k\otimes\tr_n)[(\Lambda\otimes\unit_n-S_n)^{-1}],
\end{equation}
and let
\begin{equation}
  G_n(\lambda)= \E\{H_n(\lambda)\}.
\end{equation}
Next, let $x_1, \ldots, x_r$ be a semicircular system in a
$C\cc$-probability space $(\CA,\tau)$ with $\tau$ faithful, and put
\begin{eqnarray}
  q&=&p(x_1, \ldots, x_r),\\
  s&=& a_0\otimes \unit_\CA + \sumr a_i\otimes x_i,
\end{eqnarray}
and
\begin{equation}
  G(\lambda)=(\id_k\otimes \tau)[(\Lambda\otimes\unit_\CA -s)^{-1}].
\end{equation}
Note that, according to Lemma~\ref{inv-formel}, $\Lambda\otimes\unit_\CA -s$ is
invertible. Finally, for every $\mu\in M_k(\C)$, such that
$\mu\otimes\unit_\CA -s$ is invertible, put
\begin{equation}
  \tilde G(\mu)= (\id_k\otimes \tau)[(\mu\otimes\unit_\CA -s)^{-1}].
\end{equation}

\vspace{.2cm}

\begin{lemma}\label{Rtransform of s}
  \begin{itemize}
    \item[(i)] The $\CR$-transform of $s$ w.r.t. amalgamation over
    $M_k(\C)\otimes \unit_\CA$ is given by
    \[
    \CR(z)= a_0 + \sumr a_iza_i, \qquad (z\in M_k(\C)).
    \]
    \item[(ii)] If $\mu\in M_k(\C)$ is invertible, and
    $\|\mu^{-1}\|<\frac{1}{\|s\|}$, then $ \tilde G(\mu)$ is well-defined
    and invertible, and
    \[
    a_0 + \sumr a_i \tilde G(\mu) a_i +  \tilde G(\mu)^{-1}=\mu.
    \]
     \item[(iii)] Let $\mu\in M_k(\C)$ be invertible, and let $R, T\in
     M_k(\C)$ be block diagonal matrices of the form
     \[
     R= {\rm diag}(r_1\unit_{m_1}, r_2\unit_{m_2}, \ldots, r_d\unit_{m_d}),
     \]
     \[
     T= {\rm diag}(t_1\unit_{m_1}, t_2\unit_{m_2}, \ldots, t_d\unit_{m_d}),
     \]
     where $r_1, \ldots, r_d, t_1, \ldots, t_d\in \C\setminus\{0\}$ satisfy
     \[
     r_1t_2 = r_2t_3 = \cdots = r_{d-1}t_d =r_dt_1=1.
     \]
     If $\|(R\mu T)^{-1}\|<\frac{1}{\|s\|}$, then  $ \tilde G(\mu)$ is
     well-defined and invertible, and
     \[
     a_0 + \sumr a_i \tilde G(\mu) a_i +  \tilde G(\mu)^{-1}=\mu.
     \]
   \end{itemize}
\end{lemma}

\proof (i) is essentially due to Lehner [Le]. One just have to exchange $a_i^*$ with $a_i$ in the proof of [Le, Prop.~4.1].

In order to prove (ii), note that if $\|\mu^{-1}\|<\frac{1}{\|s\|}$, then
\[
\mu\otimes\unit_\CA -s = (\mu\otimes\unit_\CA)(\unit_k\otimes\unit_\CA -
(\mu^{-1}\otimes\unit_\CA)  s),
\]
where $\|(\mu^{-1}\otimes\unit_\CA)  s\|<1$. Hence, $\mu\otimes\unit_\CA
-s$ is invertible, and $\tilde G(\mu)$ is well-defined. If in addition
$\|\mu^{-1}\|<\frac{1}{2\|s\|}$, then we get from Neumann's series that
\begin{eqnarray}
  \|\tilde G(\mu)\| &\leq & \|(\mu\otimes\unit_\CA -s)^{-1}\|\nonumber\\
   &\leq & \|\mu^{-1}\|\|(\unit_k\otimes\unit_\CA -
(\mu^{-1}\otimes\unit_\CA)  s)^{-1}\|\nonumber\\
&\leq & 2  \|\mu^{-1}\|.\label{5.13}
\end{eqnarray}
Now, $\tilde G(\mu)$ is the Cauchy transform of $s$ w.r.t. amalgamation
over $M_k(\C)\otimes\unit_\CA$ (cf. \cite{V4}, \cite{Le}). Hence, the maps
\[
z\mapsto \CR(z)+ z^{-1}
\]
and
\[
\mu\mapsto \tilde G(\mu)
\]
are inverses of each other, when $z$ and $\mu$ are invertible, and $\|z\|$
and$ \|\mu^{-1}\|$ are sufficiently small. Thus, according to (i) and
\eqref{5.13}, there is a $\delta \in \Big(0,\frac{1}{2\|s\|} \Big)$, sucht
that when $\mu\in M_k(\C)$ is invertible with $\|\mu^{-1}\|<\delta$, then
\[
 \left\{
\begin{array}{l}
\tilde G(\mu) \; {\rm is} \;{\rm invertible,}\; {\rm and}\\
a_0 + \sumr a_i\tilde G(\mu)a_i + \tilde G(\mu)^{-1} = \mu.
\end{array}
\right .
\]
This statement is equivalent to the identity
\begin{equation}\label{5.14}
\tilde G(\mu)\Big(\mu - a_0 - \sumr a_i\tilde G(\mu)a_i\Big)= \unit_k.
\end{equation}
It is easily seen that
\[
\CU = \Big\{\mu\in GL_k(\C)\,|\, \|\mu^{-1}\|<\frac{1}{\|s\|}\Big\}
\]
is an open, connected set in $M_k(\C)$. Hence, by uniqueness of analytic
continuation, \eqref{5.14} holds for all $\mu\in\CU$. Therefore, for every
$\mu\in\CU$, $\tilde G(\mu)$ is invertible with inverse
\[
\tilde G(\mu)^{-1} = \mu - a_0 - \sumr a_i\tilde G(\mu)a_i.
\]
This proves (ii). 

Finally, to prove (iii), observe that by \eqref{5.4} and \eqref{5.13},
\begin{equation}\label{5.15}
  Ra_iT=a_i, \qquad (i=0, 1, \ldots, r).
\end{equation}
If $\|(R\mu T)^{-1}\|<\frac{1}{\|s\|}$, then we get from (ii) that $\tilde
G(R\mu T)$ is well-defined and invertible, and
\begin{equation}\label{5.16}
  a_0 + \sumr a_i\tilde G(R\mu T)a_i + \tilde G(R\mu T)^{-1} = R\mu T.
\end{equation}
According to \eqref{5.15},
$(R\otimes\unit_\CA)s(T\otimes\unit_\CA)=s$. Hence,
\[
R\mu T \otimes\unit_\CA -s = (R\otimes\unit_\CA)(\mu\otimes\unit_\CA
-s)(T\otimes\unit_\CA).
\]
Then, since $R\mu T \otimes\unit_\CA -s$ is invertible, so is
$\mu\otimes\unit_\CA -s$, and
\[
(\mu\otimes\unit_\CA -s)^{-1} = (T\otimes\unit_\CA)(R\mu T \otimes\unit_\CA
-s)^{-1} (R\otimes\unit_\CA).
\]
It follows that $\tilde G(\mu)$ is well-defined, and
\[
\tilde G(\mu)=T \tilde G(R\mu T) R
\]
is invertible with inverse
\[
\tilde G(\mu)^{-1}= R^{-1} \tilde G(R\mu T)^{-1} T^{-1}.
\]
Taking \eqref{5.15} and \eqref{5.16} into account, we find that
\[
\begin{split}
  a_0 + \sumr a_i&\tilde G(\mu)a_i + \tilde G(\mu)^{-1}  \\
  & = R^{-1}\Big(Ra_0T + \sumr (Ra_iT)\tilde G(R\mu T)(Ra_iT) + \tilde
  G(R\mu T)^{-1}\Big) T^{-1}\\
  & =  R^{-1}\Big(a_0  \sumr a_i\tilde G(R\mu T)a_i + \tilde
  G(R\mu T)^{-1}\Big) T^{-1}\\
  & = R^{-1}R\mu TT^{-1}\\
  & = \mu.
\end{split}
\]
This proves (iii). $\endproof$

\vspace{.2cm}

In the following we let 
\[
\CO = \{\lambda\in M_m(\C)\,|\,\im\lambda\; {\rm is}\; {\rm positive}\;
{\rm definite}\},
\]
and as before, for $\lambda\in\CO$ we put
\[
\Lambda = \begin{pmatrix} \lambda & 0\\ 0 & \unit_{k-m}
    \end{pmatrix}.
\]

\begin{lemma}\label{Lemma5.2} There is a constant $C'$, depending only on $s= a_0 \otimes
  \unit_\CA + \sumr a_i\otimes x_i$, such that:
  \begin{itemize}
    \item[(i)] For all $\lambda\in \CO$,
      \[
      \|(\Lambda\otimes\unit_\CA-s)^{-1}\|\leq C' (1+ \|(\im\lambda)^{-1}\|).
      \]
      Moreover, $G(\lambda)$ is invertible, and
      \[
      a_0 + \sumr a_iG(\lambda)a_i + G(\lambda)^{-1}=\Lambda.
      \]
    \item[(ii)] Let $\lambda\in\CO$, and suppose that $\mu\in M_k(\C)$
    satisfies
    \[
    \|\mu-\Lambda\|<\frac{1}{2C' (1+ \|(\im\lambda)^{-1}\|)}.
    \]
    Then $\mu\otimes \unit_\CA -s$ is invertible, and
    \[
    \|(\mu\otimes \unit_\CA -s)^{-1}\|<2C' (1+ \|(\im\lambda)^{-1}\|).
    \]
    Moreover, $\tilde G(\mu)$ is invertible, and
    \[
     a_0 + \sumr a_i\tilde G(\mu)a_i + \tilde G(\mu)^{-1}=\mu.
    \]
  \end{itemize}
\end{lemma}

\proof (i) With $C_{1,1}'$ and  $C_{2,1}'$ as in Lemma~\ref{estimat of A(lambda)}, put $C'=
\max\{C_{1,1}', C_{2,1}'\}$. Then by Lemma~\ref{estimat of A(lambda)},
\[
\|(\Lambda\otimes\unit_n-s)^{-1}\|\leq C' (1+ \|(\im\lambda)^{-1}\|).
\]
Put
\[
\CO'=\{\lambda\in\CO\,|\, \|\lambda^{-1}\|<\min\{1, \|s\|^{-d}\}\}.
\]
Then $\CO'$ is a non-empty, open subset of $\CO$. At first we will show
that the remaining part of (i) holds for all $\lambda\in\CO'$. Let
$\lambda\in\CO'$, and put
\[
\alpha = \|\lambda^{-1}\|^{\frac 1d}<\min\Big\{1, \frac{1}{\|s\|}\Big\}.
\]
Then with
\[
(r_1, \ldots, r_d)=(\alpha^{d-1}, \alpha^{d-2}, \ldots, \alpha, 1)
\]
and
\[
(t_1, \ldots, t_d)= (1, \alpha^{1-d}, \alpha^{2-d}, \ldots, \alpha^{-1}),
\]
put
\[
R= {\rm diag}(r_1\unit_{m_1}, \ldots, r_d\unit_{m_d})
\]
and
\[
T= {\rm diag}(t_1\unit_{m_1}, \ldots, t_d\unit_{m_d}).
\]      
Then $r_1t_2=r_2t_3=\cdots = r_{d-1}t_d= r_dt_1=1$, and
\[
R\Lambda T= \begin{pmatrix} \alpha^{d-1}\lambda & 0\\ 0 &
  \alpha^{-1}\unit_{k-m}\end{pmatrix},
\]
whence
\[
\|(R\Lambda T)^{-1}\|= \max\{\alpha^{d-1}\|\lambda\|, \alpha\}= \alpha <
\frac{1}{\|s\|}.
\]
Therefore, by Lemma~\ref{Rtransform of s}, $G(\lambda)=\tilde G(\Lambda)$
is invertible and satisfies
\[
a_0 + \sumr a_i G(\lambda) a_i + G(\lambda)^{-1} = \Lambda.
\]
It follows that
\begin{equation}\label{5.17}
G(\lambda)\Big(\Lambda - a_0 - \sumr a_i G(\lambda) a_i\Big)= \unit_k
\end{equation}
holds for all $\lambda \in \CO'$, but then, since $\CO$ is open and
connected, it follows from uniqueness of analytic continuation that
\eqref{5.17} holds for all $\lambda\in\CO$. That is, for every such
$\lambda$, $G(\lambda)$ is invertible and satisfies
\[
      a_0 + \sumr a_iG(\lambda)a_i + G(\lambda)^{-1}=\Lambda.
\]

(ii) Suppose that $\mu\in M_k(\C)$ and
\[
    \|(\mu\otimes \unit_\CA -s)^{-1}\|<2C' (1+ \|(\im\lambda)^{-1}\|).
\]
According to (i), $\|(\Lambda\otimes \unit_\CA- s)^{-1}\|\leq C'(1 +
\|\im\lambda\|^{-1})$. Put $x= \Lambda\otimes \unit_\CA -s$ and  $y=
\mu\otimes \unit_\CA -s$. Then
\[
\|x^{-1}(x-y)\| < \frac12.
\]
Therefore, $y=x(\unit- x^{-1}(x-y))$ is invertible, and by Neumann's series,
\[
\|y^{-1}\| \leq \Big\|\sum_{n=0}^\infty (x^{-1}(x-y))^n\Big\|\|x^{-1}\|\leq
2\|x^{-1}\|.
\]
Hence,
\[
\|(\mu\otimes \unit_\CA -s)^{-1}\|\leq 2C'(1+ \|(\im\lambda)^{-1}\|).
\]
Put
\[
\CO'' = \bigcup_{\lambda\in\CO}\Big\{\mu\in M_k(\C)\,|\,
  \|\mu-\Lambda\|<\frac{1}{2C'(1+ \|(\im\lambda)^{-1}\|)} \Big\}.
\]
Since $\CO$ is connected, $\CO''$ is a connected, open subset of
$M_k(\C)$. In order to prove (ii), we must show that
\begin{equation}\label{5.18}
  \tilde G(\mu)\Big(\mu -  a_0 - \sumr a_i\tilde G(\mu)a_i\Big)=\unit_k
\end{equation}
holds for all $\mu\in\CO''$. Again, by uniqness of analytic continuation,
it suffices to show that \eqref{5.18} holds for all $\mu$ in a non-empty
open subset of $\CO''$. Choose $\lambda\in\CO$ such that $\|\lambda^{-1}\|<
\min\{1, \|s\|^{-d}\}$, and define block matrices $R, T\in M_k(\C)$ as in
the proof of (i). Then $\|(R\Lambda T)^{-1}\|<\frac{1}{\|s\|}$. Since
$x\mapsto x^{-1}$ is continuous on $GL_k(\C)$, we may choose $\delta\in
\Big(0, \frac{1}{2C'(1+ \|(\im\lambda)^{-1}\|)}\Big)$, such that if
$\|\mu-\Lambda\|<\delta$, then $\mu$ is invertible, and
\[
\|(R\mu T)^{-1}\|<\frac{1}{\|s\|}.
\]
Put
\[
\CO''' =\{\mu\in M_k(\C)\,|\, \|\mu-\Lambda\|<\delta\}.
\]
Since $\delta < \frac{1}{2C'(1+ \|(\im\lambda)^{-1}\|)}$, $\CO'''\subseteq
\CO''$. Moreover, we get from
Lemma~\ref{Rtransform of s} that when $\mu\in\CO'''$, then $\tilde G(\mu)$
is invertible, and
\[
\tilde G(\mu)^{-1} = \mu - a_0 - \sumr a_i \tilde G(\mu) a_i.
\]
That is, \eqref{5.18} holds for all $\mu\in \CO'''$ and therefore for all
$\mu\in\CO''$. $\endproof$ 

\vspace{.2cm}

Let $\lambda\in\CO$, and put $\Lambda= \begin{pmatrix} \lambda & 0 \\ 0 &
  \unit_{k-m}\end{pmatrix}$. 
According to Theorem~\ref{master ineq}, there is a constant $C_1\geq 0$ such that
\begin{equation}\label{estim1}
\Big\|\Big(\Lambda-a_0-\sumr a_iG_n(\lambda)a_i\Big)G_n(\lambda)-\unit_k \Big\|\leq \frac{C_1}{n^2}(1+\|(\im\lambda)^{-1}\|^4).
\end{equation}
Put
\[
B_n(\lambda) = \Lambda-a_0-\sumr a_iG_n(\lambda)a_i, \qquad (\lambda\in\CO).
\]
Then, by Neumann's Lemma and \eqref{estim1}, if
\begin{equation}\label{estim2}
\frac{C_1}{n^2}(1+\|(\im\lambda)^{-1}\|^4) <\frac 12,
\end{equation}
then $B_n(\lambda)G_n(\lambda)$ is invertible with $\|(B_n(\lambda)G_n(\lambda))^{-1}\|\leq 2$. Hence $G_n(\lambda)$ is invertible too with
\begin{equation}\label{estim2a}
\begin{split}
\|G_n(\lambda)^{-1}\| & \leq \|(B_n(\lambda)G_n(\lambda))^{-1}\|\|B_n(\lambda)\|\\
& \leq 2 \|B_n(\lambda)\|\\
& \leq 2\big(\|\lambda\|+1+ \|a_0\| + \sumr \|a_i\|^2\|G_n(\lambda)\|\big).
\end{split}
\end{equation}
Taking Lemma~\ref{normlemma} into account we find that for some constant $C_2\geq 0$,
\begin{equation}\label{5.21}
\|G_n(\lambda)^{-1}\|\leq C_2(\|\lambda\|+1)(1+\|(\im\lambda)^{-1}\|).
\end{equation}
By \eqref{estim1}, if $\lambda\in\CO$ satisfies \eqref{estim2}, then
\begin{equation}\label{estim4}
\begin{split}
\|\Lambda-a_0-\sumr a_iG_n(\lambda)a_i - G_n(\lambda)^{-1}\| & \leq \frac{C_1}{n^2}(1+\|(\im\lambda)^{-1}\|^4)\|G_n(\lambda)^{-1}\|\\
& \leq \frac{C_1}{n^2}(1+\|(\im\lambda)^{-1}\|^4)  C_2(\|\lambda\|+1)(1+\|(\im\lambda)^{-1}\|)\\
& \leq  \frac{C_3}{n^2}(1+\|\lambda\|)(1+\|(\im\lambda)^{-1}\|^5) 
\end{split}
\end{equation}
for some constant $C_3\geq 0$. 
For $\lambda\in\CO$ fulfilling (\ref{estim2}) define $\Lambda_n(\lambda)\in M_k(\C)$ by
\begin{equation}\label{5.23}
\Lambda_n(\lambda) = a_0 + \sumr a_iG_n(\lambda)a_i + G_n(\lambda)^{-1}.
\end{equation}
Note that 
\begin{equation}\label{estim6}
\Lambda - \Lambda_n(\lambda) = B_n(\lambda) - G_n(\lambda)^{-1},
\end{equation}
and therefore, by (\ref{estim4}),
\begin{equation}\label{estim7}
\|\Lambda - \Lambda_n(\lambda)\|\leq \frac{C_3}{n^2}(1+\|\lambda\|)(1+\|(\im\lambda)^{-1}\|^5).
\end{equation}
Let $C'$ be as in Lemma~\ref{Lemma5.2}. Then, if also
\begin{equation}\label{5.24}
  \frac{2C'C_3}{n^2}(1+\|\lambda\|)(1+\|(\im\lambda)^{-1}\|^5)(1+\|(\im\lambda)^{-1}\|)<1,
\end{equation}
then
\begin{equation}\label{5.25}
  \|\Lambda_n(\lambda)-\Lambda\|< \frac{1}{2C'(1+ \|(\im\lambda)^{-1}\|)}.
\end{equation}
Hence, by Lemma~\ref{Lemma5.2}, $\tilde G(\Lambda_n(\lambda))$ is
well-defined and invertible and satisfies
\begin{equation}\label{5.26}
a_0 + \sumr a_i\tilde G(\Lambda_n(\lambda))a_i + \tilde
G(\Lambda_n(\lambda)^{-1}= \Lambda_n(\lambda).
\end{equation}
Put
\[
C_4 = \max\{2C_1, 2C'C_3\}
\]
and
\begin{equation}\label{5.27}
  V_n = \Big\{\lambda\in\CO\,\Big|\, \frac{C_4}{n^2}(1+\|\lambda\|)(1+\|(\im\lambda)^{-1}\|^5)(1+\|(\im\lambda)^{-1}\|)< 1 \Big\}.
\end{equation}
Then for all $\lambda\in V_n$, \eqref{5.21} and \eqref{5.24} hold, and
hence also \eqref{5.25} and \eqref{5.26} hold. Observe that the set
\[
U_n = \Big\{\i t\unit_m\,\Big|\,t>0,\,\frac{C_4}{n^2}(1+t)(1+t^{-5})(1+t^{-1})< 1 \Big\}
\]
is contained in $V_n$, and that the function
\begin{equation}\label{est22}
f: t\mapsto (1+t)(1+t^{-5})(1+t^{-1})
\end{equation}
is strictly convex on $]0,\infty[$ and satisfies
\begin{equation}\label{est23}
\begin{split}
&f(t)\rightarrow \infty, \quad t\rightarrow 0^+,\\
& f(t) \rightarrow \infty, \quad t\rightarrow \infty.
\end{split}
\end{equation}
Therefore the set
\begin{equation}\label{est21}
\CI_n = \Big\{ t>0 \,\Big|\,\frac{C_4}{n^2}(1+t)(1+t^{-5})(1+t^{-1})< 1 \Big\}
\end{equation}
is either empty or an open, bounded interval. In particular, $U_n$ is arc-wise connected.

For $\lambda\in\CO$ put $\eps(\lambda)=\|(\im\lambda)^{-1}\|^{-1}$. Then,
as in [HT, Proof of Proposition~5.6], we find that $\i
\eps(\lambda)\unit_m\in U_n$ for all $\lambda\in V_n$, and that the line
segment connecting $\lambda$ and $\i \eps(\lambda)\unit_m$ is contained in
$V_n$.  Hence, either $V_n=\emptyset$ or $V_n$ is connected.

For $\lambda\in V_n$ we get from \eqref{5.23} and \eqref{5.26} that
\begin{equation}\label{5.28}
  \sumr a_i G_n(\lambda) a_i + G_n(\lambda)^{-1}= \sumr a_i \tilde
  G(\Lambda_n(\lambda)) a_i + \tilde
  G(\Lambda_n(\lambda))^{-1}.
\end{equation}
In the following we will show that \eqref{5.28} implies that $G_n(\lambda)=
\tilde G(\Lambda_n(\lambda))$ for all $\lambda\in V_n$.

\begin{lemma}\label{Lemma5.3}
Let $z,w\in GL_k(\C)$, and suppose that
\begin{equation}\label{5.29}
  \sumr a_i za_i + z^{-1} =  \sumr a_i wa_i + w^{-1}.
\end{equation}
If there exists $T\in GL_k(\C)$, such that
\begin{equation}\label{5.30}
  \sumr \|wa_i\|\|Ta_izT^{-1}\| <1,
\end{equation}
then $z=w$.
\end{lemma}

\proof By \eqref{5.29},
\[
w\Big(\sumr a_i za_i + z^{-1}\Big)z= w\Big(\sumr a_i wa_i + w^{-1}\Big)z,
\]
i.e.
\[
\sumr wa_i (z-w)a_i z = z-w.
\]
Therefore,
\[
\sumr wa_i(z-w)T^{-1}(Ta_izT^{-1})= (z-w)T^{-1},
\]
which implies that
\[
\Big(\sumr \|wa_i\|\|Ta_izT^{-1}\|\Big)\|(z-w)T^{-1}\|\geq \|
(z-w)T^{-1}\|.
\]
Hence, if \eqref{5.30} holds, then $)\|(z-w)T^{-1}\|=0$, and thus
$z=w$. $\endproof$

\vspace{.2cm}

\begin{lemma}\label{Lemma5.4}For all $n\in\N$
  and all $\lambda\in\CO$ with $\|(\im\lambda)^{-1}\|\leq 1$ there exists
  $T\in GL_k(\C)$, depending only on $\|(\im\lambda)^{-1}\|$, and a constant
  $C_5>0$, depending only on $a_0, \ldots, a_r$, such that 
  \[
  \|Ta_i G_n(\lambda)T^{-1}\| \leq C_5 \|(\im\lambda)^{-1}\|^{\frac 1d},
  \qquad (i=0,1,\ldots, r).
  \]
\end{lemma}

\proof With the same notation as in the proof of Lemma~\ref{normlemma},
\begin{eqnarray*}
  \|G_n(\lambda)\|&\leq & \E\{\|A_n(\lambda)^{-1}\|\}\\
  &\leq & \E\{\|C_n\|\}+
  \E\{\|B_n^{(1)}\|\|B_n^{(2)}\|\}\|(\im\lambda)^{-1}\|,
\end{eqnarray*}
where
\[
C_{1,1}:= \sup_{n\in\N}\E\{\|C_n\|\}<\infty,
\]
and
\[
C_{2,1}:= \sup_{n\in\N} \E\{\|B_n^{(1)}\|\|B_n^{(2)}\|\}<\infty.
\]
In the same way we get for $T\in GL_k(\C)$ that
\begin{eqnarray}
  \|Ta_iG_n(\lambda)T^{-1}\|&\leq & \E\{\|Ta_iC_nT^{-1}\|\}+
  \E\{\|Ta_iB_n^{(1)}\|\|B_n^{(2)}T^{-1}\|\}\|(\im\lambda)^{-1}\|\nonumber\\
  &\leq &\E\{\|Ta_iC_nT^{-1}\|\}+ C_{2,1}\|a_i\|\|T\|\|T^{-1}\|\|(\im\lambda)^{-1}\|.\label{5.31}
\end{eqnarray}
By Lemma~\ref{inv-formel}, $C_n$ has at most $\frac 12 d(d-1)$ non-zero
block entries and takes the form
\[
C_n=
\begin{pmatrix}
0 & 0 & 0 & 0 & 0 & \cdots & 0 \\
0 & \ast & \ast & \ast & \ast & \cdots & \ast \\
0& 0 & \ast & \ast & \ast &\cdots & \ast  \\
0 & 0 & 0 & \ast & \ast & \cdots & \ast \\
\vdots & \vdots & \vdots & \vdots & \ddots & \ddots & \vdots   \\
0 & 0 & 0 & 0 & \cdots & \ast & \ast \\
0 & 0 & 0 & 0 & \cdots & 0 & \ast \\
\end{pmatrix}.
\]
Combining this with \eqref{5.4}, we find that $a_i C_n$ is a strictly upper
triangular $d\times d$ block matrix. Let $\beta \geq 1$, and put
\begin{equation}\label{5.32}
T= {\rm diag}(\beta\unit_{m_1},\beta^2\unit_{m_2}, \ldots,
\beta^d\unit_{m_d}).
\end{equation}
Then $Ta_i C_n T^{-1}$ is obtained from $a_i C_n$ by multiplying the $(\mu,
\nu)$'th block entry by $\beta^{\mu-\nu}$. Since $a_i C_n$ is strictly
upper triangular,
\begin{eqnarray*}
\|Ta_i C_n T^{-1}\|&\leq &\sum_{\mu<\nu} \|[Ta_i C_n T^{-1}]_{\mu,\nu}\|\\
&=& \sum_{\mu<\nu} \beta^{\mu-\nu} \|[a_i C_n]_{\mu,\nu}\|\\
&\leq &  \beta^{-1} \sum_{\mu<\nu} \|[a_i C_n]_{\mu,\nu}\|,
\end{eqnarray*}
where $[x]_{\mu,\nu}$ denotes the $(\mu,\nu)$'th block entry of a matrix
$x\in M_k(\C)$. Hence,
\[
\|Ta_i C_n T^{-1}\|\leq \frac{d^2}{\beta}\|a_i C_n\|
\]
Since $\E\{\|C_n\|\}\leq C_{1,1}$, we get from \eqref{5.31} that
\begin{equation}\label{5.33}
\|Ta_i G_n(\lambda)T^{-1}\|\leq \Big(\frac{d^2}{\beta}C_{1,1}+  C_{2,1}\|T\|\|T^{-1}\|\|(\im\lambda)^{-1}\|\Big)\|a_i\|
\end{equation}
Moreover, since $\beta\geq 1$, we have that $\|T\|= \beta^d$ and
$\|T^{-1}\|= \beta^{-1}$. Now, if $\|(\im\lambda)^{-1}\|\leq 1$, put
$\beta= \|(\im\lambda)^{-1}\|^{-\frac 1d}\geq 1$. Then by \eqref{5.33},
\[
\|Ta_i G_n(\lambda)T^{-1}\|\leq  \frac{\|a_i\|}{\beta}(C_{1,1}d^2 + C_{2,1}).
\]
Put $C_5= \Big(\sumr \|a_i\|\Big)(C_{1,1}d^2 + C_{2,1})$. Then
\[
\|Ta_i G_n(\lambda)T^{-1}\|\leq \frac{C_5}{\beta} =
C_5\|(\im\lambda)^{-1}\|^{\frac 1d}. \qquad\qquad \endproof
\]

\begin{lemma}\label{Lemma5.5} There is a positive integer $N$, such that
  for all $n\geq N$,
  \[
  G_n(\lambda)=\tilde G(\Lambda_n(\lambda)), \qquad (\lambda\in V_n).
  \]
\end{lemma}

\proof Let $\lambda\in V_n$, and put $z= G_n(\lambda)$ and $w=\tilde
G(\Lambda_n(\lambda))$. According to \eqref{5.28},
\[
\sumr a_iza_i + z^{-1}= \sumr a_iwa_i + w^{-1}.
\]
Moreover, by Lemma~\ref{Lemma5.2}~(ii) and \eqref{5.25} we have that
\begin{eqnarray*}
  \|w\|&\leq  & \|(\Lambda_n(\lambda)\otimes\unit_\CA-s)^{-1}\|\\
  &\leq & 2C'(1+\|(\im\lambda)^{-1}\|).
\end{eqnarray*}
Thus, if $\|(\im\lambda)^{-1}\|<1$, then $\|w\|<4C'$. Moreover, by
Lemma~\ref{Lemma5.4} there exist a constant $C_5$ and $T\in GL_k(\C)$, such
that
\[
\|Ta_izT^{-1}\|\leq C_5\|(\im\lambda)^{-1}\|^{\frac 1d}.
\]
Hence
\[
\sumr \|wa_i\|\|Ta_izT^{-1}\| < 4C'C_5\Big(\sumr
\|a_i\|\Big)\|(\im\lambda)^{-1}\|^{\frac 1d}.
\]
Put
\[
\eps = \min\big\{\big(4C'C_5\sumr \|a_i\|\big)^{-d}, 1\big\}
\]
and
\begin{equation}\label{5.35}
  V_n' = \{\lambda\in V_n\,|\, \|(\im\lambda)^{-1}\|<\eps\}.
\end{equation}
Then for all $\lambda\in V_n'$,
\[
\sumr \|wa_i\|\|Ta_izT^{-1}\| <1,
\]
and therefore, by Lemma~\ref{Lemma5.4}, $z=w$. That is, for all $\lambda\in
V_n'$,
\begin{equation}\label{5.36}
G_n(\lambda)=\tilde G(\Lambda_n(\lambda)).
\end{equation}
Recall from the proof of Lemma~\ref{Lemma5.2} that
\[
U_n = \Big\{\i t\unit_m\,\Big|\,t>0,\,\frac{C_4}{n^2}(1+t)(1+t^{-5})(1+t^{-1})< 1 \Big\}
\]
is a subset of $V_n$. Hence, if
\begin{equation}\label{5.37}
\frac{C_4}{n^2}\Big(1+\frac{\eps}{2}\Big)\Big(1+\Big(\frac{\eps}{2}\Big)^{-5}\Big)\Big(1+\Big(\frac{\eps}{2}\Big)^{-1}\Big)< 1,
\end{equation}
then $\i\frac{\eps}{2}\unit_m\in V_n'$. Choose $N\in\N$, such that
\[
N^2 >
C_4\Big(1+\frac{\eps}{2}\Big)\Big(1+\Big(\frac{\eps}{2}\Big)^{-5}\Big)\Big(1+\Big(\frac{\eps}{2}\Big)^{-1}\Big).
\]
Then for all $n\geq N$, \eqref{5.37} holds, and hence $V_n'$ is a non-empty
open subset of $V_n$. Since $V_n$ is open, and $\lambda\mapsto
G_n(\lambda)-\tilde G(\Lambda_n(\lambda))$ is analytic, \eqref{5.36} holds
for all $\lambda\in V_n$ when $n\geq N$. $\endproof$

\vspace{.2cm}

\begin{thm}\label{Theorem5.6}
  There exist $N\in\N$ and a constant $C_6>0$, both depending only on $a_0,
  \ldots, a_r$, such that for all $\lambda\in\CO$ and all $n\geq N$,
  \begin{equation}\label{est25}
\|G_n(\lambda)-G(\lambda)\|\leq \frac{C_6}{n^2}(1+\|\lambda\|)(1+\|(\im\lambda)^{-1}\|^7).
\end{equation}
\end{thm}

\proof Let $N$ be as in Lemma~\ref{Lemma5.5}. Then for $n\geq N$ and
$\lambda\in V_n$,
\begin{eqnarray*}
  \|G_n(\lambda)-G(\lambda)\|&=& \|\tilde G(\Lambda_n(\lambda))-G(\lambda)\|\\
  & \leq & \|(\Lambda_n(\lambda)\otimes\unit_\CA -s)^{-1}-
  (\Lambda\otimes\unit_\CA -s)^{-1}\|\\
  &=& \|(\Lambda_n(\lambda)\otimes\unit_\CA
  -s)^{-1}(\Lambda-\Lambda_n(\lambda))(\Lambda\otimes\unit_\CA -s)^{-1}\|\\
  &\leq &  \|(\Lambda_n(\lambda)\otimes\unit_\CA
  -s)^{-1}\|\|\Lambda-\Lambda_n(\lambda)\|\| (\Lambda\otimes\unit_\CA
  -s)^{-1}\|.
\end{eqnarray*}
Since \eqref{5.25} holds for all $\lambda\in V_n$, we get from
Lemma~\ref{Lemma5.2} that
\[
 \|G_n(\lambda)-G(\lambda)\|\leq 2(C')^2 (1+
 \|(\im\lambda)^{-1}\|)^2\|\Lambda-\Lambda_n(\lambda)\|.
\]
Hence, by \eqref{estim7}, for all $\lambda\in V_n$,
\begin{eqnarray*}
  \|G_n(\lambda)-G(\lambda)\|&\leq &
  \frac{2(C')^2C_3}{n^2}(1+\|\lambda\|)(1+\|(\im\lambda)^{-1}\|)^2(1+\|(\im\lambda)^{-1}\|^5)\\
  &\leq & \frac{C_6^{(1)}}{n^2}(1+\|\lambda\|)(1+\|(\im\lambda)^{-1}\|^7)
\end{eqnarray*}
for some constant $C_6^{(1)}>0$. Next, if $\lambda\in \CO\setminus V_n$,
i.e.
\[
\frac{C_4}{n^2}(1+\|\lambda\|)(1 + \|(\im\lambda)^{-1}\|^5)(1+ \|(\im\lambda)^{-1}\|)\geq 1,
\]
then 
\begin{equation}\label{5.39}
 \|G_n(\lambda)-G(\lambda)\| \leq \frac{C_4}{n^2}(1+\|\lambda\|)(1 +
 \|(\im\lambda)^{-1}\|^5)(1+ \|(\im\lambda)^{-1}\|) (\|G_n(\lambda)\|+\|G(\lambda)\|).
\end{equation}
Put
$C'' = \max\{C_{1,1}, C_{2,1}, C_{1,1}', C_{2,1}'\}$, where $C_{1,1},
C_{2,1}, C_{1,1}', C_{2,1}'$ refer to the constants from
Lemma~\ref{normlemma} and Lemma~\ref{estimat of A(lambda)}. Then
\[
\|G_n(\lambda)\|+\|G(\lambda)\|\leq 2C''(1+\|(\im\lambda)^{-1}\|),
\]
and hence, by \eqref{5.39}, there is a constant $C_6^{(2)}> 0$, such
that for all $\lambda\in \CO\setminus V_n$,
\[
\|G_n(\lambda)-G(\lambda)\| \leq
\frac{C_6^{(2)}}{n^2}(1+\|\lambda\|)(1+\|(\im\lambda)^{-1}\|^7).
\]
Thus, with $C_6= \max\{C_6^{(1)}, C_6^{(2)}\}$, \eqref{est25} holds. $\endproof$

\section{The spectrum of $Q_n$.}\label{sec4b}

As in the previous we consider a fixed polynomial  $p\in
(M_m(\C)\otimes\C\<X_1, \ldots,X_r\>)_{sa}$ and define $Q_n$ and $q$ by
\eqref{defQ_n} and \eqref{defq}, respectively. For $\lambda\in\C$ with $\im\lambda>0$ put
\begin{equation}
g(\lambda)=(\tr_m\otimes\tau)[(\lambda\unit_m\otimes\unit_\CA-q)^{-1}],
\end{equation}
\begin{equation}
g_n(\lambda)=\E\{(\tr_m\otimes\tr_n)[(\lambda\unit_m\otimes\unit_n-Q_n)^{-1}]\}.
\end{equation}

By application of Proposition~\ref{inv-formel}, we find that with $E=\unit_m\oplus 0_{k-m}\in M_k(\C)$,
\[
g(\lambda) = \frac km \,\tr_k(EG(\lambda\unit_m)E),
\]
and
\[
g_n(\lambda) = \frac km \,\tr_k(EG_n(\lambda\unit_m)E).
\]

Hence for every $n\geq N',$
\begin{equation}\label{HST6.1}
|g_n(\lambda)-g(\lambda)|\leq \frac km \frac{C_7}{n^2}(1+\|\lambda\|)(1+\|(\im\lambda)^{-1}\|^7).
\end{equation}

\vspace{.2cm}

\begin{thm}\label{spec2} For every $\phi\in C_c^\infty(\R,\R)$,
\[
\E\{(\tr_m\otimes\tr_n)\phi(Q_n)\} = (\tr_m\otimes\tau)\phi(q)+O\Big(\frac{1}{n^2}\Big).
\]
\end{thm}

\vspace{.2cm}

\proof  This follows from (\ref{HST6.1}) by minor modifications of [HT,
Proof of Theorem~6.2]. $\endproof$

We are now going to prove

\begin{thm}\label{spec1} Let $\phi\in C^\infty(\R,\R)$ such that $\phi$
  is constant outside a compact subset of $\R$, and suppose
\[
\supp(\phi)\cap \sigma(q) = \emptyset.
\]

Then
\begin{eqnarray}
\E\{(\tr_m\otimes\tr_n)\phi(Q_n)\}&=&O(n^{-2}),\label{mean}\\ 
\V\{(\tr_m\otimes\tr_n)\phi(Q_n)\}&=&O(n^{-4}),\label{variance}
\end{eqnarray}
and 
\begin{equation}\label{probability}
P(|(\tr_m\otimes\tr_n)\phi(Q_n)|<n^{-\frac 43}, \;{\rm eventually}\; {\rm as}\;  n\rightarrow \infty) = 1.
\end{equation}
\end{thm}

\vspace{.5cm}

In the proof of this theorem we shall need:

\begin{prop}\label{variance estimate} Let $m,r\in \N$, let $p\in
  M_m(\C)\otimes \C\<X_1,\ldots, X_r\>$ with $p=p\cc$, and for each $n\in
  \N$, let $X_1^{(n)}, \ldots, X_r^{(n)}$ be stocastically independent random matrices from $\SGRM(n, \frac 1n)$. Then for every compactly supported $C^1$-function $\phi: \R\rightarrow \C$ there is a constant $C\geq 0$ such that
\begin{equation}
\V\{(\tr_m\otimes\tr_n)\phi(Q_n)\}\leq \frac{C}{n^2}\E\Big\{(\tr_m\otimes\tr_n)[|\phi'|^2(Q_n)]\cdot \Big(1+\sum_{j=1}^r\|X_j^{(n)}\|^{d-1}\Big)^2\Big\},
\end{equation}

where
\[
Q_n= p(X_1^{(n)}, \ldots, X_r^{(n)}),
\]

and $d= {\rm deg}(p)$.
\end{prop}

\proof Choose monomials $m_j\in\C\<X_1, \ldots, X_r\>$, $j=1,\ldots, N$, and choose $\alpha_j\in M_m(\C)$, $ j=0,\ldots, N$, such that
\begin{equation}\label{p}
p= \alpha_0\otimes 1 + \sum_{j=1}^N\alpha_j\otimes m_j.
\end{equation}

Consider a fixed $n\in\N$, and let $\CE_{r,n}$ be the real vector space
$(M_n(\C)_{sa})^r$ equipped with the Euclidean norm $\|\cdot\|_e$ (cf. [HT,
Section~3]). Then define $f:\CE_{r,n}\rightarrow\C$ by 
\[
f(v_1, \ldots, v_r)= (\tr_m\otimes\tr_n)[\phi(p(v_1, \ldots, v_r))], \qquad (v_1, \ldots, v_r\in M_n(\C)_{sa}).
\]

According to [HT, (4.4)],
\begin{equation}\label{4.4}
\V\{(\tr_m\otimes\tr_n)[\phi(Q_n)]\}\leq \frac 1n\,\E\{\|(\grad f)(X_1^{(n)},\ldots, X_r^{(n)})\|_e^2\}.
\end{equation}

Now, let $v=(v_1, \ldots, v_r)\in \CE_{r,n}$, and let $w=(w_1, \ldots, w_r)\in \CE_{r,n}$ with $\|w\|_e=1$. As in [HT, Proof of Proposition~4.7] we find that
\begin{equation}\label{4.7}
\Big|\diff f(v+tw)\Big|^2\leq \frac{1}{m^2 n^2}\|\phi'(p(v))\|_{2, \Tr_m\otimes\Tr_n}^2\Big\|\diff p(v+tw)\Big\|_{2, \Tr_m\otimes\Tr_n}^2, 
\end{equation}
 
where
\[
\|\phi'(p(v))\|_{2, \Tr_m\otimes\Tr_n}^2= mn\cdot (\tr_m\otimes\tr_n)[|\phi'|^2(p(v))].
\]

According to (\ref{p}),
\[
\begin{split}
\Big\|\diff p(v+tw)\Big\|_{2, \Tr_m\otimes\Tr_n} & \leq  \sum_{j=1}^N\Big\|\alpha_j\otimes \diff m_j(v+tw)\Big\|_{2, \Tr_m\otimes\Tr_n}\\
& = \sum_{j=1}^N\|\alpha_j\|_{2,\Tr_m}\Big\|\diff m_j(v+tw)\Big\|_{2, \Tr_n}.
\end{split}
\]

Making use of the fact that $\|w_j\|_{2, \Tr_n}^2\leq 1$, $j=1,\ldots, r$, we find that with $d_j={\rm deg}(m_j)$,

\[
\Big\|\diff m_j(v+tw)\Big\|_{2, \Tr_n} \leq d_j \cdot \max_{1\leq i\leq r}\|v_i\|^{d_j-1}\leq d\cdot\Big(1+\sum_{i=1}^r\|v_i\|^{d-1}\Big),
\]

and it follows that
\[
\begin{split}
\Big\|\diff p(v+tw)\Big\|_{2, \Tr_m\otimes\Tr_n} & \leq d\cdot \Big(1+\sum_{i=1}^r\|v_i\|^{d-1}\Big)\Big(\sum_{j=1}^N\|\alpha_j\|_{2,\Tr_m}\Big).
\end{split}
\]

Then by insertion into (\ref{4.7}),
\[
\begin{split}
\Big|\diff f(v+tw)\Big|^2\leq \frac{1}{mn}&(\tr_m\otimes\tr_n)[|\phi'|^2(p(v))]\cdot \\
& d^2\cdot \Big(1+\sum_{i=1}^r\|v_i\|^{d-1}\Big)^2\Big(\sum_{j=1}^N\|\alpha_j\|_{2,\Tr_m}\Big)^2.
\end{split}
\]

Since $w$ was arbitrary this implies that with
\[
C=\frac 1m \cdot d^2 \cdot \Big(\sum_{j=1}^N\|\alpha_j\|_{2,\Tr_m}\Big)^2, 
\]

\[
\|(\grad f) (v)\|_e^2\leq \frac{C}{n}(\tr_m\otimes\tr_n)[|\phi'|^2(p(v))]\cdot\Big(1+\sum_{i=1}^r\|v_i\|^{d-1}\Big)^2.
\]

Then by (\ref{4.4}),
\[
\V\{(\tr_m\otimes\tr_n)[\phi(Q_n)]\}\leq \frac{C}{n^2}\E\Big\{(\tr_m\otimes\tr_n)[|\phi'|^2(Q_n)]\cdot \Big(1+\sum_{i=1}^r\|X_i^{(n)}\|^{d-1}\Big)^2\Big\},
\]

as desired. $\endproof$

\vspace{.2cm}

\begin{prop}\label{meanofnorms} There exist universal constants
  $\gamma(k)\geq 0$, $k\in\N$, such that for every $n\in\N$ and for every
  $X_n\in \SGRM(n,\frac 1n)$,
  \begin{equation}\label{UH2}
  \E\{1_{(\|X_n\|>3)}\cdot \|X_n\|^k\}\leq \gamma(k) n\e^{-\frac n2}.
  \end{equation}
\end{prop}

\vspace{.2cm}

\proof Let $k,\, n\in\N$ and $X_n\in \SGRM(n,\frac 1n)$. Define
$F:[0,\infty[\rightarrow [0,1]$ by
\[
F(t)=P(\|X_n\|\leq t),\qquad (t\geq 0).
\]

Recall from [S, Proof of Lemma~6.4] that for all $\eps >0$ one has that
\begin{equation}\label{UH1}
  1-F(2+\eps)\leq 2n\exp\Big(-\frac{n\eps^2}{2}\Big).
\end{equation}

Integrating by parts as in [Fe, Lemma~V.6.1] we get that
\begin{eqnarray*}
 \E\{1_{(\|X_n\|>3)}\cdot \|X_n\|^k\} &=& \int_3^\infty t^k \d F(t)\\
 & = & 3^k(1-F(3)) + k\,\int_3^\infty t^{k-1} (1-F(t))\d t.
\end{eqnarray*}

According to (\ref{UH1}), $1-F(3)\leq 2n\e^{-\frac n2}$ and
\begin{eqnarray*}
  \int_3^\infty t^{k-1} (1-F(t))\d t & = & \int_0^\infty
  (3+t)^{k-1}(1-F(3+t))\d t\\
  & \leq & 2n\, \int_0^\infty (3+t)^{k-1}\exp\Big(-\frac n2 (1+t)^2\Big)\d t\\
  & = & 2n \e^{-\frac n2}\int_0^\infty (3+t)^{k-1}\exp\Big(-\frac n2
  (2t+t^2)\Big)\d t.
\end{eqnarray*}

Hence (\ref{UH2}) holds with
\[
\gamma(k) = 2\cdot 3^k + 2\,\int_0^\infty (3+t)^{k-1}\exp\Big(-\frac n2
  (2t+t^2)\Big)\d t<\infty. \qquad \endproof
\]

\vspace{.5cm}

{\it Proof of Theorem~\ref{spec1}}. (\ref{mean}) follows from
Theorem~\ref{spec2} as in [HT, Proof of Lemma~6.3]. To prove
(\ref{variance}), note that by  Proposition~\ref{variance estimate},
\begin{equation}\label{UH3}
\V\{(\tr_m\otimes\tr_n)\phi(Q_n)\}\leq \frac{C}{n^2} \E\Big\{(\tr_m\otimes\tr_n)[|\phi'|^2(Q_n)]\cdot (r+1)\Big(1+\sum_{i=1}^r\|X_i^{(n)}\|^{2d-2}\Big)\Big\}.
\end{equation}

Let $\Omega_i= \{\omega\in\Omega \,|\, \|X_i^{(n)}(\omega)\|\leq 3\}$,
$i=1,\ldots,r$. Then by Proposition~\ref{meanofnorms},
\[
  \begin{split}
\E\{(\tr_m\otimes\tr_n)[|\phi'|^2(Q_n)]&\|X_i^{(n)}\|^{2d-2}\}\\
& \leq \int_{\Omega_i}3^{2d-2}(\tr_m\otimes\tr_n)[|\phi'|^2(Q_n)]\d P +
\int_{\Omega\setminus\Omega_i}\|\phi'\|_\infty^2\|X_i^{(n)}\|^{2d-2}\d P\\
& \leq 3^{2d-2}\E\{(\tr_m\otimes\tr_n)[|\phi'|^2(Q_n)]\}+
\|\phi'\|_\infty^2\| \gamma(2d-2)n\e^{-\frac n2}.
\end{split}
\]

Applying (\ref{mean}) to $|\phi'|^2$ we get that
\[
\E\{(\tr_m\otimes\tr_n)[|\phi'|^2(Q_n)]\} = O(n^{-2}),
\]

and hence by (\ref{UH3}),
\[
\V\{(\tr_m\otimes\tr_n)\phi(Q_n)\} = O\Big(n^{-4}+\frac 1n \e^{-\frac
  n2}\Big) = O(n^{-4}),
\]

which proves (\ref{variance}).

Finally, (\ref{probability}) follows from (\ref{mean}) and (\ref{variance})
as in [HT, Proof of Lemma~6.3]. $\endproof$

\vspace{.2cm}

As in [HT, proof of Theorem~6.4], (\ref{probability}) implies the following:

\begin{thm}\label{spectrum} 
For any $\eps>0$ and for almost all $\omega\in\Omega$,
\[
\sigma(Q_n(\omega))\subseteq \sigma(q)\,+\,]-\eps,\eps[,
\]
eventually as $n\rightarrow \infty$.

\end{thm}

\section{No projections in $C_{red}^*(\F_r)$ -- a new proof.}\label{sec5}

\begin{thm} \label{noprojections} (\cite{V3}, \cite{PV}). Let $m,r\in\N$, let $x_1,\ldots,x_r$ be a
  semicircular system in $(\CA,\tau)$, and let $e$ be a projection in
  $M_m(C\cc(\unit_\CA,x_1,\ldots,x_r))$. Then \mbox{$(\Tr_m\otimes\tau)e\in
  \N_0$.} In particular, $C_{red}^*(\F_r)$ contains no projections but the
  trivial ones, i.e. $\CP(C_{red}^*(\F_r))=\{0, \unit\}$.
\end{thm}

\proof Choose $p\in M_m(\C)\otimes\C\<X_1,\ldots,X_r\>$, such that $p=p\cc$ and
\[
\|e-p(x_1,\ldots,x_r)\|<\frac 18.
\]

Put $q= p(x_1,\ldots,x_r)$. By [Da, Proposition~2.1] the Hausdorff distance
between the spectra $\sigma(e)$ and $\sigma(q)$ is at most $\|e-q\|$. Hence
$\sigma(q)\subseteq\,]-\frac 18, \frac 18[\,\cup\,]\frac 78, \frac 98[$.

Choose $\phi\in\Ccinf$ such that $0\leq \phi\leq 1$, $\phi|_{]-\frac 14,\frac 14[}=0$ and $\phi|_{]\frac 34,\frac 54[}=1$. $\phi(q)$ is a projection, and 
\[
\|\phi(q)-q\|<\frac 18.
\]

 Consequently, 
\[
\|\phi(q)-e\|<\frac 14 <1,
\]

implying that $\phi(q)$ is equivalent to $e$. In particular, 
\[
(\Tr_m\otimes\tau)e =(\Tr_m\otimes\tau)\phi(q).  
\]

For each $n\in\N$, let $X_1^{(n)}, \ldots, X_r^{(n)}$ be stochastically independent random matrices from $\SGRM(n, \frac 1n)$, and put
\[
Q_n = p(X_1^{(n)}, \ldots, X_r^{(n)}).
\]

We know from Theorem~\ref{spectrum} that there is a $P$-null set
$N\subseteq \Omega$ such that for all $\omega\in\Omega\setminus N$,
\[
\sigma(Q_n(\omega))\subseteq \sigma(q)+\,\big]{\textstyle-\frac 18, \frac 18}\big[\,\subseteq\, \big]{\textstyle-\frac 14, \frac 14}\big[\,\cup\,\big]{\textstyle\frac 34, \frac 54}\big[
\]

holds eventually as $n\rightarrow \infty$.

In particular, when $\omega\in \Omega\setminus N$, there is an
$N(\omega)\in\N$ such that $\phi(Q_n(\omega))$ is a projection for all $n\geq N(\omega)$, and therefore
\begin{equation}\label{noproj1}
(\Tr_m\otimes\Tr_n)\phi(Q_n(\omega))\in\Z.
\end{equation}

Put
\[
Z_n(\omega)=(\tr_m\otimes\tr_n)\phi(Q_n(\omega))-(\tr_m\otimes\tau)\phi(q),
\qquad (\omega\in\Omega).
\]

According to Theorem~\ref{spec2}, $\E\{Z_n\}=
O\Big(\frac{1}{n^2}\Big)$. Moreover, since $\phi'$ vanishes in a
neighbourhood of $\sigma(q)$, we get, as in the
proof of Theorem~\ref{spec1}, that
\[
\V\{Z_m\} = \V\{(\tr_m\otimes\tr_n)\phi(Q_n)\} = O\Big(\frac{1}{n^4}\Big).
\]

As previously noted this implies that 
\[
P(|Z_n|<n^{-\frac 43}, \; {\rm eventually \; as \;} n\rightarrow \infty)=1.
\]

Hence, we may assume that
\begin{equation}\label{noproj2}
(\Tr_m\otimes\Tr_n)\phi(Q_n(\omega))= n(\Tr_m\otimes\tau)\phi(q)+
O(n^{-\frac 13}),
\end{equation}

holds for almost all $\omega\in\Omega\setminus N$ as well.

Now choose $\omega\in\Omega\setminus N$ and $n_0\in\N$ such that
(\ref{noproj1}) and (\ref{noproj2}) hold when $n\geq n_0$. Take $C\geq 0$
such that for all $n\in\N$,
\[
|(\Tr_m\otimes\Tr_n)\phi(Q_n(\omega))- n(\Tr_m\otimes\tau)\phi(q)|\leq
C\cdot n^{-\frac 13}.
\]

Then
\[
\begin{split}
  &{\rm dist}(n(\Tr_m\otimes\tau)\phi(q),\Z)\leq C\cdot n^{-\frac 13},\\
  &{\rm dist}((n+1)(\Tr_m\otimes\tau)\phi(q),\Z)\leq C\cdot (n+1)^{-\frac
  13},
\end{split}
\]

and hence by subtraction,
\[
{\rm dist}((\Tr_m\otimes\tau)\phi(q),\Z)\leq C(n^{-\frac 13}+(n+1)^{-\frac13})
\]

for all $n\geq n_0$. This implies that $(\Tr_m\otimes\tau)\phi(q)\in \Z$.

The last statement of Theorem~\ref{noprojections} follows from this and the fact that $C_{red}^*(\F_r)$ has a unital trace-preserving embedding into $\CA_0 = C\cc(\unit_\CA, x_1\ldots, x_r)$ (cf. [HT, Lemma~8.1]). $\endproof$

\vspace{.2cm}

\begin{remark} The last statement of Theorem~\ref{noprojections} was
  originally proved in \cite{PV} by application of methods from
  K-theory, and also the first statement of  Theorem~\ref{noprojections}
  may be obtained using K-theory. Indeed, in \cite{V3} it was shown that $K_0(\CA_0) = \Z\,
  [\unit_\CA]_0$, where $\CA_0 = C\cc (\unit_\CA, x_1, \ldots, x_r)$.   

\end{remark}
  
\section{Gaps in the spectrum of $q$.}\label{sec6}

As in the previous sections, consider a semicircular system $x_1, \ldots, x_r$ in $(\CA,\tau)$. Take $p\in M_m(\C)\otimes\C\<X_1,\ldots,X_r\>$, such that $p=p\cc$ and put 
\[
q=p(x_1, \ldots, x_r). 
\]

The following is an easy consequence of Theorem~\ref{noprojections}:

\begin{prop}\label{sec6prop1} $\sigma(q)$ is a union of at most $m$  disjoint connected
  sets, each of which is a compact interval or a one-point set.
\end{prop}

\proof $\R\setminus\sigma(q)$ is a union of disjoint open intervals. If $\R\setminus\sigma(q)$ had more than $m+1$ connected components, one could choose $m+1$ non-zero orthogonal projections $e_1, \ldots, e_{m+1}\in M_m(C\cc(\unit_\CA,x_1,\ldots,x_r))$. Since
\[
(\Tr_m\otimes\tau)e_j \in \{1,\ldots, m\}, \qquad (1\leq j\leq m+1),
\]

we would get that
\[
m = (\Tr_m\otimes\tau)(\unit_m\otimes\unit)\geq \sum_{j=1}^{m+1}(\Tr_m\otimes\tau)e_j \geq m+1
\]

-- a contradiction. Consequently, $\R\setminus\sigma(q)$ has at most $m+1$
   connected components, and $\sigma(q)$ is a union of at most $m$ disjoint
   non-empty compact intervals.$\endproof$

\vspace{.2cm}

Now, for each $n\in\N$, let $X_1^{(n)}, \ldots, X_r^{(n)}$ be stochastically independent random matrices from $\SGRM(n, \frac 1n)$, and put
\[
Q_n = p(X_1^{(n)}, \ldots, X_r^{(n)}).
\]

\begin{thm}\label{holes1} Let $\eps_0$ denote the smallest distance between disjoint connected components of $\sigma(q)$, let $\CJ$ be a connected component of $\sigma(q)$, let $0<\eps<\frac 13\eps_0$, and let $\mu_q\in {\rm Prob}(\R)$ denote the distribution of $q$ w.r.t. $\tr_m\otimes\tau$. Then $\mu_q(\CJ)=\frac km$ for some $k\in\{1, \ldots, m\}$, and for almost all $\omega\in \Omega$, the number of eigenvalues of $Q_n(\omega)$ in $\CJ + ]-\eps, \eps[$ is $k\cdot n$, eventually as $n\rightarrow \infty$.
\end{thm}

\proof Take $\phi\in\Ccinf$ such that $0\leq \phi\leq 1$, $\phi|_{\CJ + ]-\eps, \eps[}=1$, and $\phi|_{\R\setminus (\CJ + ]-2\eps, 2\eps[)}=0$. Then, $\phi(q)\in M_m(C\cc(\unit_\CA,x_1,\ldots,x_r))$, and hence, by Theorem~\ref{noprojections},
\[
\mu_q(\CJ) = (\tr_m\otimes\tau)1_\CJ(q) = (\tr_m\otimes\tau)\phi(q) = \frac km
\]

for some  $k\in\{1, \ldots, m\}$.

As in the proof of Theorem~\ref{noprojections} there is a $P$-null set $N\subset\Omega$ such that for all $\omega\in \Omega\setminus N$, 
\[
\sigma(Q_n(\omega))\subseteq \sigma(q)\, + \, ]-\eps, \eps[,
\]

eventually as $n\rightarrow \infty$, and
\begin{equation}\label{eq1}
(\tr_m\otimes\tr_n)\phi(Q_n(\omega))=\frac km + O(n^{-\frac43}).
\end{equation}

In particular, for all $\omega\in \Omega\setminus N$ there exists $N(\omega)\in \N$ such that $\phi(Q_n(\omega))$ is a projection for all $n\geq N(\omega)$. 

For $\omega\in \Omega\setminus N$ and $n\geq N(\omega)$ take $k_n(\omega)\in \{0, \ldots, m\cdot n\}$ such that 
\[
(\tr_m\otimes\tr_n)\phi(Q_n(\omega))=\frac{k_n(\omega)}{m\cdot n}.
\]

Note that $k_n(\omega)$ is the number of eigenvalues of $Q_n(\omega)$ in $\CJ+]-\eps, \eps[$. (\ref{eq1}) implies that 
\[
k_n(\omega)= k\cdot n + O(n^{-\frac13}),
\]

and hence $k_n(\omega) = k\cdot n$ for $n$ sufficiently big. $\endproof$

\section{The real and symplectic cases.}\label{sec7}

In the sections 9 and 10 we will generalize the results of section 4-6 to Gaussian random matrices with real or sympletic entries. The case of polynomials of degree 1 was treated by the second named author in \cite{S}.

The sympletic numbers $\mathbb{H}$ can be expressed as
\[
\mathbb{H}=\mathbb{R}+j\mathbb{R}+k\mathbb{R}+l\mathbb{R},
\]
where $j^2=k^2=l^2=1$ and
\[
jk=-kj=l,\quad kl=-lk=j,\quad lj=-jl=k.
\]
$\mathbb{H}$ can be realized as a subring of $M_2(\mathbb{C})$ with unit $\unitH$ by putting
\[
j=\j,\quad%
k=\k, \quad\mathrm{and}\quad
l=\l.
\]
By this realization of $\mathbb{H}$, the complectification of $\mathbb{H}$
becomes $\mathbb{H^C=H+}i\mathbb{H}=M_2(\mathbb{C})$.

Following the notation of \cite{S} we will consider the following random matrix ensembles:
\begin{enumerate}
\renewcommand{\labelenumi}{(i)}\item
$\GRMR(n,\sigma^2)$ is the set of random matrices $Y$: $\Omega\rightarrow M_n(\mathbb{R})$ fulfilling that the entries of $Y$, $Y_{uv}$, $1\le u,v\le n$, constitute a set of $n^2$ i.i.d. random variables with distribution $N(0,\sigma^2)$.

\renewcommand{\labelenumi}{(ii)}\item
$\GRM^{\mathbb{H}}(n,\sigma^2)$ is the set of random matrices $Y\rightarrow M_n(\mathbb{H})$ of the form
\[
Y=1\otimes Y^{(1)}+j\otimes Y^{(2)}+k\otimes Y^{(3)}+l\otimes Y^{(4)}
\]\noindent
where $Y^{(1)},Y^{(2)},Y^{(3)},Y^{(4)}$ are stocastically independent random matrices from $\GRMR(n,\frac{\sigma^2}{4})$.

\renewcommand{\labelenumi}{(iii)}\item
The ensemble $\GOE(n,\sigma^2)$ (resp. $\GOES(n,\sigma^2)$) from \cite{S} can be described as the set of selfadjoint random matrices, which have the same distribution as
\[
\frac{1}{\sqrt{2}}(Y+Y^*),\quad \textrm{ (resp. }\frac{1}{i\sqrt{2}}(Y-Y^*))
\]\noindent
where $Y\in \GRMR(n,\sigma^2)$.

\renewcommand{\labelenumi}{(iv)}\item
The ensemble $\GSE(n,\sigma^2)$ (resp. $\GSES(n,\sigma^2)$) from \cite{S}
can be described as the set of selfadjoint random matrices  having the same distribution as
\[
\frac{1}{\sqrt{2}}(Y+Y^*),\textrm{ (resp. }\frac{1}{i\sqrt{2}}(Y-Y^*))
\]
where $Y\in \GRM^{\mathbb{H}}(n,\sigma^2)$.
\end{enumerate}
We shall prove the formulas (1.2), (1.4), (1.5), (1.6) for selfadjoint polynomials of arbitrary degree in $r+s$ stocastically independent selfadjoint random matrices
\[
X_1^{(n)},\ldots,X_{r+s}^{(n)},\quad (r,s\ge 0,\quad r+s\ge 1),
\]
where in the real case
\begin{equation}
X_1^{(n)},\ldots, X_r^{(n)}\in \GOE(n,{\textstyle \frac{1}{n}}),\quad
  X_{r+1}^{(n)},\ldots, X_{r+s}^{(n)}\in \GOES(n,{\textstyle \frac{1}{n}}),
\end{equation}
and in the symplectic case
\begin{equation}
X_1^{(n)},\ldots, X_r^{(n)}\in \GSE(n,,{\textstyle\frac{1}{n}}),\quad
X_{r+1}^{(n)},\ldots, X_{r+s}^{(n)}\in \GSES(n,,{\textstyle \frac{1}{n}}).
\end{equation}
The symplectic case of (1.2), (1.4), (1.5) and (1.6) can easily be reduced
to the real case by use of the methods from [S, Section~7]. Therefore, in
the following we will only consider the real case.

Let $r,s\in\mathbb{N}_0$ with $r+s\ge 1$, and for each $n\in\mathbb{N}$, let $X_1^{(n)},\ldots,X_{r+s}^{(n)}$ be independent random matrices such that $X_1^{(n)},\ldots,X_r^{(n)}\in \GOE(n,\frac{1}{n})$ and $X_{r+1}^{(n)},\ldots,X_{r+s}^{(n)}\in GOE^*(n,\frac{1}{n})$. As in the previous sections, we let $p\in(M_m(\mathbb{C})\otimes C\langle X_1,\ldots,X_{r+s}\rangle)_{sa}$ and define random matrices $(Q_n)_{n=1}^\infty$ by
\begin{equation}
Q_n(\omega)=p(X_1^{(n)}(\omega),\ldots,X_{r+s}^{(n)}(\omega)),\qquad (\omega\in\Omega).
\end{equation}\noindent
With $d=\deg(p)$ we may choose $m_1,\ldots,m_{d+1}\in\mathbb{N}$ with $m=m_1=m_{d+1}$ and polynomials $u_j\in M_{m_j,m_{j+1}}(\mathbb{C})\otimes\mathbb{C}\langle X_1,\ldots,X_{r+s}\rangle$ of first degree, $j=1,\ldots,d$, such that $p=u_1u_2\cdots u_d$. For each $n\in\mathbb{N}$ define random matrices $u_j^{(n)},j=1,\ldots,d,$ by
\[
u_j^{(n)}(\omega)=u_j(X_1^{(n)}(\omega),\ldots,X_{r+s}^{(n)}(\omega)),\qquad (\omega\in\Omega).
\]\noindent
Since $Q_n(\omega)$ is self-adjoint, $\lambda\otimes \unit_n-Q_n(\omega)$ is invertible for every $\lambda\in M_m(\mathbb{C})$ with Im$\lambda > 0$. Then, according to Proposition 2.3, the random matrix
\begin{equation}
A_n(\lambda)=\left(\begin{array}{cccccc}
\lambda\otimes \unit_n & -u_1^{(n)} & 0 & 0 & \cdots & 0\\
0 & \unit_{m_2}\otimes \unit_n & -u_2^{(n)} & 0 & \cdots & 0\\
0 & 0 & \unit_{m_3}\otimes \unit_n & -u_3^{(n)} & \cdots & 0\\
\vdots & \vdots & \vdots & \ddots & \ddots & \vdots\\
0 & 0 & 0 & \cdots & \unit_{m_{d-1}}\otimes \unit_n & -u_{d-1}^{(n)}\\
-u_d^{(n)} & 0 & 0 & \cdots & 0 & \unit_{m_d}\otimes \unit_n
\end{array}\right)
\end{equation}\noindent
is (point-wise) invertible in $M_k(\mathbb{C})$, where
$k=\sum_{i=1}^{d}m_i$.

Choose $a_0,\ldots,a_{r+s}\in M_k(\mathbb{C})$ taking the form
\[
a_i=\begin{pmatrix}
0 & a_{i_1} & 0 & \cdots & 0\\
0 & 0 & a_{_i2} & \cdots & 0\\
\vdots & \vdots & \ddots & \ddots & \vdots\\
0 & 0 & \cdots & 0 & a_{i_{ d-1}}\\
a_{i_d} & 0 & \cdots & 0 & 0
\end{pmatrix},
\]
such that with
\[
S_n=a_0\otimes \unit_n+\sum_{i=1}^{r+s}a_i\otimes X_i^{(n)},
\]
and $\Lambda=\lambda\oplus 1_{k-m}$ we have:
\[
A_n(\lambda)=\Lambda\otimes \unit_n-S_n.
\]
As in section~6 put
\[
\CO=\{\lambda\in M_m(\mathbb{C})\ |\ \mathrm{Im\lambda}\ \textrm{is positive definite}\}
\]
For $\lambda\in\CO$ we put
\begin{equation}\label{eq9-25}
H_n(\lambda)=(\id_k\otimes \tr_n)[(\Lambda\otimes \unit_n-S_n)^{-1}]
\end{equation}
and
\begin{equation}
G_n(\lambda)=\mathbb{E}\{H_n(\lambda)\}.
\end{equation}
By [S, Lemma~6.4], Lemma~3.1 also holds in the real case, possibly with new
constants $C_{1,p}$ and $C_{2,p}$. Hence, $G_n(\lambda)$ is well-defined
and it is easy to check, that $\lambda\mapsto G_n(\lambda)$ is an
analytic map from $\CO$ to $M_k(\mathbb{C})$.

As in [S], we will let $A^{-t}$ denote the transpose of the inverse of an invertible matrix $A$.

\begin{thm}
There is a constant $\widetilde{C}_1\ge 0$, such that for every $n\in\mathbb{N}$ and for all $\lambda\in\CO$,
\begin{equation}\label{eq9-2}
\bigg\|\sum_{j=1}^{r+s}a_jG_n(\lambda)a_iG_n(\lambda)+(a_0-\Lambda)G_n(\lambda)+\mathbf{1}_k+\frac{1}{n}R_n(\lambda)
\bigg\|\le\frac{\widetilde{C}_1}{n^2}(1+\|(\mathrm{Im\lambda})^{-1}\|^{4}),   
\end{equation}
where
\begin{equation}
R_n(\lambda)=\sum_{j=1}^{r+s}\sum_{u,v=1}^k \eps_j a_j
e_{uv}^{(k)}\mathbb{E}\{
(\mathrm{id}_k\otimes\mathrm{tr}_n)[(\Lambda\otimes
\mathbf{1}_n-S_n)^{-t}(e_{uv}^{(k)}a_j\otimes \unit_n)(\Lambda\otimes
\unit_n-S_n)^{-1}]\}, 
\end{equation}
$\eps_j=1$, $1\leq j\leq r$, and $\eps_j=-1$, $r+1\leq j\leq r+s$.
\end{thm}

\proof 
We may assume that
\[
\left.\begin{array}{rl}
\displaystyle X_j^{(n)}=\frac{1}{\sqrt{2}}(Y_j^{(n)}+Y_j^{(n)^*}), & (1\le j\le r),\vspace{2.5mm}\\
\displaystyle X_j^{(n)}=\frac{1}{i\sqrt{2}}(Y_j^{(n)}-Y_j^{(n)^*}), & (r+1\le j\le r+s),
\end{array}\right.
\]\noindent
where $Y_1^{(n)},\ldots,Y_{r+s}^{(n)}$ are $r+s$ stocastically independent random matrices from $\GRMR(n,\frac{1}{n})$. Then
\[
S_n=a_0\otimes\mathbf{1}_n+\sum_{j=1}^{r+s}(b_j\otimes Y_j^{(n)}+c_j\otimes Y_j^{(n)^*}),
\]
where
\begin{equation}
\left.\begin{array}{rl}
b_j=c_j=\frac{1}{\sqrt{2}}a_j, & (1\le j\le r),\vspace{2.5mm}\\
b_j=-c_j=\frac{1}{i\sqrt{2}}a_j, & (r+1\le j\le r+s).
\end{array}\right.
\end{equation}\noindent
Following now the proof of [S, Theorem~2.1] we get that
\begin{equation}\label{eq9-1}
\mathbb{E}\big\{(a_0-\lambda)H_n(\lambda)+\sum_{i=1}^{r+s}(b_jH_n(\lambda)c_jH_n(\lambda)+ c_jH_n(\lambda)b_jH_n(\lambda))+\mathbf{1}_k\big\} = -\frac{1}{n}R_n(\lambda),
\end{equation}
where
\[
\left.\begin{array}{rcl}
R_n(\lambda) & = & \displaystyle\sum_{j=1}^{r+s}\sum_{u,v=1}^{k}b_je_{uv}^{(k)}\mathbb{E}\{ (\mathrm{id}_m\otimes \mathrm{tr}_n)(\lambda\otimes \mathbf{1}_n-S_n)^{-t}(e_{uv}^{(k)}b_j\otimes \mathbf{1}_n)(\lambda\otimes \mathbf{1}_n-S_n)^{-1}\}\\
& + & \displaystyle\sum_{j=1}^{r+s}\sum_{u,v=1}^{k}c_j e_{uv}^{(k)}\mathbb{E}\{ (\mathrm{id}_m\otimes\mathrm{tr}_n)(\lambda\otimes\mathbf{1}_n-S_n)^{-t} (e_{uv}^{(k)}c_j\otimes\mathbf{1}_n)(\lambda\otimes\mathbf{1}_n-S_n)^{-1}\}
\end{array}\right.
\]
Hence by \eqref{eq9-1},
\begin{equation}\label{eq9-16}
\mathbb{E}\{((a_0-\lambda)H_n(\lambda)+\sum_{j=1}^{r+s}a_jH_n(\lambda)a_jH_n(\lambda)+ \mathbf{1}_k\}=-\frac{1}{n}R_n(\lambda)
\end{equation}
where
\begin{equation}\label{eq9-3}
R_n(\lambda)=\sum_{j=1}^{r+s}\sum_{u,v=1}^k\eps_j a_j e_{uv}^{(k)}\mathbb{E}\{(\mathrm{id}_m\otimes \mathrm{tr}_n)(\lambda\otimes\mathbf{1}_n-S_n)^{-t}(e_{uv}^{(k)}a_j\otimes\mathbf{1}_n) (\lambda\otimes\mathbf{1}_n-S_n)^{-1}\}
\end{equation}\noindent
and where $\eps_j=1$, $1\le j\le r$ and $\varepsilon_j=-1$, $r+1\le j\le r+s$. Now, combining the method of proof from [S, proof of Theorem~2.4] with the proof of Theorem~4.3 of this paper, one finds that
\[
\bigg\|\sum_{j=1}^{r+s}a_jG_n(\lambda)a_jG_n(\lambda)+(a_0-\lambda)G_n(\lambda)+\mathbf{1}_k + \frac{1}{n}R_n(\lambda)\bigg\| \le \frac{\widetilde{C}}{n^2}(C_{1,4}+C_{2,4}\| (\mathrm{Im\lambda})^{-1}\|^{4})
\]
for some constant $\widetilde{C}>0$ depending only on
$a_0,\ldots,a_r$. Hence \eqref{eq9-2} holds with
$\widetilde{C}_1=\widetilde{C}\cdot (C_{1,4}+C_{2,4})$. $\endproof$

\begin{cor} There is a constant $\widetilde{C}_2\ge 0$, such that for every $\lambda\in\CO$,
\[
\bigg\|\sum_{j=1}^{r+s}a_jG_n(\lambda)a_jG_n(\lambda)+(a_0-\Lambda)G_n(\lambda)+ \mathbf{1}_k\bigg\| \le \frac{\widetilde{C}_2}{n}(1+\|(\mathrm{Im\lambda})^{-1}\|^4).
\]
\end{cor}

\proof 
By \eqref{eq9-3} and Lemma~3.1 (for the real case)
\begin{equation}\label{eq9-9}
\left.\begin{array}{rcl}
\|R_n(\lambda)\| & \le & \displaystyle k^2\bigg(\sum_{j=1}^{r+s}\|a_j\|^2\bigg)\mathbb{E} \{\|(\Lambda\otimes\mathbf{1}_n-S_n)^{-1}\|^2\}\\
& \le & \displaystyle k^2\bigg(\sum_{j=1}^{r+s}\| a_j\|^2\bigg)(C_{2,1}+C_{2,2}\| (\mathrm{Im\lambda})^{-1}\|^2)\\
& \le & C''(\mathbf{1}+\|(\mathrm{Im\lambda})^{-1}\|^2)
\end{array}\right.
\end{equation}\noindent
for a constant $C''$ depending only on $a_0,\ldots,a_r$. Hence by Theorem~9.1,
\[
\left.\begin{array}{rcl}
\displaystyle\bigg\|\sum_{i=1}^{r+s}a_iG_n(\lambda)a_iG_n(\lambda)+(a_0-\Lambda)G_n(\lambda)+\mathbf{1}_k\| & \!\!\le\!\! & \displaystyle \frac{C''}{n}(\mathbf{1}+\|(\mathrm{Im\lambda})^{-1}\|^2)+\frac{\widetilde{C}_1}{n^2} (1+\|(\mathrm{Im\lambda})^{-1}\|^4)\\
& \!\!\le\!\! & \displaystyle\frac{\widetilde{C}_2}{n}(\mathbf{1}+\|(\mathrm{Im\lambda})^{-1}\|^4)
\end{array}\right.
\]\noindent
for a constant $\widetilde{C}_2\ge 0$. $\endproof$

\vspace{.2cm}

Let $(x_1,\ldots,x_{r+s})$ be a semicircular system in a $C^*$-probability space $(\CA,\tau)$, where $\tau$ is a faithfull state on $\CA$. Put
\begin{equation}
s=a_0\otimes\mathbf{1}_{\CA}+\sum_{j=1}^{r+s}a_j\otimes x_j,
\end{equation}
\begin{equation}\label{eq9-4}
G(\lambda)=(\mathrm{id}_k\otimes\tau)[(\Lambda\otimes\mathbf{1}_\CA -s)^{-1}],\qquad (\lambda\in\CO),
\end{equation}
and put
\begin{equation}
\widetilde{G}(\mu)=(\mathrm{id}_k\otimes\tau)[(\mu\otimes\mathbf{1}_\CA -s)^{-1}]
\end{equation}
for all $\mu\in M_k(\mathbb{C})$ with $\mu\otimes\mathbf{1}_\CA -s$ is invertible.

\begin{thm}
There is an $N\in\mathbb{N}$ and a constant $\widetilde{C}_3$ both
depending only on $a_0,\ldots,a_r$ such that for all $\lambda\in\CO$ 
\[
\|G_n(\lambda)-G(\lambda)\|\le\frac{\widetilde{C}_3}{n} (1+\|\lambda\|)
(1+\|(\mathrm{Im\lambda})^{-1}\|^7),
\]
where $G(\lambda)$ is definded by \eqref{eq9-4}.
\end{thm}

\proof This follows from Corollary~9.2 exactly as Theorem~5.6 followed from
Theorem~4.3. One just has to replace $n^2$ by $n$ in the proofs in
Section~5. $\endproof$

\begin{remark}
From the proof of Theorem~9.3, i.e. from section~5 with $n^2$ replaced by
$n$ (cf. the formulas (5.18) through (5.33) and Lemma~5.5), it follows
that there exist positive constants $C_2,C_3$ and $C_4$ such that when $V_n$ denotes the set
\begin{equation}
V_n= \Big\{\lambda\in\CO\,\Big|\, \frac{C_4}{n^2}(1+\|\lambda\|)(1+\|(\im\lambda)^{-1}\|^5)(1+\|(\im\lambda)^{-1}\|)< 1 \Big\},
\end{equation}
then for all $\lambda\in V_n$, $G_n(\lambda)$ is invertible and the following estimate holds:
\begin{equation}\label{eq9-99}
\|G_n(\lambda)^{-1}\|\le C_2(1+\|\lambda\|)(1+\|(\mathrm{Im\lambda})^{-1}\|).
\end{equation}
Moreover, if one defines $\Lambda_n(\lambda)$ by
\begin{equation}\label{eq9-5}
\Lambda_n(\lambda)=a_0+\sum_{j=1}^{r+s}a_jG_n(\lambda)a_j+G_n(\lambda)^{-1},
\end{equation}
then for $\lambda\in V_n$,
\begin{equation}\label{eq9-11}
\|\Lambda_n(\lambda)-\Lambda\| \le \frac{C_3}{n^2}(1+\|\lambda\|)(1+\|(\mathrm{Im\lambda})^{-1}\|^5),
\end{equation}
and
\begin{equation}
\|\Lambda_n(\lambda)-\Lambda\| \le \frac{1}{2C'(1+\|(\mathrm{Im\lambda})^{-1}\|)},
\end{equation}
where $C'$ is the constant from Lemma~5.2. Finally, $\widetilde{G}(\Lambda_n(\lambda))$ is well-defined, and
\begin{equation}\label{eq9-12}
\widetilde{G}(\Lambda_n(\lambda))=G_n(\lambda),\qquad (\lambda\in V_n).
\end{equation}
\end{remark}

As above, consider a semicircular system $x_1,\ldots,x_{r+s}$ in a
$C^*$-probability space $(\CA, \tau)$, where $\tau$ is a faithful state. It
is no loss of generality to assume that 
\[
\CA=C\cc(x_1,\ldots,x_{r+s}).
\]

Note that the random matrices $X_1^{(n)},\ldots,X_r^{(n)}\in \GOE(n,\frac{1}{n})$ are symmetric, whereas $X_{r+1}^{(n)},\ldots,X_{r+s}^{(n)}$ $\in \GOES(n,\frac{1}{n})$ are skew-symmetric. This is the reason for the following choice of ``transposition'' in $\CA$ and $M_k(\CA)$.

\begin{lemma}
\begin{enumerate}
\renewcommand{\labelenumi}{(\theenumi)}

\item
There is a unique bounded linear map $a\mapsto a^t$ of $\CA$ onto itself such that
\begin{enumerate}
\renewcommand{\labelenumii}{(a)}\item $x_j^t=x_j,\qquad (1\le j\le r)$,
\renewcommand{\labelenumii}{(b)}\item $x_j^t=-x_j,\qquad (r+1\le j\le r+s)$,
\renewcommand{\labelenumii}{(c)}\item $(ab)^t=b^ta^t,\qquad (a,b\in\CA)$. 
\end{enumerate}\noindent
Moreover, $(a^t)^t=a$ and $\|a^t\|=\| a\|$ for all $a\in\CA$.

\item
Define a map $a\rightarrow a^t$ of $M_k(\CA)$ onto itself by
\[
((a_{uv})_{u,v=1}^k)^t=(a_{vu}^t)_{u,v=1}^k
\]\noindent
Then $(ab)^t=b^ta^t$, ($a,b\in M_k(\CA)$). Moreover, $(a^t)^t=a$ and $\| a^t\|=\| a\|$ for all $a\in M_k(\CA)$.
\end{enumerate}
\end{lemma}

\proof
(1) By the proof of [S, lemma 5.2(ii)] the is $\tau$-prepreserving $^*$-automorphism $\psi$ of $\CA =C^*(x_1,\ldots,x_r,\mathbf{1})$, such that
\[
\left.\begin{array}{rl}
\psi(x_j)=x_j, & 1\le j\le r\\
\psi(x_j)=-x_j, & r+1\le j\le r+s
\end{array}\right.
\]\noindent
Moreover, by [S, lemma 5.2(i)], there is a canjugate linear $^*$-isomorphism $\varphi$ of $\CA$ such that $\tau\circ\varphi=\bar{\tau}$ and $\varphi(x_j)=x_j,\ 1\le j\le r+s$. Put now
\[
a^t=\psi\circ\varphi(a^*),\qquad a\in\CA
\]\noindent
Then it is clear, that $a\rightarrow a^t$ satisfies all the conditions of (1). Also $a\rightarrow a^t$ is unique by boundedness and (a), (b), (c).\\
(2) It is elementary to check that the map on $M_k(\CA)$ defined in (2) is involutive and reverses the product. Hence it is a $^*$-isomorphism of $M_k(\CA)$ on the opposite algebra $M_k(\CA)^{op}$, in particular it is an isometri.
\hfill{$\blacksquare$}\\
\\
For an invertible element $a\in M_k(\CA)$ we let $a^{-t}$ denote the operator $(a^{-1})^t=(a^t)^{-1}$. In analogy with (9.8) we can now put
\begin{equation}\label{eq9-39}
R(\lambda)=\sum_{j=1}^{r+s}\varepsilon_j a_j e_{uv}^{(k)}(\mathrm{id}_k\otimes\tau) [(\Lambda\otimes 1_{\CA}-s)^{-t} (e_{uv}^{(k)}a_j\otimes\mathbf{1}_\CA)(\Lambda\otimes\mathbf{1}_\CA -s)^{-1}]
\end{equation}\noindent
Note that by lemma 3.2
\begin{equation}\label{eq9-13}
\|R(\lambda)\|\le\widetilde{C}''(1+\|(\mathrm{Im\lambda})^{-1}\|^2),\qquad \lambda\in\CO
\end{equation}\noindent
for a constant $\widetilde{C}''\ge 0$.

\begin{thm} There is a constant $\widetilde{C}_4\ge 0$ such that for all $\lambda\in\CO$ and all $n\in\mathbb{N}$,
\[
\| R_n(\lambda)-R(\lambda)\| \le \frac{\widetilde{C}_4}{n}(1+\|(\mathrm{Im\lambda})^{-1}\|^{14}).
\]
\end{thm}

\vspace{.2cm}

Before proving Theorem~9.6 we will show how the main result of this section
(Theorem~9.7 below) can be derived from Theorem~9.6 and the previous
results of this section.

As in [S, section~4] we put
\[
L(\lambda)=(\mathrm{id}_m\otimes\tau)[(\lambda\otimes\mathbf{1}_\CA -s)^{-t} (R(\lambda)G(\lambda)^{-1}\otimes\mathbf{1}_\CA) (\lambda\otimes\mathbf{1}_\CA -s)^{-1}].
\]

\begin{thm} There is a constant $\widetilde{C}_5\ge 0$ such that for all $\lambda\in\CO$ and all $n\in\mathbb{N}$,
\[
\Big\| G_n(\lambda)-G(\lambda)+\frac{1}{n}L(\lambda)\Big\| \le \frac{\widetilde{C}_5}{n^2}(1+\|\lambda\|) (1+\|(\mathrm{Im\lambda})^{-1}\|^{17}).
\]
for all $\lambda\in\CO$ and all $n\in\mathbb{N}$.
\end{thm}

\proof The proof follows the proof of [S, Theorem~4.4]. As in Remark~9.4 we put
\[
V_n= \Big\{\lambda\in\CO\,\Big|\,
\frac{C_4}{n^2}(1+\|\lambda\|)(1+\|(\im\lambda)^{-1}\|^5)
(1+\|(\im\lambda)^{-1}\|)< 1 \Big\},  
\]
and we let $\Lambda_n(\lambda)$ be given by \eqref{eq9-5}. Then by
Theorem~9.1 and (9.18), for all $\lambda\in\Lambda_n$, 
\begin{eqnarray}
&
&\Big\|\Lambda_n(\lambda)-\Lambda+\frac{1}{n}R_n(\lambda)
G_n(\lambda)^{-1}\Big\|\nonumber\\  
&= & \Big\|\Big(\sum_{j=1}^{r+s}a_jG_n(\lambda)a_jG_n(\lambda) +
(a_0-\Lambda)G_n(\lambda)+\mathbf{1}_k
+\frac{1}{n}R_n(\lambda)\Big)G_n(\lambda)^{-1}\Big\|\nonumber\\ 
&\le & \frac{\widetilde{C}_5}{n^2}(1+\|\lambda\|)(1+\|(\mathrm{Im\lambda})^{-1}\|^6),\label{eq9-10}
\end{eqnarray}
for some constant $\widetilde{C}_5\ge 0$. By Lemma~5.2,
\begin{equation}\label{eq9-6}
\| G(\lambda)\| \le C'(\mathbf{1}+\|(\mathrm{Im\lambda})^{-1}\|,\qquad (\lambda\in\CO)
\end{equation}
Moreover, by Lemma~5.2, $G(\lambda)$ is invertible and
\begin{equation}
G(\lambda)^{-1}=\Lambda-a_0-\sum_{i=1}^{r+s}a_iG(\lambda)a_i.
\end{equation}
Hence, by \eqref{eq9-6},
\begin{equation}\label{eq9-8}
\| G(\lambda)^{-1}\| \le \widetilde{C}'(1+\|\lambda\|)(\mathbf{1}+\|(\mathrm{Im\lambda})^{-1}\|,\qquad (\lambda\in\CO)
\end{equation}
for some constant $\widetilde{C}'\ge 0$. Applying Theorem~9.3, \eqref{eq9-99} and \eqref{eq9-8}, we have for $\lambda\in V_n$,
\[
\left.\begin{array}{rcl}
\| G_n(\lambda)^{-1}-G(\lambda)^{-1}\| & = &
\| G_n(\lambda)^{-1}(G(\lambda)-G_n(\lambda))G(\lambda)^{-1}\|\\
& \le & \displaystyle\frac{\widetilde{C}_6}{n}(1+\|\lambda\|)^3(\mathbf{1}+\|(\mathrm{Im\lambda})^{-1}\|^9)
\end{array}\right.
\]
for a constant $\widetilde{C}_6\ge 0$. Combining this with Theorem~9.6 and \eqref{eq9-9}, we now have
\begin{eqnarray*}
\| R_n(\lambda)G_n(\lambda)^{-1}-R(\lambda)G(\lambda)^{-1}\|&\leq &
 \| R_n(\lambda)\|\| G_n(\lambda)^{-1}-G(\lambda)^{-1}\| + \|R_n(\lambda)-R(\lambda)\|\|G(\lambda)^{-1}\|\\
&\le & \displaystyle\frac{\widetilde{C}_7}{n}(1+\|\lambda\|)^3(1+\|(\mathrm{Im\lambda})^{-1}\|^{15})
\end{eqnarray*}
for a constant $\widetilde{C}_7\ge 0$. Hence by\eqref{eq9-10}, 
\begin{equation}
\Big\|\Lambda_n(\lambda)-\Lambda+\frac{1}{n}R(\lambda)G(\lambda)^{-1}\Big\|\le \frac{\widetilde{C}_8}{n^2}(1+\|\lambda\|)^3(\mathbf{1}+\|(\mathrm{Im\lambda})^{-1}\|^{15})
\end{equation}
for a constant $\widetilde{C}_8\ge 0$. By lemma 5.2 and (9.21), we have for
$\lambda\in V_n$,
\begin{equation}
\|(\Lambda\otimes\mathbf{1}_\CA -s)^{-1}\|\le
C'(\mathbf{1}+\|(\mathrm{Im\lambda})^{-1}\|),
\end{equation}
\begin{equation}
\|(\Lambda_n(\lambda)\otimes\mathbf{1}_\CA -s)^{-1}\|\le
2C'(1+\|(\mathrm{Im\lambda})^{-1}\|).
\end{equation}
Proceeding as in [S, (4.24) and (4.25)], one now gets,  using
\eqref{eq9-11}, \eqref{eq9-12}, \eqref{eq9-13} and the above estimates,
that
\[
\| G_n(\lambda)-G(\lambda)-\frac{1}{n}L(\lambda)\| \le
\frac{1}{n^2}(1+\|\lambda\|)^3 P_\mathbf{1}(\|\mathrm{Im\lambda}\|^{-1}) 
\]
for a polynomical $P_\mathbf{1}$ of degree 17. Finally, if
$\lambda\in\CO\backslash V_n$, then one obtains exactly as in [S, proof of
Theorem~4.4] that 
\[
\Big\|G_n(\lambda)-G(\lambda)-\frac{1}{n}L(\lambda)\Big\| \le (1+\|\lambda\|)^2 P_2(\|(\mathrm{Im\lambda})^{-1}\|)
\]
for a polynomial $P_2$ of degree 13. Put $P=P_1+P_2$. Then $P$ is of degree
17, and 
\[
\Big\| G_n(\lambda)-G(\lambda)-\frac{1}{n}L(\lambda)\Big\|\le
(1+\|\lambda\|)^3 P(\|(\mathrm{Im\lambda})^{-1}\|)
\]
for all $\lambda\in\CO$. This proves Theorem~9.7. $\endproof$

We now return to the proof of Theorem~9.6. The proof will be devided into a
series of lemmas. The first lemma is a simpel but very useful observation: 

\begin{lemma}
Let $A$ be a unital algebra, and let $x, z\in\mathrm{GL}(A)$ and $y\in A$. Then
\[
\begin{pmatrix}
x & y\\
0 & z
\end{pmatrix}
\]
is invertible in $M_2(A)$ with inverse
\[
\begin{pmatrix}
x^{-1} & -x^{-1}yz^{-1}\\
0 & z^{-1}
\end{pmatrix}.
\]
\end{lemma}

\vspace{.2cm}

Let $\lambda\in\CO$ and $x\in M_k(\mathbb{C})$ and as usual put
\[
\Lambda=\begin{pmatrix}
\lambda & 0\\
0 & \mathbf{1}_{k-m}
\end{pmatrix}.
\]
Moreover, we put
\begin{eqnarray}
\pi_n(\lambda,x) & = & \displaystyle\begin{pmatrix}
\Lambda^t\otimes\mathbf{1}_n-S_n^t & x\otimes\mathbf{1}_n \\
0 & \Lambda\otimes\mathbf{1}_n-S_n
\end{pmatrix},\vspace{2mm}\label{eq9-16a}\\
H_n(\lambda,x) & = & (\mathrm{id}_{2k}\otimes\mathrm{tr}_n)[\pi_n(\lambda,x)^{-1}],\vspace{2mm}\label{eq9-17}\\
G_n(\lambda,x) & = & \mathbb{E}\{H_n(\lambda,x)\}\label{eq9-18},
\end{eqnarray}
and
\begin{eqnarray}
\pi(\lambda,x) & = & \displaystyle\begin{pmatrix}
\Lambda^t\otimes\mathbf{1}_\CA -s^t & x\otimes\mathbf{1}_\CA\\
0 & \Lambda\otimes\mathbf{1}_\CA -s
\end{pmatrix},\vspace{2mm}\\
G(\lambda,x) & = & (\mathrm{id}_{2k}\otimes\tau)[\pi(\lambda,x)^{-1}]\label{eq9-23}
\end{eqnarray}
Finally, we put
\[
\widehat{s}=\begin{pmatrix}
s^t & 0\\
0 & s
\end{pmatrix}
\]\noindent
and
\begin{equation}\label{eq9-31}
\widetilde{G}(\mu)=(\mathrm{id}_{2k}\otimes\tau)((\mu\otimes\mathbf{1}_\CA -\hat{s})^{-1})
\end{equation}
whenever $\mu\otimes\mathbf{1}_\CA -\widehat{s}$ is invertible in $M_{2k}(\CA)$. Note that
\begin{equation}\label{eq9-14}
G(\lambda,x)=\widetilde{G}\begin{pmatrix}
\Lambda^t & x\\
0 & \Lambda
\end{pmatrix}.
\end{equation}\noindent
The idea is now to estimate $\| G_n(\lambda,x)-G(\lambda,x)\|$ by the
methods of section 5 (with $n^2$ replaced by $n$). The estimate we obtain
in Lemma~9.15 below combined with Lemma~9.8 will then complete the proof of
Theorem~9.6. 

\begin{lemma}
\begin{enumerate}
\renewcommand{\labelenumi}{(i)}\item
The $R$-transform of $\widehat{s}$ with respect to amalgamation over $M_{2k}(\mathbb{C})$ is
\[
\widehat{R}(z)=\hat{a}_0+\sum_{i=1}^{r+s}\hat{a}_i z\hat{a}_i,\quad z\in M_{2k}(\mathbb{C})
\]\noindent
where
\[
\widehat{a}_i=\begin{pmatrix}
a_i^t & 0\\
0 & a_i
\end{pmatrix},\qquad (i=0,\ldots,r+s).
\]

\renewcommand{\labelenumi}{(ii)}\item
Let $\mu\in M_{2k}(\mathbb{C})$. If $\mu$ is invertible and $\|\mu^{-1}\| <
\frac{1}{\| S\|}$, then $\widetilde{G}(\mu)$ is well-defined and invertible. Moreover,
\[
\widehat{a}_0+\sum_{i=1}^{r+s}\widehat{a}_0\widetilde{G}(\mu)\widehat{a}_i+\widetilde{G}(\mu)^{-1}=\mu.
\]

\renewcommand{\labelenumi}{(iii)}\item
Let $\mu\in M_{2k}(\mathbb{C})$. If $\mu$ is invertible and if
\[
\left\|\left(\begin{pmatrix}
T^t & 0\\
0 & R
\end{pmatrix}\mu\begin{pmatrix}
R^t & 0\\
0 & T
\end{pmatrix}\right)^{-1}\right\| < \frac{1}{\| s\|},
\]
for some choice of block diagonal matrices $R$ and $T$ of the form
\[
\left.\begin{array}{rcl}
R & = & \mathrm{diag}(r_1\mathbf{1}_{m_{1}}, r_2\mathbf{1}_{m_{2}},\ldots, r_d\mathbf{1}_{m_{k}}), \vspace{2mm}\\
T & = & \mathrm{diag}(t_1\mathbf{1}_{m_{1}}, t_2\mathbf{1}_{m_{2}},\ldots, t_d\mathbf{1}_{m_{k}}),
\end{array}\right.
\]
where $r_1,\ldots,r_d, t_1,\ldots,t_d\in\mathbb{C}\backslash\{ 0\}$ satisfy
\[
r_1t_2=r_2t_3=\ldots=r_{d-1}t_d=r_dt_1=1,
\]\noindent
then $\widetilde{G}(\mu)$ is well-defined and invertible and satisfies
\[
\widehat{a}_0+\sum_{i=1}^{r+s}\widehat{a}_i\widetilde{G}(\mu)\widehat{a}_i+\widetilde{G}(\mu)^{-1}=\mu.
\]
\end{enumerate}
\end{lemma}

\proof Observe, that with $R$ and $T$ as in (iii),
\[
\begin{pmatrix}
T^t & 0\\
0 & R
\end{pmatrix}\widehat{a}_i\begin{pmatrix}
R^t & 0\\
0 & T
\end{pmatrix} = \begin{pmatrix}
(Ra_iT)^t & 0\\
0 & Ra_iT
\end{pmatrix} = \widehat{a}_i
\]
for $i=0,\ldots,r+s$ because $Ra_iT=a_i$ by (5.14). The rest of the proof
of Lemma~9.10 is a straightforward generalization of the proof of
Lemma~5.1. $\endproof$

\vspace{.2cm}

Let $B$ denote the open unitball in $M_k(\mathbb{C})$ i.e.
\[
B=\{ x\in M_k(\mathbb{C})\ |\ \| x\| < 1 \}.
\]\noindent

\begin{lemma} There is a constant $\widetilde{C}$ depending only on $a_0,\ldots,a_{r+s}$, such that:
\begin{enumerate}
\renewcommand{\labelenumi}{(i)}\item
For all $\lambda\in\CO$ and $x\in B$
\[
\|\pi(\lambda,x)^{-1}\|\le\widetilde{C}(\mathbf{1}+\|\mathrm{(Im\lambda)}^{-2}\|).
\]
Moreover, for all such $\lambda$ and $x$, $G(\lambda,x)$ is invertible, and
\[
\widehat{a}_0+\sum_{i=1}^{r}\widehat{a}_iG(\lambda,x)\widehat{a}_i + G(\lambda,x)^{-1}=
\begin{pmatrix}
\Lambda^t & x\\
0 & \Lambda
\end{pmatrix}.
\]

\renewcommand{\labelenumi}{(ii)}\item
Let $(\lambda,x)\in\CO\times B$ and assume that $\mu\in M_{2k}(\mathbb{C})$ satisfies
\[
\left\| \mu -\begin{pmatrix}
\Lambda^t & x\\
0 & \Lambda
\end{pmatrix}\right\| < \frac{1}{2\widehat{C}(1+\|(\mathrm{Im\lambda})^{-2}\|)}.
\]
Then $\mu\otimes\mathbf{1}_\CA -\hat{s}$ is invertible, and
\[
\|(\mu\otimes\mathbf{1}_\CA -\hat{s})^{-1}\| < 2\widehat{C}(1+\|(\mathrm{Im\lambda})^{-2}\|).
\]
Moreover, $\widetilde{G}(\mu)$ is invertible and
\[
\hat{a}_0+\sum_{i=1}^{r+s}\hat{a}_i\widetilde{G}(\mu)\hat{a}_i+\widetilde{G}(\mu)^{-1}=\mu.
\]\noindent
\end{enumerate}
\end{lemma}

\proof Since $a\mapsto a^t$ is an isometry of $M_k(\CA)$ and since
$(a^{-1})^t=(a^t)^{-1}$, when $a$ is invertible, we have that
\[
\|(\Lambda^t\otimes\mathbf{1}_\CA -s^t)^{-1}\| =\|(\Lambda\otimes\mathbf{1}_\CA -s)^{-1}\|.
\]
Hence for $\lambda\in\CO$ and $x\in B$, we get by Lemma~9.8 and Lemma~3.2 that
\[
\left.\begin{array}{rcl}
\|\pi(\lambda,x)^{-1}\| & = & \displaystyle\left\|\begin{pmatrix}
\Lambda^t\otimes\mathbf{1}_\CA -s^t & x\otimes\mathbf{1}_\CA\\
0 & \Lambda\otimes\mathbf{1}_\CA -s
\end{pmatrix}^{-1}\right\|\vspace{2mm}\\
& \le & \|(\Lambda\otimes\mathbf{1}_\CA -s)^{-1}\| +\|(\Lambda\otimes\mathbf{1}_\CA -s)^{-1}\|^2\| x\|\vspace{2mm}\\
& \le & C'_{1,1}+C'_{2,1}\|(\mathrm{Im\lambda})^{-1}\| + C'_{2,1}+C'_{2,2}\|(\mathrm{Im\lambda})^{-1}\|^2\vspace{2mm}\\
& \le & \widehat{C}(\mathbf{1}+\|(\mathrm{Im\lambda}^{-1}\|^2)
\end{array}\right.
\]
for a constant $\widetilde{C}$ depending only on $C'_{i,j}$, $i,j=1,2$. Put
\[
\CO'=\{\lambda\in\CO\ |\ \|\lambda^{-1}\| < \min\{1,(2\| s\|)^{-d}\}\},
\]
and for a fixed $\lambda\in\CO'$ put
\[
\alpha=\|\lambda^{-1}\|^{\frac{1}{d}}<\min\left\{1,\frac{1}{2\| s\|}\right\}.
\]
Next, let
\[
\left.\begin{array}{rcl}
(r_1,\ldots,r_d) & = & (\alpha^{d-1},\alpha^{d-2},\ldots,\alpha,\mathbf{1}),\vspace{2mm}\\
(t_1,\ldots,t_d) & = & (\mathbf{1},\alpha^{1-d},\alpha^{2-d},\ldots,\alpha^{-1}),
\end{array}\right.
\]
and
\[
\left.\begin{array}{rcl}
R & = & \mathrm{diag}(r_1\mathbf{1}_{m_{1}},\ldots,r_d\mathbf{1}_{m_{d}}),\vspace{2mm}\\
T & = & \mathrm{diag}(t_1\mathbf{1}_{m_{1}},\ldots,t_d\mathbf{1}_{m_{d}}).
\end{array}\right.
\]
Then, as in the proof of Lemma~5.2, we get that
\[
\|(R\Lambda T)^{-1}\|=\alpha <\frac{1}{2\| s\|}.
\]
Hence by Lemma~9.8, we have for $x\in B$ that
\[
\left.\begin{array}{rcl}
\left\|\left(\begin{pmatrix}
T^t & 0\\
0 & R
\end{pmatrix}\begin{pmatrix}
\Lambda^t & x\\
0 & \Lambda
\end{pmatrix}\begin{pmatrix}
R^t & 0\\
0 & T
\end{pmatrix}\right)^{-1}\right\| & = & \left\|\begin{pmatrix}
(R\Lambda T)^t & T^t\times T\\
0 & R\Lambda T
\end{pmatrix}^{-1}\right\|\vspace{2mm}\\
& \le & \|(R\Lambda T)^{-1}\| + \|(R\Lambda T)^{-1}\|^2\| T\|^2\| x \|\vspace{2mm}\\
& \le & \displaystyle \frac{1}{2\| s\|} + \frac{\| T\|^2\| x\|}{4\| s\|^2}.
\end{array}\right.
\]
Moreover, $\| T\| =\alpha^{1-d}\le\alpha^{-d}=\|\lambda^{-1}\|^{-1}$. Thus, if $\| x\| < 2\| s\|\ \|\lambda^{-1}\|^2$, then
\[
\left\|\begin{pmatrix}
T^t & 0\\
0 & R
\end{pmatrix}\begin{pmatrix}
\Lambda^t & x\\
0 & \Lambda
\end{pmatrix}\begin{pmatrix}
R^t & 0\\
0 & T
\end{pmatrix}^{-1}\right\| < \frac{1}{\| s\|},
\]
so by \eqref{eq9-14} and Lemma~9.10~(iii),
\begin{equation}\label{eq9-15}
\hat{a}_0+\sum_{i=1}^{r}\hat{a}_iG(\lambda,x)\hat{a}_i+G(\lambda,x)^{-1}=\begin{pmatrix}
\Lambda^t & x\\
0 & \Lambda
\end{pmatrix}.
\end{equation}
Since
\[
\{(\lambda,x)\in\CO'\times B\ |\ \| x\| < 2\| s\|\ \|\lambda^{-1}\|^2\}
\]\noindent
is a non-empty open subset of $\CO\times B$, we can use uniqueness of
analytic continuation as in the proof of Lemma~5.2 and obtain that
$G(\lambda,x)$ is invertible for all $(x,\lambda)\in\CO\times B$ and that
these all satisfy \eqref{eq9-15}. This proves (i). The proof of (ii) is a
straightforward generalization of the proof of Lemma~5.2~(ii). $\endproof$

\vspace{.2cm}

For $\lambda\in\CO$ and $x\in M_k(\mathbb{C})$ we let $\pi_n(\lambda,x),\
H_n(\lambda,x)$ and $G_n(\lambda,x)$ be given by \eqref{eq9-16a},
\eqref{eq9-17} and \eqref{eq9-18}, respectively. Moreover, we put
\begin{equation}\label{eq9-20}
R_n(\lambda,x)=\sum_{i=1}^{r+s}\sum_{u,u=1}^{2k}\hat{a}_i e_{uv}^{(2k)}
\mathbb{E}\{(\mathrm{id}_{2k}\otimes\mathrm{tr}_n)[\pi_n(x,\lambda)^{-t} (e_{uv}^{(2k)}\hat{a}_i\otimes\mathbf{1}_n)\pi_n(x,\lambda)^{-1}]\}.
\end{equation}

\begin{lemma} There is a constant $\widehat{C}_1\ge 0$ only depending on $a_0,\ldots,a_{r+s}$, such that for all $\lambda\in\CO$ and all $x\in B$,
\begin{equation}
\begin{split}
\|\sum_{i=1}^{r+s}\hat{a}_iG_n(\lambda,x)\hat{a}_iG_n(\lambda,x)+\Big(\hat{a}_0-\begin{pmatrix}
\Lambda^t & x\\
0 & \Lambda
\end{pmatrix}\Big)G_n(\lambda,x)+\mathbf{1}_n\Big)+&\frac{1}{n}R_n(\lambda,x)\|\\
& \le\frac{\widehat{C}_1}{n^2}(\mathbf{1}+\|(\mathrm{Im\lambda})^{-1}\|^8)
\end{split}.
\end{equation}
\end{lemma}

\proof This is a fairly straightforward generalization of the proof of
Theorem~9.1. The master equation \eqref{eq9-16} now becomes
\begin{equation}\label{eq9-21}
\mathbb{E}\left\{\sum_{i=1}^{r+s}\hat{a}_iH_n(\lambda,x)\hat{a}_iH_n(\lambda,x) + \left(\hat{a}_0-\begin{pmatrix}
\Lambda^t & x\\
0 & \Lambda
\end{pmatrix}\right) H_n(\lambda,x)+\mathbf{1}_n\right\} = -\frac{1}{n}R_n(\lambda,x).
\end{equation}
Since $\| x\| < 1$, we get from Lemma~9.8 and Lemma~3.1 that
\begin{equation}\label{eq9-22}
\left.\begin{array}{rcl}
\mathbb{E}\{\|\pi_n(\lambda,x)^{-1}\|^p\} & = & \displaystyle\mathbb{E}\left\{\left\|\begin{pmatrix}
\Lambda^t\otimes\mathbf{1}_n-s_n^t & x\otimes\mathbf{1}_n\\
0 & \Lambda\otimes\mathbf{1}_n-s_n
\end{pmatrix}^{-1}\right\|^p\right\}\vspace{2mm}\\
& \le & \mathbb{E}\{(\|\Lambda\otimes\mathbf{1}_n-s_n\|+ \|\Lambda\otimes\mathbf{1}_n-s_n\|^2)^p\}\vspace{2mm}\\
& \le & 2^p\mathbb{E}\{\|\Lambda\otimes\mathbf{1}_n-s_n\|^p + \|\Lambda\otimes\mathbf{1}_n-s_n\|^{2p}\}\vspace{2mm}\\
& \le & 2^p(C_{1,p}+C_{2,p}\|(\mathrm{Im\lambda})^{-1}\|^p+C_{1,2p}+C_{2,2p}\|(\mathrm{Im\lambda})^{-1} \|^{2p})\vspace{2mm}\\
& \le & C'_p(1+\|(\mathrm{Im\lambda})^{-1}\|^{2p}),
\end{array}\right.
\end{equation}
for a constant $C'_p$ which only depends on the constants $C_{i,j},\
i=1,2,\ j=p,2p$. Applying now \eqref{eq9-22} for $p=4$, Lemma~9.11 follows
from \eqref{eq9-21} exactly as in the proof of Theorem~9.1. $\endproof$

\begin{cor}
There is a constant $\widehat{C}_2\ge 0$, such that for all $\lambda\in\CO$ and all $x\in B$,
\[
\left\|\sum_{i=1}^{r+s}\hat{a}_iG_n(\lambda,x) \hat{a}_iG_n(\lambda,x)+\left(\hat{a}_0-\begin{pmatrix}
\Lambda^t & x\\
0 & \Lambda
\end{pmatrix}\right) G_n(\lambda,x)+\mathbf{1}_n\right\| \le\frac{\widehat{C}_2}{n}(1+\|(\mathrm{Im\lambda})^{-1}\|^8).
\]
\end{cor}

\proof By \eqref{eq9-20} and \eqref{eq9-22} (for $p=2$), we have that 
\[
\| R_n(\lambda,x)\| \le \widehat{C}'' (1+\|(\mathrm{Im\lambda})^{-1}\|^4),
\]
for some constant $\widehat{C}''\ge 0$. The corollary now follows
immediately from Lemma~9.12 (cf. the proof of Corollary~9.2). $\endproof$

\vspace{.2cm} 

Note that by \eqref{eq9-17}, \eqref{eq9-18} and \eqref{eq9-22},
\begin{equation}\label{eq9-39a}
\| G_n(\lambda,x)\| \le C'_p(1+\|(\mathrm{Im\lambda})^{-1}\|^2),\qquad ((\lambda,x)\in\CO\times B),
\end{equation}
and by \eqref{eq9-23} and Lemma~9.11~(i),
\begin{equation}\label{eq9-40}
\| G(\lambda,x)\| \le \widehat{C}(1+\|(\mathrm{Im\lambda})^{-1}\|^2),\qquad ((\lambda,x)\in\CO\times B).
\end{equation}
Proceeding now as in (5.18)-(5.23) with $n^2$ replaced by $n$ (see also
Remark~9.4), one finds that after suitable changes of the exponents, that there exists positive constants $C_2,C_3,C_4$ and $C_6$, such that when $\widehat{V}_n$ denotes the set
\begin{equation}\label{eq9-38}
\widehat{V}_n=\Big\{\lambda\in\CO\ |\ \frac{C_4}{n}(1+\|\lambda\|) (1+\|(\mathrm{Im\lambda})^{-1}\|^{10})(1+\|(\mathrm{Im\lambda})^{-1}\|^2)<1\Big\},
\end{equation}
then for $(\lambda,x)\in V_n\times B$, $G_n(\lambda,x)$ is invertible and
\begin{equation}
\| G_n(\lambda)^{-1}\| \le C_2(1+\|\lambda\|) (1+\|(\mathrm{Im\lambda})^{-1}\|^2).
\end{equation}
Moreover, if one defines $\Lambda_n(\lambda,x)$ by
\begin{equation}\label{eq9-27}
\Lambda_n(\lambda,x)=\hat{a}_0+\sum_{i=1}^{r+s}\hat{a}_iG_n(\lambda,x)\hat{a}_i+G_n(\lambda,x)^{-1},
\end{equation}
then for $(\lambda,x)\in \widehat{V}_n\times\CO$,
\begin{equation}\label{eq9-36}
\left\|\Lambda_n(\lambda,x)-\begin{pmatrix}
\Lambda^t & x\\
0 & \Lambda
\end{pmatrix}\right\| \le \frac{C_3}{n^2}(1+\|\lambda\|) (1+\|(\mathrm{Im\lambda})^{-1}\|^{10}),
\end{equation}
and
\begin{equation}
\left\|\Lambda_n(\lambda,x)-\begin{pmatrix}
\Lambda^t & x\\
0 & \Lambda
\end{pmatrix}\right\| \le \frac{1}{2\widehat{C}(1+\|(\mathrm{Im\lambda})^{-1}\|^2)},
\end{equation}
where $\widehat{C}$ is the constant from Lemma~9.10~(ii). Hence by Lemma~9.10~(ii), $\Lambda_n(\lambda,x)\otimes\mathbf{1}_\CA -\hat{s}$ is invertible, and
\begin{equation}\label{eq9-32}
\|(\Lambda_n(\lambda,x)\otimes\mathbf{1}_\CA -\hat{s})^{-1}\| \le 2\widehat{C}(1+\|(\mathrm{Im\lambda})^{-1}\|^2).
\end{equation}
Moreover, $\widetilde{G}(\Lambda_n(\lambda,x))$ is well-defined, invertible and
\begin{equation}\label{eq9-28}
\hat{a}_0+\sum_{i=1}^{r+s}\hat{a}_i\widetilde{G}(\Lambda_n(\lambda,x))\hat{a}_i + \widetilde{G}(\Lambda_n(\lambda,x))^{-1}=\Lambda_n(\lambda,x).
\end{equation}
Recall that for $\lambda\in\CO$,
\[
G_n(\lambda)=\mathbb{E}\{(\mathrm{id}_k\otimes\mathrm{tr}_n)((\Lambda\otimes\mathbf{1}_n-s_n)^{-1})\}. 
\] 

\begin{lemma} There is a constant $C_5\ge 0$, independent of $\lambda$ and $n$, such that when $\lambda\in\CO$ and $\|(\mathrm{Im\lambda})^{-1}\|\le 1$, there exist $R,S\in GL(k,\mathbb{C})$, such that
\begin{equation}\label{eq9-26}
\| RG_n(\lambda)a_iR^{-1}\|\le C_5\|(\mathrm{Im\lambda})^{-1}\|^{\frac{1}{d}},
\end{equation}
and
\begin{equation}\label{eq9-24}
\| Sa_iG_n(\lambda)S^{-1}\|\le C_5\|(\mathrm{Im\lambda})^{-1}\|^{\frac{1}{d}},
\end{equation}
for $n\in\mathbb{N}$ and $1\le i\le r+s$.
\end{lemma}

\proof Lemma~3.1 holds in the real case too (possibly with change of
constants). Therefore Lemma~5.4 also holds in the real case, which proves
\eqref{eq9-24}. Moreover, by the proof of Lemma~5.4, the matrix $S$ given
by
\begin{equation}\label{eq9-29}
S=\mathrm{diag}(\beta\mathbf{1}_{m_{1}},\beta^2\mathbf{1}_{m_{2}},\ldots,\beta^d\mathbf{1}_{m_{d}}),
\end{equation}
where $\beta=\|(\mathrm{Im\lambda})^{-1}\|^{\frac{1}{d}}\ge 1$, satisfies
\eqref{eq9-24}. As in the proof of Lemma~3.1, write
\[
(\lambda\otimes\mathbf{1}_n-S_n)^{-1}=C_n+B_n^{(1)}(\lambda\otimes\mathbf{1}_n-Q_n)^{-1}B_n^{(2)}.
\]
Then by \eqref{5.4} and the positions of the non-zero entries of $C_n$, we get, that $G_n(\lambda)a_i$ is a $d\times d$ block matrix of the form
\[
G_n(\lambda)a_i=\begin{pmatrix}
0 & 0 & \cdots & \cdots & 0\\
* & 0 & * & \cdots & *\\
* & 0 & 0 & \ddots & \vdots\\
\vdots & \vdots & \vdots & \ddots & *\\
* & 0 & 0 & \cdots & 0
\end{pmatrix}.
\]
Let
\begin{equation}\label{eq9-30}
R=\mathrm{diag}(\beta^d\mathbf{1}_{m_{1}},\beta\mathbf{1}_{m_{2}},\beta^2\mathbf{1}_{m_{5}},\ldots, \beta^{d-1}\mathbf{1}_{m_{d}}),
\end{equation}
where as before $\beta =\|(\mathrm{Im\lambda})^{-1}\|^{-\frac{1}{d}}$. Then
the map $G_n(\lambda)a_i\rightarrow RG_n(\lambda)a_iR^{-1}$ multiplies the
upper diagonal entries of $[G_n(\lambda)a_i]_{uv},\ 2\le u<v\le d$, by
$\beta^{u -v}$ and it multiplies the entries $[G_n(\lambda)a_i]_{u 1},\
2\le u\le d$ by $\beta^{u-1-d}$. Thus,
\[
\| [RG_n(\lambda)a_iR^{-1}]_{u v}\| \le\beta^{-1}\| [G_n(\lambda)a_i]_{u v}\|
\]\noindent
for all $u, v\in \{1,\ldots,d\}$. The rest of the proof of \eqref{eq9-26}
is now a simple modification of the proof of Lemma~5.4. $\endproof$

\begin{lemma}
There is a positive integer $N$, such that for all $n\ge N$,
\[
\widetilde{G}(\Lambda_n(\lambda,x))=G_n(\lambda,x),\qquad (\lambda\in\widehat{V}_n,\ x\in B).
\]
\end{lemma}

\proof 
Let $(\lambda,x)\in\widehat{V}_n\times B$, and at first assume that $\|(\mathrm{Im\lambda})^{-1}\|\le 1$. Put
\begin{equation}
z=G_n(\lambda,x)\qquad\textrm{and}\qquad w=\widetilde{G}(\Lambda_n(\lambda,x)).
\end{equation}
By \eqref{eq9-27} and \eqref{eq9-28}, $z$ and $w$ are invertible, and
\begin{equation}
\sum_{i=1}^{r+s}\hat{a}_iz\hat{a}_i+z^{-1}=\sum_{i=1}^{r+s}\hat{a}_iw\hat{a}_i+w^{-1}.
\end{equation}
Put
\begin{equation}
T=\begin{pmatrix}
R^{-t} & 0\\
0 & S
\end{pmatrix},
\end{equation}
where $R,S\in GL(k,\mathbb{C})$ are the matrices from Lemma~9.13 given by
\eqref{eq9-29} and \eqref{eq9-30}. We will show that if $\|(\mathrm{Im\lambda})^{-1}\|$ and $\| x\|$ are sufficiently small, then
\begin{equation}\label{eq9-35}
\sum_{i=1}^{r+s}\| w\hat{a}_i\|\| T\hat{a}_izT^{-1}\| <1,
\end{equation}
and thus, by the proof of Lemma~5.3, it follows, that $z=w$.

By Lemma~9.8 we have for $\lambda\in\CO$ and $x\in B$ that
\begin{equation}\label{eq9-38a}
\left.\begin{array}{rcl}
G_n(\lambda,x) & = & \displaystyle\mathbb{E}\Big\{(\mathrm{id}_{2k}\otimes\mathrm{tr}_n) \left[\begin{pmatrix}
\lambda^t\otimes\mathbf{1}_n-S_n^t & x\otimes \unit_n\\
0 & \lambda\otimes\mathbf{1}_n-S_n
\end{pmatrix}^{-1}\right]\Big\}\vspace{2mm}\\
& = & \begin{pmatrix}
G_n(\lambda)^t & K_n(\lambda,x)\\
0 & G_n(\lambda)
\end{pmatrix}
\end{array}\right.,
\end{equation}
where
\begin{equation}\label{eq9-33}
K_n(\lambda,x)=-\mathbb{E}\{(\mathrm{id}_k\otimes\mathrm{tr}_n) [(\lambda^t\otimes\mathbf{1}_n-s_n^t)^{-1}(x\otimes\mathbf{1}_n)(\lambda\otimes\mathbf{1}_n-s_n)^{-1}]\}.
\end{equation}
With $R,S$ and $T$ as above we have that
\[
T\hat{a}_iG_n(\lambda,x)T^{-1}=\begin{pmatrix}
(RG_n(\lambda)a_iR^{-1})^t & R^{-t}a_i^tK_n(\lambda,x)S^{-1}\\
0 & Sa_iG_n(\lambda)S^{-1}
\end{pmatrix}.
\]
Hence, by Lemma~9.13,
\[
\| T\hat{a}_iG_n(\lambda,x)T^{-1}\|\le C_5\|(\mathrm{Im\lambda})^{-1}\|^\frac{1}{d} + \|R^{-1}\|\|a_i^tK_n(\lambda,x)\|\| S^{-1}\|.
\]
By \eqref{eq9-29} and \eqref{eq9-30}, $\| R^{-1}\| =\| S^{-1}\| =
\frac{1}{\beta}\le 1$. \eqref{eq9-33} and Lemma~3.1 imply that
\[
\| K_n(\lambda,x)\|\le (C_{1,2}+C_{2,2}\|(\mathrm{Im\lambda})^{-1}\|^2)\| x\|
\]
Hence
\begin{equation}\label{eq9-34}
\| T\hat{a}_iG_n(\lambda,x)T^{-1}\|\le C_5\|(\mathrm{Im\lambda})^{-1}\|^\frac{1}{d} + \| x\|\| a_i\| (C_{1,1}+C_{2,2}\|(\mathrm{Im\lambda})^{-1}\|^2).
\end{equation}
Moreover, by \eqref{eq9-31} and \eqref{eq9-32},
\begin{equation}
\|\widetilde{G}(\Lambda_n(\lambda,x))\|\le 2\widehat{C}(1+\|(\mathrm{Im\lambda})^{-1}\|^2.
\end{equation}
By \eqref{eq9-33} and \eqref{eq9-34}, there is a $\delta\in(0,1)$, such that when $\|(\mathrm{Im\lambda})^{-1}\| <\delta$ and $\| x\| <\delta$, then for all $n\in\mathbb{N}$,
\[
\sum_{i=1}^{r+s}\|\widetilde{G}(\Lambda_n(\lambda,x))\hat{a}_i\|\| T\hat{a}_iG_n(\lambda,x)T^{-1}\| < 1.
\]
That is, \eqref{eq9-35} holds, and therefore $z=w$, which shows that
$G_n(\lambda,x)=\widetilde{G}(\Lambda_n(\lambda,x))$ when $(\lambda,x)$
belongs to the set
\[
\CU_n=\{\lambda\in\widehat{V}_n\ |\ \|(\mathrm{Im\lambda})^{-1}\| <\delta\}\times\{ x\in M_k(\mathbb{C})\ |\ \| x\| <\delta\}.
\]\noindent
Exactly as for the sets $V_n$ in section~5, we can prove that
$\widehat{V}_n$ is connected and that there exists $N\in\mathbb{N}$, such
that $\{\lambda\in\widehat{V}_n\ |\ \|(\mathrm{Im\lambda})^{-1}\|
<\delta\}$ is non-empty for all $n\ge\mathbb{N}$. Lemma~9.14 now follows by
uniqueness of analytic continuation. $\endproof$

\begin{lemma} There is a constant $C_6\ge 0$ such that for $\lambda\in\CO$, $x\in B$ and $n\in\mathbb{N}$
\[
\| G_n(\lambda,x)-G(\lambda,x)\|\le\frac{C_6}{n}(1+\|\lambda\|)(1+\|(\mathrm{Im\lambda})^{-1}\|^{14}).
\]
\end{lemma}

\proof At first assume that $\lambda\in\widehat{V}_n$. Put
\[
\mu_n=\Lambda_n(\lambda,x)\qquad\textrm{and}\qquad\mu=\begin{pmatrix}
\Lambda^t & x\\
0 & \Lambda
\end{pmatrix}.
\]
According to Lemma~9.14 and \eqref{eq9-14} we then have that
\[
\widetilde{G}(\mu_u)=G_n(\lambda,x)\qquad\textrm{and}\qquad \widetilde{G}(\mu)=G(\lambda,x).
\]
Hence,
\begin{equation}\label{eq9-37}
\left.\begin{array}{rcl}
\| G_n(\lambda,x)-G(\lambda,x)\| & = & \|\widetilde{G}(\mu_n)-\widetilde{G}(\mu)\|\vspace{2mm}\\
& \le & \|(\mu_n\otimes\mathbf{1}_\CA-\hat{s})^{-1}- (\mu\otimes\mathbf{1}_\CA-\hat{s})^{-1}\|\vspace{2mm}\\
& \le & \|(\mu_n\otimes\mathbf{1}_\CA -s)^{-1}\|\|\mu_n-\mu\| \|(\mu\otimes\mathbf{1}_\CA-\hat{s})^{-1}\|
\end{array}\right..
\end{equation}
By Lemma~9.10~(i),
\[
\|(\mu\otimes\mathbf{1}_\CA-\hat{s})^{-1}\|= \|\pi(\lambda,x)^{-1}\|\le\widehat{C}(1+\|(\mathrm{Im\lambda})^{-1}\|^2),
\]
and by \eqref{eq9-32} and \eqref{eq9-36},
\[
\|(\mu_n\otimes\mathbf{1}_\CA-\hat{s})^{-1}\|\le 2\widehat{C}(1+\|(\mathrm{Im\lambda})^{-1}\|^2),
\]
and
\[
\|\mu_n-\mu\|\le\frac{C_3}{n}(1+\|\lambda\|)(1+\|(\mathrm{Im\lambda}^{-1}\|^{10}).
\]
Inserting these estimates in \eqref{eq9-37} we get that
\[
\| G_n(\lambda,x)-G(\lambda,x)\|\le \frac{C_6^{(1)}}{n}(1+\|\lambda\|)(1+\|(\mathrm{Im\lambda})^{-1}\|^{14})
\]
for some constant $C_6^{(1)}\ge 0$. If $\lambda\notin\widehat{V}_n$, then
by \eqref{eq9-38},
\[
1\le\frac{C_4}{n}(1+\|\lambda\|)(1+\|(\mathrm{Im\lambda})^{-1}\|^{10}) (1+\|(\mathrm{Im\lambda})^{-1}\|^2).
\]
Therefore
\[
\left.\begin{array}{rl}
& \| G_n(\lambda,x)-G(\lambda,x)\|\vspace{2mm}\\
\le & \displaystyle\frac{C_4}{n}(1+\|\lambda\|)(1+\|(\mathrm{Im\lambda})^{-1}\|^{10}) (1+\|(\mathrm{Im\lambda)}^{-1}\|^2) (\| G_n(\lambda,x)\|+\| G(\lambda,x)\|).
\end{array}\right.
\]
Taking (9.44) and (9.45) into account we then get the estimate
\[
\|G_n(\lambda,x)-G(\lambda,x)\| \le
\frac{C_6^{(2)}}{n}(1+\|\lambda\|)(1+\|(\mathrm{Im\lambda})^{-1}\|^{14})
\]\noindent
for a constant $C_6^{(2)}\ge 0$. This proves Lemma~9.15 with $C_6=\max\{
C_6^{(1)},C_6^{(2)}\}$. $\endproof$

{\it Proof of Theorem~9.6} Let $\lambda\in\CO$ and $x\in B$. By
\eqref{eq9-38a} and \eqref{eq9-33},
\[
G_n(\lambda,x)=\begin{pmatrix}
G_n(\lambda)^t & K_n(\lambda,x)\\
0 & G_n(\lambda)
\end{pmatrix},
\]
where $K_n(\lambda,x)$ is given by \eqref{eq9-33}. Similarly,
\[
G(\lambda,x)=\begin{pmatrix}
G(\lambda)^t & K(\lambda,x)\\
0 & G(\lambda)
\end{pmatrix},
\]
where
\begin{equation}
K(\lambda,x)=-(\mathrm{id}_k\otimes\tau) [(\Lambda^t\otimes\mathbf{1}_\CA-s^t)^{-1}(x\otimes\mathbf{1}_\CA)(\Lambda\otimes\mathbf{1}_\CA-s)^{-1}]
\end{equation}\noindent
Hence, by Lemma~9.15,
\[
\| K_n(\lambda,x)-K(\lambda,x)\|\le\frac{C_6}{n}(1+\|\lambda\|) (1+\|(\mathrm{Im\lambda})^{-1}\|^{14}).
\]
But $K_n(\lambda,x)$ and $K(\lambda,x)$ are well-defined for all $x\in M_k(\mathbb{C})$. Moreover, since they are linear in $x$, it follows that
\begin{equation}\label{eq9-40a}
\| K_n(\lambda,x)-K(\lambda,x)\| \le \frac{C_6}{n}(1+\|\lambda\|)(1+\|(\mathrm{Im\lambda})^{-1}\|^{14})\| x\|
\end{equation}
for all $\lambda\in\CO$ and all $x\in M_k(\mathbb{C})$. By \eqref{eq9-3}
and \eqref{eq9-39}, 
\[
R_n(\lambda)=\sum_{j=1}^{r+s}\sum_{u,v=1}^{k}\varepsilon_j a_j e_{uv}^{(k)}K_n(\lambda,e_{uv}^{(k)}a_j),
\]
and
\[
R(\lambda)=\sum_{j=1}^{r+s}\sum_{u,v=1}^{k}\varepsilon_j a_j e_{uv}^{(k)}K(\lambda,e_{uv}^{(k)}a_j).
\]
Hence by \eqref{eq9-40a}
\[
\| R_n(\lambda)-R(\lambda)\| \le \frac{\widetilde{C}_4}{n}(1+\|\lambda\|) (1+\|(\mathrm{Im\lambda})^{-1}\|^{14})
\]\noindent
for same constant $\widetilde{C}_4\ge 0$. $\endproof$

\section{The spectrum of $Q_n$ -- the real case.}\label{sec8}

With the same notation as in the previous section, $x_1,\ldots,x_{r+s}$ is a semicircular system in a $C^*$-probability space $(\CA,\tau)$ with $\tau$ faithful, $X_1^{(n)},\ldots,X_{r+s}^{(n)}$ are stocasticly independent random matrices, for which,
\[
X_1^{(n)},\ldots,X_r^{(n)}\in\mathrm{GOE}\Big(n,\frac{1}{n}\Big)\qquad\textrm{and}\qquad
X_{r+1}^{(n)},\ldots,X_{r+s}^{(n)}\in\mathrm{GOE}^*\Big(n,\frac{1}{n}\Big).
\]\noindent
Let $p\in M_m(\mathbb{C})\otimes\mathbb{C}<X_1,\ldots,X_{r+s}>$ and put
\[
q=p(x_1,\ldots,x_{r+s}),\quad Q_n=p(X_1^{(n)},\ldots,X_{r+s}^{(n)}).
\]\noindent
Moreover define $g,g_n : \mathbb{C}\setminus\mathbb{R}\rightarrow\mathbb{C}$ by
\[
g(\lambda)=(\tr_m\otimes\tau)[(\lambda\unit_m\otimes\unit_\CA-q)^{-1}],\qquad
(\lambda\in\C\setminus \R),
\]
\[
g_n(\lambda)=\E\{(\tr_m\otimes\tr_n)[(\lambda\unit_m\otimes\unit_n-Q_n)^{-1}]\},\qquad
(\lambda\in\C\setminus \R).
\]

Let $E= \unit_m \oplus O_{k-m}\in M_k(\C)$. Then for $\lambda\in \C\setminus\R$, we have
\[
g(\lambda) = \frac{k}{m}\,\tr_k(EG(\lambda\unit_m)E),
\]

and
\[
g_n(\lambda) = \frac{k}{m}\,\tr_k(EG_n(\lambda\unit_m)E).
\]

Now, define $l:\C\setminus\R\rightarrow \C$ by
\[
l(\lambda) = \frac{k}{m}\,\tr_k(EL(\lambda\unit_m)E), \qquad (\lambda\in
\C\setminus \R).
\]

Then $l$ is analytic, and applying Theorem~9.7 we find that
there is a constant $C\geq 0$ such that
\begin{equation}\label{GOE8}
 \Big|g(\lambda)-g_n(\lambda)+\frac 1n l(\lambda)\Big|\leq
  \frac{C}{n^2}(1+|\lambda|)^2(1+|\im\lambda|^{-17})
\end{equation}  
when $\mathrm{Im\lambda}>0$. Moreover, by arguing as in the proof of [S, Theorem~4.5] the inequality (\ref{GOE8}) also holds when $\mathrm{Im\lambda}<0$.
\vspace{.2cm}

\begin{lemma}\label{distribution}
There is a distribution $\Delta \in \Distc$ with $\supp(\Lambda)\subseteq \sigma(q)$, such that for any $\phi \in \Ccinf$,
\begin{equation}\label{dist}
\Delta(\phi) = \lim_{y\rightarrow 0^+}\frac{\i}{2\pi}\int_\R \phi(x)[l(x+\i y)-l(x-\i y)]\d x.
\end{equation}
\end{lemma}

\proof At first we prove that $l$ has an analytic continuation to
$\C\setminus\sigma(q)$. We know that for any
$\lambda\in\C\setminus\sigma(q)$, $\lambda\unit_m\otimes\unit_\CA-q$ is
invertible. Thus with $\Lambda = (\lambda\unit_m) \oplus
\unit_{k-m}\in M_k(\C)$ we know from lemma~3.2 that $\Lambda\otimes \unit_\CA
-s$ is invertible. But then $\Lambda^t\otimes\unit_\CA-s^t$ is also
invertible. It follows that $\lambda\rightarrow R(\lambda\unit_m)$
and $\lambda\rightarrow G(\lambda\unit_m)$ have
analytic continuations to $\C\setminus\sigma(q)$. Moreover,
$G(\lambda\unit_m)$ is invertible for all
$\lambda\in\C\setminus\sigma(q)$. Indeed, $\C\setminus\sigma(q)$ is
connected, and we have seen that for all $\lambda$ belonging to some open
non-empty subset of $\C\setminus\sigma(q)$, the identity
\begin{equation}\label{dist1}
\Big(a_0 + \sumrs a_iG(\lambda\unit_m)a_i - \Lambda\Big)G(\lambda\unit_m) +
\unit_k= 0.
\end{equation}  

holds. Then, by uniqueness of analytic continuation, (\ref{dist1}) must
hold for all $\lambda\in\C\setminus\sigma(q)$. In particular,
$G(\lambda\unit_m)$ is invertible for such $\lambda$. We conclude that $l$ is well-defined and analytic in all of  $\C\setminus\sigma(q)$.

The next step is to prove that $l$ satisfies (a) and (b) of [S, Theorem~5.4].
Let $\lambda\in \C\setminus\R$, and put $\Lambda = (\lambda\unit_m) \oplus
\unit_{k-m}\in M_k(\C)$. According to the proof of lemma~3.2,
\[
(\Lambda\otimes \unit_\CA -s)^{-1} = C + B(\lambda),
\]

where $(Q\otimes\unit_\CA)C=C(Q\otimes\unit_\CA)=0$, and
\[
\|B(\lambda)\|  \leq C'_{2,1}\|(\lambda\unit_m\otimes\unit_\CA-q)^{-1}\|\leq C'_{2,1}|\im\lambda|^{-1}.
\]

Moreover, if $|\lambda|> \|q\|$, then
\[
\|B(\lambda)\|\leq \frac{C'_{2,1}}{|\lambda|-\|q\|}.
\]

Now,
\[
EL(\lambda\unit_m)E = (\id_k\otimes\tau)[(E\otimes\unit_\CA)(\Lambda\otimes \unit_\CA -s)^{-1}(R(\lambda\unit_m)G(\lambda\unit_m)^{-1}\otimes\unit_\CA)(\Lambda\otimes \unit_\CA -s)^{-1}(E\otimes\unit_\CA)],
\]

implying that
\begin{equation}\label{normQLQ}
\|EL(\lambda\unit_m)E\|\leq C_{2,1}^2|\im\lambda|^{-2}\|R(\lambda\unit_m)G(\lambda\unit_m)^{-1}\|,
\end{equation}

and if $|\lambda|>\|q\|$, then
\begin{equation}\label{normQLQa}
\|EL(\lambda\unit_m)E\|\leq \Bigg(\frac{C_{2,1}}{|\lambda|-\|q\|}\Bigg)^2\|R(\lambda\unit_m)G(\lambda\unit_m)^{-1}\|.
\end{equation}

We have seen that
\begin{equation}\label{GOE4}
\|R(\lambda\unit_m)\| \leq   k^2\sum_{i=1}^{r+s}\|a_i\|^2\|(\Lambda\otimes\unit_\CA - s)^{-1}\|^2, 
\end{equation}

where
\begin{equation}\label{GOE7}
\|(\Lambda\otimes\unit_\CA - s)^{-1}\| \leq C'_{1,1}+C'_{2,1}|\im\lambda|^{-1},
\end{equation}

and if $|\lambda|>\|q\|$, then
\begin{equation}\label{GOE6}
\|(\Lambda\otimes\unit_\CA - s)^{-1}\| \leq C'_{1,1}+\frac{C'_{2,1}}{|\lambda|-\|q\|}.
\end{equation}

Also, there is a constant $C_1\geq 0$
such that
\begin{eqnarray}
\|G(\lambda\unit_m)^{-1}\|&\leq &
\|\Lambda\|+\|a_0\|+\sum_{i=1}^{r+s}\|a_i\|\|G(\lambda\unit_m)\|\|a_i\| \\
&\leq & \|\Lambda\|+\|a_0\|+\sum_{i=1}^{r+s}\|a_i\|\|(\Lambda\otimes\unit_\CA - s)^{-1}\|\|a_i\|\\
& \leq & C_1(1+|\lambda|)(1+|\im\lambda|^{-1}),\label{normG-1}
\end{eqnarray}

and if $|\lambda|>\|q\|$, then
\begin{eqnarray}
\|G(\lambda\unit_m)^{-1}\|&\leq & C_1(1+|\lambda|)\Bigg(1+\frac{1}{|\lambda|-\|q\|}\Bigg).\label{normG-1a}
\end{eqnarray}

(\ref{normQLQa}), (\ref{GOE4}), (\ref{GOE6}) and (\ref{normG-1a}) imply
that
\begin{equation}\label{conditiona}
  |l(\lambda)|\leq O\Big(\frac{1}{|\lambda|}\Big) \quad {\rm as} \quad
   |\lambda|\rightarrow \infty.
\end{equation}

Combining (\ref{normQLQ}), (\ref{GOE4}), (\ref{GOE7}) and (\ref{normG-1})
we find that for some constant $C_2\geq 0$,
\begin{equation}\label{norml}
|l(\lambda)| \leq \|EL(\lambda\unit_m)E\|\leq C_2(1+|\lambda|)( |\im\lambda|^{-2} + |\im\lambda|^{-5}).
\end{equation}

Choose $a, b\in \R$, $a<b$, such that $\sigma(q)\subseteq [a,b]$. Put
$K= [a-1,b+1]$ and
\[
D= \{\lambda\in\C|0<{\rm dist}(\lambda,K)\leq 1\}.
\]

By (\ref{norml}), there is a constant $C_3\geq 0$ such that for any
$\lambda\in D$,
\[
|l(\lambda)| \leq C_3\cdot \max\{1,({\rm dist}(\lambda,K))^{-5} \} =
 C_3\cdot({\rm dist}(\lambda,K))^{-5}.
\]

(\ref{conditiona}) implies that $l$ is bounded on $\C\setminus
D$. Therefore $C_3$ may be chosen such that for all $\lambda\in \C\setminus \sigma(q)$
\begin{equation}\label{conditionb}
  |l(\lambda)| \leq C_3\cdot \max\{1,({\rm dist}(\lambda,K))^{-5} \}.
\end{equation}

By (\ref{conditiona}) and (\ref{conditionb}), $l$ satisfies (a) and (b) of
[S, Theorem~5.4], and the lemma follows. $\endproof$

\vspace{.2cm}

Knowing that (\ref{GOE8}) holds we are now able to prove:

\begin{thm}\label{GOE9}
Let $\phi\in \Ccinf$. Then
\[
\E\{(\tr_m\otimes\tr_n)\phi(Q_n)\} = (\tr_m\otimes\tau)\phi(q) + \frac 1n
\Delta(\phi)+O\Big(\frac{1}{n^2}\Big).
\]
\end{thm}

\proof
The result follows from a simple modification of the proof of [S, Theorem~5.6].
$\endproof$

\vspace{.2cm}

\begin{lemma}\label{rmeanofnorms}
Let $n\in \N$, and let $X_n\in \GOE(n,\frac 1n)\cup  \GOES(n,\frac
1n)$. Then
\begin{itemize}
  \item[(i)] for all $\eps >0$, $P(\|X_n\|>\sqrt 2
  (2+\eps))<2n\,\exp\Big(-\frac{n\eps^2}{2}\Big)$,
  \item[(ii)] there is a sequence of constants not depending on $n$, $(\gamma'(k))_{k=1}^\infty$, such that for all $k\in\N$,
  \[
  \E\{1_{(\|X_n\|>3\sqrt 2)}\|X_n\|^k\}\leq \gamma'(k)n\e^{-\frac n2}.
  \]
\end{itemize}
\end{lemma}

\vspace{.2cm}

\proof We may assume that $X_n = \frac{1}{\sqrt 2}(Y_n+\overline Y)$ or
$X_n = \frac{1}{\sqrt 2\i}(Y_n-\overline Y)$ for some $Y_n\in\SGRM(n,\frac
1n)$. Hence, by [S, Proof of Lemma~6.4],
\[
P(\|X_n\|>\sqrt 2(2+\eps))\leq P(\|Y_n\|> 2+\eps)\leq 2n\,
\exp\Big(-\frac{n\eps^2}{2}\Big)
\]

holds for all $\eps >0$. Then by application of the proof of
Proposition~\ref{meanofnorms} to the random variable $\frac{1}{\sqrt 2}\|X_n\|$ we get
that 
 \[
  \E\Big\{1_{\big(\frac{1}{\sqrt 2}\|X_n\|>3\big)}\Big(\frac{1}{\sqrt 2}\|X_n\|\Big)^k\Big\}\leq \gamma(k)n\e^{-\frac n2}.
  \]

Hence (ii) holds with $\gamma'(k)=2^{k/2}\gamma(k)$. $\endproof$

\vspace{.2cm}

\begin{lemma}\label{1/n-term} Let $\Delta\in \Distc$ be as in
  Lemma~\ref{distribution}. Then $\Delta(1)=0$.
\end{lemma}

\proof With $d={\rm deg}(p)$, $x_0 = \unit_\CA$, and $X_0^{(n)}=\unit_n$ we
may choose $c_{i_1,\ldots, i_d}\in M_m(\C)$, $0\leq i_1,\ldots, i_d\leq
r+s$, such that
\[
q= p(x_1, \ldots, x_{r+s})=\sum_{0\leq i_1,\ldots, i_d\leq r+s}c_{i_1,\ldots,
  i_d}\otimes x_{i_1}\cdots x_{i_d}
\]

and
\[
Q_n = p(X_1^{(n)}, \ldots, X_{r+s}^{(n)})=\sum_{0\leq i_1,\ldots, i_d\leq r+s}c_{i_1,\ldots,
  i_d}\otimes X_{i_1}^{(n)}\cdots X_{i_d}^{(n)}.
\]

Put
\[
R = (3\sqrt2)^d \sum_{0\leq i_1,\ldots, i_d\leq r+s}\|c_{i_1,\ldots,
  i_d}\|.
\]

Then
\[
(\|Q_n\| >R) \subseteq (\sum_{0\leq i_1,\ldots, i_d\leq r+s}\|c_{i_1,\ldots,
  i_d}\|\|X_{i_1}^{(n)}\|\cdots \|X_{i_d}^{(n)}\| >R) \subseteq
  \bigcup_{i=1}^{r+s}(\|X_i^{(n)}\| >3\sqrt2 ),
\]

implying that
\[
P(\|Q_n\| >R) \leq r\cdot P(\|X_1^{(n)}\|
>3\sqrt2) + s \cdot P(\|X_{r+1}^{(n)}\| >3\sqrt2).
\]

Now, by Lemma~\ref{rmeanofnorms},
\[
 P(\|X_i^{(n)}\|>3\sqrt2) \leq 2n\cdot \exp\Big(-\frac{n}{2}\Big), \qquad (i=1,
 \ldots, r+s),
\]

and thus
\[
P(\|Q_n\| >R) \leq 2(r+s)n\cdot \exp\Big(-\frac n2\Big).
\]

Consequently,
\begin{equation}\label{indicator}
\E\{(\tr_m\otimes \tr_n)1_{]-\infty,R[\cup]R,\infty[}(Q_n)\} \leq
P(\|Q_n\|>R)\leq 2(r+s)n\cdot \exp\Big(-\frac n2\Big).
\end{equation}

Now, let $\phi\in \Ccinf$ such that $0\leq \phi \leq 1$ and
$\phi|_{[-R,R]}=1$. Then $\phi(x) = 1$ for all $x$ in a neighbourhood of
$\sigma(q)\supseteq \supp(\Lambda)$. Hence
$\Delta(\phi) = \Delta(1)$, and we have that
\[
\E\{(\tr_m\otimes\tr_n)\phi(Q_n)\} = 1 +
\frac 1n \Delta(1) + O\Big(\frac{1}{n^2}\Big),
\]

where
\[
\E\{(\tr_m\otimes\tr_n)\phi(Q_n)\} =
1 + \E\{(\tr_m\otimes\tr_n)(\phi-1)(Q_n)\}
\]

and
\[
|\E\{(\tr_m\otimes\tr_n)(\phi-1)(Q_n)\}|
 \leq \E\{(\tr_m\otimes \tr_n)1_{]-\infty,R[\cup]R,\infty[}(Q_n)\} \leq 2(r+s)n\cdot \exp\Big(-\frac n2\Big).
\]

Altogether we have that
\begin{equation}\label{phi-1}
\E\{(\tr_m\otimes\tr_n)(\phi-1)(Q_n)\} =
\frac 1n \Delta(1) + O\Big(\frac{1}{n^2}\Big),
\end{equation}

where the left hand side is of the order $n\cdot \exp(-\frac n2)= O\Big(\frac{1}{n^2}\Big)$. Hence, the $\frac 1n$-term
appearing on the right hand side of (\ref{phi-1}) must be zero. $\endproof$

\vspace{0.2cm}

\begin{prop}\label{HST6.2}
Let $\phi\in C^\infty(\R,\R)$ such that $\phi$ is constant outside a
compact subset of $\R$. Suppose that
\[
\supp(\phi)\cap\sigma(q)=\emptyset.
\]

Then
\[
\V\{(\tr_m\otimes\tr_n)\phi(Q_n)\}= O(n^{-4}),
\]

and
\[
P(|(\tr_m\otimes\tr_n)\phi(Q_n)|\leq n^{-\frac 43}, \; {\rm eventually \;
  as }\; n\rightarrow \infty)=1.
\]
\end{prop}

\vspace{.2cm}

\proof Taking Lemma~\ref{rmeanofnorms} (ii) into account this result
follows as in the complex case (cf. proof of Theorem~\ref{spec1}). $\endproof$

Taking Theorem~\ref{GOE9},
Lemma~\ref{1/n-term}, Proposition~\ref{HST6.2} and Proposition~\ref{variance estimate} into account,
we find, as in Section~\ref{sec4b}:

\begin{thm}\label{GOE10} 
Let $\eps>0$. Then for almost every $\omega\in\Omega$,
\[
\sigma(Q_n(\omega))\subseteq \sigma(q)+\,]-\eps,\eps[,
\]
eventually as $n\rightarrow\infty$.
\end{thm}

\begin{remark}
By [S, Section~7] and the remarks in the beginning of section~9, Theorem~10.2, Proposition~10.5 and Theorem~10.6 can easily be generalized to the symplectic case i.e. to the case where $X_1^{(n)},\ldots,X_{r+s}^{(n)}$ are stochasticly independent random matrices for which
\[
X_1^{(n)},\ldots,X_r^{(n)}\in\mathrm{GSE}\Big(n,\frac{1}{n}\Big)\qquad\textrm{and}\qquad
X_{r+1}^{(n)},\ldots,X_{r+s}^{(n)}\in\mathrm{GSE}^*\Big(n,\frac{1}{n}\Big).
\]\noindent
\end{remark}

\section{Gaps in the spectrum of $q$ --- the real case}
\label{sec9}
\setcounter{equation}{0}

In this section we shall prove that Theorem~\ref{holes1} holds in the $\GOE
\cup\GOE^*$-case as well. That is, if $p\in (M_m(\bC)\otimes\bC\langle
X_1,\dots,X_{r+s}\rangle)_\sa$, if $x_1,\dots,x_{r+s}$ is a semicircular
system in $(\cA,\tau)$, and if for each $n\in\bN$,
$X_1^{(n)},\dots,X_{r+s}^{(n)}$ are stochastically independent random
matrices from $\GOE(n,\frac1n)\cup\GOE^*(n,\frac1n)$ as in Section~\ref{sec7}, then with $q=p(x_1,\dots,x_{r+s})$ and $Q_n=
p(X_1^{(n)},\dots,X_{r+s}^{(n)})$ we have:

\begin{thm}
\label{thm8-1}
Let $\eps_0$ denote the smallest distance between disjoint connected
components of $\sigma(q)$, let $\cJ$ be a connected component of
$\sigma(q)$, let $0<\eps<\frac13\eps_0$, and let $\mu_q\in \Prob(\bR)$ denote
the distribution of $q$ w.r.t. $\tr_m\otimes\tau$. Then $\mu_q(\cJ)=\frac
km$ for some $k\in \{1,\dots,m\}$, and for almost all $w\in\Omega$, the
number of eigenvalues of $Q_n(w)$ in $\cJ+]-\eps,\eps[$ is $k\cdot n$,
eventually as $n\to\infty$.
\end{thm}

\proof
Take $\phi\in C_c^\infty(\bR)$, such that $\phi_{|\cJ+]-\eps,\eps[}=1$ and
$\phi_{|\bR\backslash(\cJ+]-2\eps,2\eps[)}=0$. Then $\phi(q)$ is a non-zero
projection in $M_m(C^*(\unit_\CA,x_1,\dots,x_n))$ and hence, by Theorem~\ref{noprojections},
\begin{equation}
\label{eq8-1}
\mu_q(\cJ) = \int_\bR\varphi\, d\mu_q= (\tr_m\otimes\tau)\phi(q)=\frac km
\end{equation}
for some $k\in\{1,\dots,m\}$. By Theorem~\ref{GOE10} there is a $P$-null
set $N\subset\Omega$, such that for all $\omega\in\Omega\backslash N$
\[
\sigma(Q_n(\omega)) \subseteq\sigma(q)+]-\eps,\eps[,
\]
eventually as $n\to\infty$. In particular, for all $\omega\in\Omega\backslash N$
there exists $N(\omega)\in\bN$ such that $\phi(\Q_n(\omega))$ is a projection for
all $n\ge N(w)$. For $\omega\in\Omega\backslash N$ and $n\ge N(\omega)$ take
$K_n(\omega)\in \{0,\dots,m\cdot n\}$, such that
\begin{equation}
\label{eq8-2}
(\tr_m\otimes\tr_n)\varphi(Q_n(\omega)) = \frac{K_n(\omega)}{m\cdot n}.
\end{equation}
Let $\Delta\colon C_c^\infty(\bR)\to\bC$ be the distribution from Lemma~\ref{distribution}, and put
\[
Z_n =
(\tr_m\otimes\tr_n)\phi(Q_n)-(\tr_m\otimes\tau)\phi(q)-\frac1n\Delta(\phi).
\]
Then by Theorem~\ref{GOE9}, $\bE(Z_n)=O(\frac{1}{n^2})$. Moreover, since
$\varphi'$ vanishes in a neighbourhood of $\sigma(q)$, we get as in the
proof of Theorem~\ref{spec1} that
\[
\bV(Z_n) = O(\tfrac{1}{n^4})
\]
and
\[
Z_n = O(n^{-\frac43})\quad\mbox{almost surely.}
\]
Hence there exists a $P$-null set $N'\subseteq N$, such that
\[
(\tr_m\otimes\tr_n)\phi(Q_n(\omega)) = (\tr_m\otimes\tau)\phi(q) + \frac1n
\Delta(\phi) +O(n^{-\frac43})
\]
holds for all $\omega\in\Omega\setminus N'$.

Taking \eqref{eq8-1} and \eqref{eq8-2} into account, we get after
multiplication by $mn$ that for $\omega\in \Omega\setminus N$ and $n\geq N(\omega)$,
\begin{equation}
\label{eq8-3}
K_n(\omega) = nk + m\Delta(\varphi)+O(n^{-\frac13}).
\end{equation}
Therefore there exists $C>0$ such that $\dist(m\Delta(\varphi),\bZ)\le C\,
n^{-\frac13}$ for all $n\ge N(\omega)$, which implies that
$m\Delta(\varphi)\in\bZ$.

We next use an argument based on homotopy to show that $\Delta(\varphi)=0$. By
definition
\[
Q_n=p(X_1^{n)},\dots,X_{r+s}^{(n)})
\]
where $X_1^{(n)},\dots,X_r^{(n)}\in\GOE(n,\frac1n)$,
$X_{r+1}^{(n)},\dots,X_{r+s}^{(n)} \in\GOE^*(n,\frac1n)$ form a set of
$r+s$ independent random matrices. We may without loss of generality assume
that there exist $Y_1^{(n)},\dots,Y_r^{(n)}\in\GOE^*(n,\frac1n)$ and
$Y^{(n)}_{r+1},\dots,Y^{(n)}_{r+s}\in\GOE(n,\frac1n)$ such that
$X_1^{(n)},\dots,X^{(n)}_{r+s}$, $Y_1^{(n)},\dots,Y^{(n)}_{r+s}$ form a set of
$2(r+s)$ independent random matrices. For $j=1,\dots,r+s$ put
\[
X_j^{(n)}(t)=\cos t \; X_j^{(n)}+\sin t \; Y_j^{(n)},\quad (0\le t\le
\tfrac{\pi}{4}).
\]
It is a simple observation that if $Z\in \GOE(n,\frac1n)$ and
$W\in\GOE^*(n,\frac1n)$, then
$\frac{1}{\sqrt{2}}(Z+W)\in\SGRM(n,\frac1n)$. Hence
\[
(X_1^{(n)}(t),\dots,X^{(n)}_{r+s}(t)),\quad (0\le t\le\tfrac{\pi}{4})
\]
defines a path which connects the given set of random matrices
$X_1^{(n)},\dots,X^{(n)}_{r+s}$ (at $t=0$) with a set of $r+s$ independent
$\SGRM(n,\frac1n)$ random matrices (at $t=\frac{\pi}{4}$). Put
\[
Q_n(t) = q(X_1^{(n)}(t),\dots,X_{r+s}^{(n)}(t)),\quad (0\le
t\le\tfrac{\pi}{4}).
\]
Let $x_1,\dots,x_{r+s},y_1,\dots,y_{r+s}$ be a semicircular system in a
$C^*$-probability space $(\cA,\tau)$ with $\tau$ faithful. Put
\[
x_j(t)=\cos t\; x_j+\sin t\; y_j,\quad (0\le t\le\tfrac{\pi}{4}).
\]
Since an orthogonal transformation of a semicircular system is again a
semicircular system (cf. \cite[Proposition~5.12]{VDN}), $x_1(t),\dots,x_{r+s}(t)$
is a semicircular system for each $t\in [0,\frac{\pi}{4}]$. Hence the
operators
\[
q(t)=q(z_1(t),\dots,x_{r+s}(t)),\quad (0\le t\le\tfrac{\pi}{4})
\]
form a norm continuous path in $\cA$ for which
$\sigma(q(t))=\sigma(q)$. Moreover
\begin{equation}
\label{eq8-4}
(\lambda,t)\to(\lambda \unit_m\otimes \unit_\cA-q(t))^{-1}
\end{equation}
is norm continuous on $(\bC\backslash \sigma(q))\times [0,\frac{\pi}{4}]$. For $t\in
[0,\frac{\pi}{4}]$, $Q_n(t)$ can be expressed as a polynomial in
$X_1^{(n)},\dots,X_{r+s}^{(n)},Y_1^{(n)},\dots,Y_{r+s}^{(n)}$, and $q(t)$ can be
expressed as the same polynomial in
$x_1,\dots,x_{r+s},y_1,\dots,y_{r+s}$. Hence, by Lemma~\ref{distribution} and Theorem~\ref{GOE9},
there exists for each $t\in[0,\frac{\pi}{4}]$ a distribution
$\Lambda_t\colon C^\infty_r(\bC)\to\bC$, such that for all $\psi\in
C_c^\infty(\bR)$:
\[
\bE\{(\tr_m\otimes\tr_n)\psi(Q_n(t))\}=(\tr_m\otimes\tau)\psi(q(t))+\frac1n\Lambda_t(\psi)+
O(\tfrac{1}{n^2}). 
\]
Since $\sigma(q(t))=\sigma(q)$, $0\le t\le\frac{\pi}{4}$, we get by the
first part of this proof that
\begin{equation}
\label{eq8-5}
m\Delta_t(\varphi)\in\bZ,\quad (0\le t\le\tfrac{\pi}{4}),
\end{equation}
where  $\varphi$ is the function chosen in the beginning of the
proof. Moreover, by Theorem~\ref{spec2},
\[
\bE\{(\tr_m\otimes\tr_n)\varphi(Q_n(\tfrac{\pi}{4})\} =
(\tr_m\otimes\tau)\varphi(q(\tfrac{\pi}{4}))+O(\tfrac{1}{n^2}),
\]
which implies that
\begin{equation}
\label{eq8-6}
\Delta_{\pi/4}(\varphi)=0.
\end{equation}
We next prove that $t\to \Delta_t(\varphi)$ is a continuous function:

Let
\[
\ell_t(\lambda) = \Delta_t\Big(\frac{1}{\lambda -x}\Big),\quad
(\lambda\in\bC\backslash \sigma(q))
\]
be the Stieltjes transformation of $\Delta_t$, $0\le t\le\frac{\pi}{4}$
(cf. \cite[Lemma~5.4]{S}). By a simple modification of the proof of
\cite[Lemma 5.6]{S}, we get that
\begin{equation}
\label{eq8-7}
\Delta_t(\varphi) = \frac{1}{2\pi i}\int_{\partial R}
\ell_t(\lambda)d\lambda
\end{equation}
where $\partial R$ is the boundary of the rectangle
\[
R = (\cJ + [-\varepsilon,\varepsilon])\times[-1,1]
\]
with counter clockwise orientation. Since $(\lambda,t)\to(\lambda
\unit_m\otimes\unit_\CA-q(t))^{-1}$ is norm continuous on $(\bC\backslash \sigma(q))\times
[0,\frac{\pi}{4}]$, it is obvious from the explicit formula for the Stieltjes transform $\ell_t(\lambda)$
(cf. (\ref{dist})) that
$(\lambda,t)\to\ell_t(\lambda)$ is a continuous function on $(\bC\backslash
\sigma(q))\times [0,\frac{\pi}{4}]$. Hence by \eqref{eq8-7},
$\Delta_t(\varphi)$ is a continous funtion of
$t\in[0,\frac{\pi}{4}]$. Together with \eqref{eq8-5} and \eqref{eq8-6} this
shows that $\Delta(\varphi)=\Delta_0(\varphi)=0$. Hence by \eqref{eq8-3} we
have for all $\omega\in \Omega\setminus N'$ and all $n\ge N(\omega)$ that
\[
K_n(\omega)=nk+O(n^{-\frac13}),
\]
and since $K_n(\omega)\in\bN$, it follows that $k_n(\omega)=nk$ eventually as
$n\to\infty$.  
\eproof

\begin{remark}
Using again [S, Section~7], Theorem~11.1 can also be generalized to the symplectic case (cf. remark~10.7).
\end{remark}

Department of Mathematics and Computer Science\\
University of Southern Denmark\\
Campusvej 55, 5230 Odense M\\
Denmark\\
\texttt{
haagerup@imada.sdu.dk\\
schultz@imada.sdu.dk\\
steenth@imada.sdu.dk
}
\end{document}